\documentclass[12pt]{amsart}
\usepackage{amsmath,amscd,latexsym,verbatim,amssymb}
\usepackage{times}       

 \newcommand{\resp}{{\it resp.} }
\newcommand{\cf}{{\it cf.} }
\newcommand{\ie}{{\it i.e.} }
\newcommand{\eg}{{\it e.g.} }

\newcommand{\loccit}{{\it loc. cit.} }

\newcommand{\F}{\mathbf{F}}
\newcommand{\un}{\mathbf{1}}
\newcommand{\la}{\lambda}
\newcommand{\La}{\Lambda}
\newcommand{\N}{\mathbb{N}}
\newcommand{\Q}{\mathbb{Q}}
\newcommand{\R}{\mathbb{R}}
\newcommand{\C}{\mathbb{C}}
  \newcommand{\Z}{\mathbb{Z}}
\newcommand{\bA}{\mathbb{A}}
\newcommand{\bG}{\mathbb{G}}

\newcommand{\bP}{\mathbb{P}}

  \newcommand{\sA}{{\mathcal{A}}}

\newcommand{\sC}{{\mathcal{C}}}
\newcommand{\sD}{{\mathcal{D}}}

\newcommand{\sG}{{\mathcal{G}}}

\newcommand{\sF}{{\mathcal{F}}}

\newcommand{\sN}{{\mathcal{N}}}
\newcommand{\sO}{{\mathcal{O}}}
\newcommand{\sP}{{\mathcal{P}}}

\newcommand{\sR}{{\mathcal{R}}}

\newcommand{\sT}{{\mathcal{T}}}

\newcommand{\sY}{{\mathcal{Y}}}

\newcommand{\inj}{\hookrightarrow}

\newcommand{\surj}{\rightarrow\!\!\!\!\!\rightarrow}

   \newcommand{\ord}{{\rm{ord}}}

  \newcommand{\coim}{\operatorname{coim}}
\newcommand{\coker}{\operatorname{coker}}
\newcommand{\Ker}{\operatorname{Ker}}
\newcommand{\IM}{\operatorname{Im}}
\newcommand{\Coim}{\operatorname{Coim}}
\newcommand{\Coker}{\operatorname{Coker}}
\newcommand{\rk}{\operatorname{rk}}
\newcommand{\sk}{\operatorname{sk}}
 \newcommand{\gr}{\operatorname{gr}}
\newcommand{\im}{\operatorname{im}}
 \newcommand{\Spec}{\operatorname{Spec}}
\newcounter{spec}

\swapnumbers

\newtheorem{thm}{Theorem}[subsection]
\newtheorem{lemma}[thm]{Lemma}
\newtheorem{prop}[thm]{Proposition}

\newtheorem{cor}[thm]{Corollary}

\theoremstyle{definition}
 \newtheorem{defn}[thm]{Definition}
 
\newtheorem{ex}[thm]{Example}
\newtheorem{exs}[thm]{Examples}

\newtheorem{rem}[thm]{Remark}
\newtheorem{rems}[thm]{Remarks}

\numberwithin{equation}{section}

 at10pt
 at12pt

\renewcommand{\qed}{\hfill $\square$\medskip}

\begin{document}

\title[Slope filtrations]{Slope filtrations}

\author{Yves
Andr\'e}
 
   \address{D\'epartement de Math\'ematiques et Applications, \'Ecole Normale Sup\'erieure  \\ 
45 rue d'Ulm,  75230
  Paris Cedex 05\\France.}
\email{andre@dma.ens.fr}
   \date{\today}

 \maketitle

 \begin{sloppypar}

  \begin{abstract}    Many  slope filtrations occur in algebraic geometry, asymptotic analysis, ramification theory, p-adic theories, geometry of numbers... These functorial filtrations, which are indexed by rational (or sometimes real) numbers, have a lot of common properties. 

 We propose a unified abstract treatment of slope filtrations, and survey how new ties between different domains have been woven by dint of deep correspondences between different concrete slope filtrations. 
  \end{abstract}  
 
 \bigskip
 
{{\tableofcontents}}

 \bigskip
 
\section*{Introduction}  

\medskip
\subsection{}   Slope filtrations occur in algebraic and analytic geometry, in asymptotic analysis, in ramification theory, in p-adic theories, in geometry of numbers... Five basic examples are the Harder-Narasimhan filtration of vector bundles over a smooth projective curve, the Dieudonn\'e-Manin filtration of  $\rm F$-isocrystals over the $p$-adic point, the Turrittin-Levelt filtration of formal differential modules,  the Hasse-Arf filtration of finite Galois representations of local fields, and the Grayson-Stuhler filtration of euclidean lattices. 

Despite the variety of their origins, these filtrations share a lot of similar features.

\smallskip In this paper of bourbachic inspiration, we develop a unified and systematic abstract treatment of slope filtrations, with the aim of freeing the ``yoga of stability" from any ad hoc property of the underlying category. This should not only clarify the analogies, but also allow to replace the pervasive adaptations of arguments from one context to another by a single formal argument.

Such an argument may supplant some quite nonformal arguments in the literature. For instance, it is sometimes considered that proving that the slopes of subobjects are bounded from above is a required preliminary step in the construction of a slope filtration.  An a priori proof of boundness may be difficult in specific instances (\cf \eg \cite{Bru}), but the general theory shows that it is unnecessary: boundness rather appears as a corollary.  

\medskip \subsection{}\label{0.2}  Loosely speaking, (descending) slope filtrations are filtrations of  objects $M$ of a given additive category $\sC$ by 
 subobjects $F^{\geq \la} \, M $ indexed by real numbers. 
 The filtration $F^{\geq .}M$ is supposed to be functorial in $M$, and to be left-continuous and locally constant in $\la$: it comes from a finite flag  
 
$$ 0 \subset  F^{\geq \la_1} \, M  \subset \ldots \subset    F^{\geq \la_r} \, M = M  $$ 
where the $\la_1> \ldots > \la_r$ are the {\it breaks} of the filtration.
     
\smallskip On the other hand, it is assumed that objects of $\sC$ have a well-defined {\it rank} in $\mathbb N$ (typically they are linear objects with some extra structure, and the rank refers to the underlying linear structure). 
    This allows to attach to any object $M$ its {\it Newton polygon}: the polygon which lies below the  concave piecewise linear curve $Np(M)$ emanating from the origin, whose breaks (including end-points) are at the abscissa $0, \ldots , \rk F^{\geq \la_i} \, M, \ldots, \rk M$, and which has slope $\la_i$ between the abscissa $\rk F^{\geq \la_{i-1}}$ and $\rk F^{\geq \la_i}$.

The ``principle" is that, in the presence of slope filtrations, {\it one can ``unscrew" objects $M$ according to their Newton polygons, functorially in $M$}. In almost all ``natural examples", this principle is enhanced by the combinatorial constraints coming from the fact that the coordinates of the vertices of Newton polygons are integers.

When the underlying category is tannakian, this is a powerful tool to compute tannakian groups (see \eg how N. Katz \cite{Ka2} uses the Turrittin-Levelt slope filtration to compute differential Galois groups).

\medskip \subsection{}\label{0.3} The {\it degree} $\deg M$ is the ordinate of the right end-point of $Np(M)$ (with abscissa $\rk M$).

 The  degree function $\deg$ on $Ob \,\sC$ which is attached in this way to $F^{\geq.}$ satisfies some simple axioms (\cf \ref{sf} below). We show that, conversely, any function on $Ob \,\sC$ satisfying these axioms is the degree function attached to a unique slope filtration on $\sC$ ({\it theorem} \ref{T1}). 
  
 \smallskip This general fact synthesizes (and supersedes) the numerous constructions of concrete slope filtrations of Harder-Narasimhan type found in the literature.   
   
\medskip\subsection{} In most examples, the category $\sC$ is additive, but quite often non-abelian. We show that the right context is that of {\it quasi-abelian categories}: additive categories with kernels and cokernels in which $Ext( -, -)$ is bifunctorial (this notion goes back to Yoneda \cite{Y}).  However, the categories of hermitian coherent sheaves which occur in the context of Arakelov geometry are not additive, and we have to introduce a non-additive version of quasi-abelian categories (which we call {\it proto-abelian categories}) in order to deal with these examples on an equal footing.

 We also analyse in detail the exactness properties of slope functions ({\it theorem} \ref{T2}), and we indicate how slope filtrations are related to stability structures on triangulated categories (in the sense of Bridgeland). 

\medskip\subsection{} Usually, the underlying quasi-abelian category $\sC$ is endowed with a natural tensor product $\otimes$. This leads us to develop the notion of a {\it quasi-tannakian category}. 

One can distinguish two radically different behaviours of slope filtrations wih respect to $\otimes$. 

In the first type of slope filtrations, the breaks of $M\otimes N$ are the sums of a break of $M$ and a break of $N$ ($\otimes$-multiplicative filtrations: \eg Harder-Narasimhan, Dieudonn\'e-Manin).

 In the second type, the breaks of $M\otimes N$ are bounded by the maximum of breaks of $M$ and $N$ ($\otimes$-bounded filtrations: \eg Turrittin-Levelt, Hasse-Arf). 
 We analyse these two types of slope filtrations in general ({\it theorems} \ref{DETFI}, \ref{T3}, \ref{T4}).
 
\medskip \subsection{} The paper begins with a review of the five slopes filtrations mentioned above, and its last portion consists in a reasoned catalogue of slope filtrations in a variety of mathematical domains, underlining a number of links between them. 
 
It ends with a review of some semicontinuity results for Newton polygons in families (with respect to the Harder-Narasimhan, Dieudonn\'e-Manin, Turrittin-Levelt filtrations respectively).

\medskip  We hope that this unified setting will inspire some further transfers of ideas from one domain to another.

 \newpage

\addtocontents{toc}{{\bf I. General theory of slope filtrations.}\hfill\thepage}
\
\bigskip
\begin{center}
\large\bf I. General theory of slope filtrations.
\end{center}
\bigskip

   \section{Brief review of five basic examples.}
   
    \medskip
  \subsection{Harder-Narasimhan filtration of vector bundles.}\label{hnf} 
 
 Let $X$ be the smooth connected projective curve over $\C$. The classification of vector bundles of given rank and degree on $X$ is not straightforward: in order to construct nice moduli schemes, one should either rigidify them, or consider only those vector bundles that Mumford called (semi)stable.
 
 \smallskip Let $N$ be non-zero vector bundle. Its degree $\deg N$ is the degree of its determinant line bundle.
 Its slope is the ratio $\mu(N)=\frac{\deg N}{\rk N}$.
 
  $N$ is {\it stable} (\resp {\it semistable}) if and only if for any non-zero subbundle $M$, $\mu(M)< \mu(N)$ (\resp $\mu(M)\leq \mu(N)$).
  
  Any semistable bundle $N$ of slope $\la $ is a successive extension of stable bundles of slope $\la$. 
  
  Any bundle $N$ is a successive extension of semistable bundles of increasing slopes: more precisely 
  $N$ has a unique descending filtration - the Harder-Narasimhan filtration \cite{HN} -
   $$ 0 \subset   F^{\geq \la_1} \, N  \subset \ldots \subset       F^{\geq \la_r} \, N  = N$$ 
for which $\la_1> \ldots > \la_r$, and the graded pieces $ gr^{\la_i}\,N = F^{\geq \la_i}\,N/F^{> \la_i}\,N$ are semistable bundles of slope $\la_i$. Moreover, $\deg N$ coincides with the degree attached to this filtration in the sense of \ref{0.3}.
 
\smallskip Narasimhan and Seshadri \cite{NS} have described stable bundles in terms of monodromy representations; stable bundles of degree $0$ correspond irreducible unitary representations of $\pi_1(X(\C))$.

 \bigskip
   \subsection{Dieudonn\'e-Manin filtration of $\phi$-modules.}\label{dmf} 
  
In Dieudonn\'e theory of formal groups and crystalline cohomology, one encounters finite dimensional vectors spaces over a $p$-adic field, endowed with an injective semilinear endomorphism. The classification of these objects is due to Dieudonn\'e and Manin  \cite{Di}\cite{Man}. 

\smallskip Let $K$ be a complete valued field of characteristic $0$, with residue field $k$ of characteristic $p>0$. Let $\phi$ be a lifting of some fixed positive power of the Frobenius endomorphism of $k$. In particular, $\phi$ is an isometric endomorphism of $K$.

Let $N$ be a $\phi$-module\footnote{also called $\rm F$-isocrystals (over the point), after Grothendieck.} over $K$, \ie a finite dimensional $K$-vector space $N$ endowed with an isomorphism $\, \Phi_N:\,N\otimes_{K, \phi}K\overset{\cong}{\to} N$ . The determinant $\det N$ is a rank one $\phi$-module; in a given basis, $\Phi_{\det N}$ is given by an element $a\in K^\times$, well-defined up to multiplication by an element of the form $b/\phi(b), \,b\in K^\times$.  Thus the valuation of $\Phi_{\det N}$ is well-defined (\ie as the valuation of $a$).

Let us define $\mu(N)$ to be $-\frac{v(\Phi_{\det N})}{\rk N}\,$ \footnote{unlike the usual convention, we have put a sign in order to get a descending filtration, which fits into the general convention of this paper to deal with descending filtrations. See \ref{asc} and \ref{asc2} for the easy dictionary between descending and ascending slope filtrations.}, 
and say that  $N$ is {\it isoclinic}  if and only if for any non-zero $\phi$-submodule $M$, $\mu(M)= \mu(N)$.
     
\smallskip There is a unique descending filtration - the Dieudonn\'e-Manin filtration\footnote{or, rather, the descending version of the original Dieudonn\'e-Manin filtration.}-
  $$ 0 \subset   F^{\geq \la_1} \, N  \subset \ldots \subset       F^{\geq \la_r} \, N  = N$$ 
for which $\la_1> \ldots > \la_r$, and  $gr^{\la_i}\,N$ is {isoclinic of slope $\la_i$}.

\smallskip Moreover, if $k$ is perfect, $\phi$ is bijective and the Dieudonn\'e-Manin filtration splits\footnote{although the category of $\phi$-modules need not be semisimple.}.
  If $k$ is algebraically closed, simple $\phi$-modules $N$ can be described explicitly: $\varpi^{\Phi_{\det N}}$ is prime to $\rk N$, and $N$ admits a cyclic basis (with respect to $\Phi$) such that the image of the last vector is the first vector multiplied by $\varpi^{\Phi_{\det N}}$ (where $\varpi$ denotes an uniformizer of $K$).

\medskip
  \subsection{Turrittin-Levelt filtration of formal differential modules.}\label{tlf} 
  
 In the field of analytic linear differential equations, the classical opposition singular versus irregular singularities goes back to Fuchs. 
 
The derivation $\partial = x\frac{d}{dx}$ acts on ${{K}} = \C((x))$, respecting the valuation $\ord_x$.   
 A linear differential  operator $\displaystyle P = \partial^n - a_{n-1}\partial^{n-1}-\ldots -a_0$  is {\it regular}
if the ``Fuchs number" 
\begin{equation}\label{fn} {\rm ir}\,P ={ \max}(0, {\max} (-{\ord}_x a_i)) \end{equation}
 is zero. Actually, this number depends only on the associated differential module\footnote{a differential module over ${K}$ is a ${K}\langle \partial \rangle$-module of finite length (equivalently, of finite ${K}$-dimension).} $N = {K}\langle \partial \rangle/ {K}\langle \partial \rangle  P$, and is called the {\it irregularity} of $N$ and denoted by ${\rm ir}\,N$. 
 
  Let us define $\mu(N)$ to be $ \frac{ir(M)}{\rk N}\,$ and say that $N$ is {\it isoclinic}  if and only if for any non-zero differential submodule $M$, $\mu(M)= \mu(N)$.
     
\smallskip Any regular differential module is a successive extension of rank one (regular) differential modules of the form ${K}\langle \partial \rangle/ {K}\langle \partial \rangle  (\partial -c),\;c\in \C$.

 \smallskip  Any differential module $N$ has a unique descending filtration - the Turrittin-Levelt filtration -
   $$ 0 \subset   F^{\geq \la_1} \, N  \subset \ldots \subset       F^{\geq \la_r} \, N  = N$$ 
for which $\la_1> \ldots > \la_r$, and  $gr^{\la_i}\,N$ is isoclinic of slope $\la_i$.  In fact, the filtration splits canonically (cf \cite{Tu}\cite{L}, and also \cite{Mal}\cite[\S 2]{A5}).  

Moreover, ${\rm ir}\, N$ coincides with the degree attached to this filtration in the sense of \ref{0.3}; the highest slope of the Newton polygon is called the {\it Poincar\'e-Katz rank} of $N$.

\smallskip Simple differential modules $N$ can be described explicitly: ${\rm ir}\, N$ is prime to  $r=\rk N$, and $N$ is induced by a rank one differential module over ${K}' = \C((x^{1/r}))$, of the form 
$${K}'\langle \partial \rangle/ {K}'\langle \partial \rangle  (\partial -f),\;f\in {K}'[x^{-1/r}],\; \deg_{x^{-1/r}} f= {\rm ir}\, N.$$

\bigskip
 \subsection{Hasse-Arf filtration of local Galois representations.}\label{haf} 
 
Let $  ({K}, v)$ be a complete discretely valued field with perfect residue field $k$, and let $G_{{K}}=Gal({{K}}^{sep}/{{K}})$ be its absolute Galois group. By analysing the ``norm" of $g - id$ acting on finite extensions ${ L}/{{K}}$, ramification theory provides a decreasing sequence of normal subgroups
 $$G_{{K}}^{(\la)}\triangleleft G_{{K}}, \, \la \in \mathbb Q_+.$$ 
Let $F$ be a field of characteristic zero, and let $M$ be a $F$-linear representation of $G_{{K}}$ with finite image. Then the filtration $G_{{K}}^{(\la)}$ gives rise to a descending filtration of $M$ indexed by rational numbers  - the Hasse-Arf filtration.  In fact, the filtration splits canonically. 

 The degree attached to this filtration in the sense of \ref{0.3} is the so-called {\it Swan conductor} of $M$. This is an integer (Hasse-Arf theorem  \cite{Ha}\cite{Arf}, \cf also \cite[IV, VI]{Se}).

\smallskip To be more concrete, consider the case ${{K}}=k((x))$. If 
 ${\rm char}\, k=0$, then elements of ${{K}}^{sep}$ are just Puiseux series, and the Hasse-Arf filtration is trivial.  If ${\rm char}\, k=p>0$,  Puiseux series\footnote{with $p$-integral exponents, by separability.}  form only the so-called {\it tame} part ${{K}}^{tame}$ of ${{K}}^{sep}$ (for instance, there is no solution of the Artin-Schreier equation $y^p-y=1/x$ in terms of Puiseux series); the {\it wild} subgroup of $G_{{K}}$, $$Gal({{K}}^{sep}/{{K}}^{tame})= \bigcup_{\la>0}\,G_{{K}}^{(\la)} , $$ is a pro-$p$-group.  
 
 \bigskip
   \subsection{Grayson-Stuhler filtration of euclidean lattices.}\label{gsf}

   Let $N$ be a euclidean lattice, \ie a $\Z$-lattice together with a scalar product on its real span.  Its degree is defined by
   \begin{equation}\,\deg \, N = - \log \,vol(N\otimes \R/N) .\end{equation} 
   
   If $N\neq 0$, its slope is the ratio 
   $\, \mu(N)= \frac{\deg N}{ \rk N}\, $.
   
  Any euclidean lattice $N$ is a successive extension of semistable lattices of increasing slopes: more precisely, 
  $N$ has a unique descending filtration - the Grayson-Stuhler filtration \cite{Gra}\cite{St} -
   $$ 0 \subset   F^{\geq \la_1} \, N  \subset \ldots \subset       F^{\geq \la_r} \, N  = N$$ 
for which $\la_1> \ldots > \la_r$, and the graded pieces $ gr^{\la_i}\,N = F^{\geq \la_i}\,N/F^{> \la_i}\,N$ are semistable bundles of slope $\la_i$.  The breaks are related to the successive minima in the sense of Minkovski's geometry of numbers \cite{Bo}.
 
\medskip In the previous examples, the underlying categories were additive (vector bundles, $\phi$-modules, differential modules, Galois representations). Here, this is no longer the case:  in the underlying category of euclidean lattices, morphisms are additive maps of norm $\leq 1$.

   \smallskip \bigskip
 
\section{Proto-abelian and quasi-abelian categories.}\label{sII1}  

 Our aim is to study slope filtrations independently of the particular context in which they arise.

In order to do so, one is at once faced to the problem of chosing a class of categories which covers the majority of examples in the literature, without being too general. 
 As is shown by the first basic example, abelian categories are not enough (vector bundles on a curve do not form an abelian category).   
 
 A convenient class of additive categories to work with is the class of {\it quasi-abelian categories} (\cf \ref{qac}). It allows to treat all concrete examples of slope filtrations on additive categories. 
 
 However, it is too restrictive, since it excludes the category of euclidean lattices with contracting morphisms, and other non-additive categories arising in Arakelov geometry. 
 
 A close inspection of the logical network involved in each example shows that it is indeed possible to drop additivity, and that the right class of categories to consider in order to develop a general theory of slope filtrations is the class of  {\it proto-abelian categories} (a non-additive version of quasi-abelian categories, \cf \ref{pac}).

\bigskip\subsection{Kernels and cokernels.}\label{kcok}

We begin with three reminder subsections, using MacLane's terminology \cite[p. 191]{ML}. 

 Let $\sC$ be a category with a {\it null object} $0$, \ie an object that is both initial and terminal. 
 For any pair $M,N$ of objects of $\sC$, the {\it $0$ morphism} is the composed morphism $M\to 0\to N$. 
 
\medskip For any morphism $M\overset{f}{\to} N$, a {\it kernel} $\ker f$ of $f$ is a morphism with codomain $M$ such that $f\circ \ker f=0$, that is universal for this property (hence unique up to unique isomorphism). 

By common abuse of language, one also calls ``kernel of $f$" the domain of $\ker f$ (which we denote by $\Ker f$ in order to prevent confusion).

Any monic\footnote{\ie left cancellable.} has kernel $0$ (the converse is not true in general).

Any kernel is monic, and is called {\it strict monic}\footnote{some authors say ``admissible" or ``normal" instead of ``strict".}; its domain is called a {\it strict subobject} of its codomain.

 \medskip 
 Dually, a cokernel $\coker f$ of $f$ is a morphism with domain $M$ such that $(\coker f) \circ f=0$, that is universal for this property (hence unique up to unique isomorphism). 
 By common abuse of language, one also calls ``cokernel of $f$" the codomain of $\coker f$ (which we denote by $\Coker f$).

Any epi\footnote{\ie right cancellable.} has cokernel $0$.

Any cokernel is epi, and is called {\it strict epi}; its domain is called a {\it strict quotient} of its domain.

\medskip A \emph{short exact sequence}\footnote{some authors say ``strictly exact".}, denoted by
$$ 0\to M  \overset{f}{\to} N \overset{g}{\to} P \to 0,$$
is a pair $(f,g)$ of composable morphisms such that
$$f= \ker g,\;\; g=\coker f.$$ One says that $N$ is an \emph{extension} of $P$ by $M$, and one writes $P= N/M$.
 
 \medskip A functor is  \emph{exact} if it preserves short exact sequences.

 \medskip\subsection{Categories with kernels and cokernels.}\label{ckc} 

Let $\sC$ be a {\it category with kernels and cokernels}, \ie with a null object, and such that any morphism has a kernel and a cokernel.
 For a morphism $M\overset{f}{\to} N$, one sets  $$\coim f = \coker\ker f,\;\;\im f= \ker\coker f$$ and one denotes the codomain of $\,\coim f $ by $\Coim f$ and the domain of $\im f$ by $\IM f$ or $f(M)$.
 
\smallskip One then has (\cite[p. 193]{ML})
$$ \coker\im f=\coker f,\;\;\ker\coim g= \ker g,$$ whence the equivalence, for morphisms $M\overset{f}{\to} N,  \,N \overset{g}{\to} P $, between
\begin{itemize}
\item 
$f$ is strict monic (\resp $g$ is strict epi)
\item $f=\im f$ (\resp $g=\coim g$),
\item there is a morphism $N\overset{g}{\to} P$ (\resp $M\overset{f}{\to} N$) such that  $ \, 0\to M  \overset{f}{\to} N \overset{g}{\to} P \to 0 \,$   is a short exact sequence.  \end{itemize}

\medskip Any  $f$ has a unique factorization as
$$f = \im f \circ \bar f \circ \coim f $$ 
(where $\bar f$ may have non-zero kernel or cokernel in general\footnote{this occurs even in the additive case, \cf \cite{Rum}.}). For any factorization
$$f=  m \circ f' \circ e$$ where $e$ is strict epi and $m$ is strict monic, there are unique factorizations (\cf  \cite[p. 193]{ML})
$$\coim f = eg,\, \im f= hm, \; \bar f = h f' g .$$

\medskip Let $M\overset{f}{\to} N  \overset{g}{\to} P$ be composable morphisms. One has:
\begin{itemize}\item $\ker g=0\Rightarrow \ker gf =\ker f\,$ and $\,\coim gf=\coim f$,
\item  $\coker f=0\Rightarrow \coker gf=\coker g\,$ and $\,\im gf=\im g$.
 \end{itemize}

   \bigskip
\subsection{Pull-backs and push-outs.}\label{pbpo} 

\medskip  Let $P\overset{f}{\to}Q \overset{g}{\leftarrow} N$ be a pair of morphisms with a common codomain. A {\it pull-back} square (or cartesian square) is a commutative square 
\begin{equation}\label{CD1}\begin{CD} M @>f'>>  N  \\ @Vg'VV @VVgV   \\    P  @>f>>   Q,  \end{CD}\end{equation}
built on $(f,g)$, that is universal. One says that $f'$ (\resp $g'$) is the pull-back of $f$ by $g$ (\resp $f$). Pull-back squares may be composed.

Dually, for a pair $P \overset{g'}{\leftarrow} M \overset{f'}{\to} N$ of morphisms with a common domain, one has the notion of {\it push-out}. 

\smallskip Let us assume that $\sC$ has kernels and cokernels. 
 Then for any pull-back square \eqref{CD1}, the natural morphisms 
$$ \Ker f' \to \Ker f ,\;\; \Ker g'\to \Ker g$$ are isomorphisms (inverses are provided by the universal property). 

  Dually, for any push-out square \eqref{CD1}, the natural morphisms 
$$ \Coker f' \to \Coker f ,\;\; \Coker g'\to \Coker g$$ are isomorphisms.

\medskip The {\it pull-back of a strict monic $f$ always exists} and is strict monic. 

Indeed, take $f'= \ker ((\coker f) \circ g)$ (which is strict monic); then because $f=\ker\coker f$, there is a canonical factorization $g\circ f'=f\circ g'$, and any morphism $h:L\to N$ such that $gh$ factors through $f$ satisfies $(\coker f) \circ g \circ h=0$, hence factors through $f'$. This shows that $f'$ is the pull-back of $f$.
  One has $\, \ker gf'= \ker g'$ and $\, \coim gf' =\coim g'$.
  
\medskip  If $g$ is also strict monic, one writes $M= N\cap P$. If $Q\inj Q'$ is monic and $N, P$ are strict subobjects of $Q'$, the pull-back of $P {\to}Q'  {\leftarrow} N$ is $N\cap P$.

One has
$$N\cap P = \Ker (N \to Q/P)= \Ker (P\to Q/N) .$$

\medskip Dually, the {\it push-out of a strict epi $g'$ always exists} and is strict epi: $g= \coker(f'\circ ker g')$. One has $\, \coker gf'= \coker g'$ and $\, \im gf' =\im g'$.
  
  If $g'$ is also strict epi, so that $N=M/N',\, P=M/P'$, one writes (abusively) $Q=M/(N'+P')$. 
    One has
$$M/(N'+P') =  \Coker (N' \to M/P')= \Coker (P'\to M/N') .$$

\begin{lemma}\label{mmmeee} Assume that 
$\sC$ has kernels and cokernels. \begin{enumerate}
\item 
Let $M\overset{f}{\to} N  \overset{g}{\to} P$ be composable morphisms. 
\\ If $gf$ is strict monic and $g$ is monic, then $f$ is strict monic.
\\ If $gf$ is strict epi and $f$ is epi, then $g$ is strict epi.
\item  Any pull-back square \eqref{CD1}
 in which $g$ is strict epi and $g'$ is epi is also a push-out square. Dually, any push-out square \eqref{CD1}
 in which $f'$ is strict monic and $f$ is monic is also a pull-back square. 
\end{enumerate}
\end{lemma} 

\proof (1) By duality, it suffices to treat the first case. Let us consider the pull-back square 
\begin{equation}\begin{CD} L @>h'>>  N  \\ @Vg'VV @VVgV   \\    M  @>gf>>   P.  \end{CD}\end{equation}
 Since  $g$ is monic, $gh'= gf  g' $ implies $h' = fg'$, and since $h'$ is monic, so is $g' $.  On the other hand, applying the pull-back property to $(id_M, f)$, one gets a right-inverse to $g'$, hence $g'$ is an isomorphism. Therefore $f$ is strict monic like $h' = fg'$. 
 
 (2) If \eqref{CD1} is a pull-back square, the natural morphism  $ Ker g'\overset{f''}{\to} Ker g$ is an isomorphism. In particular, for any pair of morphisms $P\overset{u}{\to}Q' \overset{v}{\leftarrow} N$ such that $ug'=vf'$, the composition $v \circ \ker g = vf' \circ \ker g' \circ (f'')^{-1} $ is $0$, hence $v$ factors through $  \Coker \ker  g$, which is $g$ since $g$ is strict epi. Let us write $v= w  g$ and set $u'= w   f$. Then $u'  g'= v   f' = ug'$, and since $g'$ is epi, $u= u' = wf$. This shows that \eqref{CD1} is a push-out square.
  \qed

 \begin{exs}\label{exproto} 1) The category of groups has kernels and cokernels, and even pull-backs and push-outs. One has:
$\,$ monic = injective, $\,$  strict epi = epi = surjective, $\,$ strict subobject = normal subgoup (\cf \cite[ex. 5, p. 21]{ML}).
 For composable morphisms $G\overset{f}{\to} G'  \overset{g}{\to} G''$, \begin{itemize}
  \item if $gf$ is strict monic, $f$ is monic but not strict in general,
\item if $f$ and $g$ are strict monic, $gf$ is monic but not strict in general,
\item if $gf$ is strict monic and $g$ is monic, $g$ is not strict in general.
 \end{itemize}

\smallskip\noindent 2) The category of hermitian (finite-dimensional real or complex) vector spaces,  with linear maps of norm $\leq 1$ as morphisms,  has kernels and cokernels. One has: $\,$  monic = injective, $\,$ epi = surjective. A subobject (\resp quotient) is strict if its norm is the induced (\resp quotient) norm.

This category has finite coproducts (the usual orthogonal sum) and even push-outs. But the self product of a non-zero object does not exist (since the diagonal map has norm $> 1$); a fortiori, pull-backs do not exist in general in this category. 

 \smallskip\noindent 3)  The category of euclidean lattices (with additive maps of norm $\leq 1$ as morphisms)  has kernels and cokernels. One has: $\,$  monic = injective, $\,$ epi = surjective on the real span. A subobject is strict if it is cotorsionfree and if its norm is the induced  norm, a quotient is strict if its norm is the quotient norm. Any epi-monic is the composition (in either order) of an isometric epi-monic and a morphism which is identity on the underlying lattice.
    \end{exs}

  \medskip
\subsection{Proto-abelian categories.}\label{pac} 

\begin{defn} A category $\sC$ with kernels and cokernels is \emph{proto-abelian} if 
\begin{enumerate}
  {\it\item any morphism with zero kernel (\resp cokernel) is monic (\resp epi),
\smallskip\item the pull-back of a strict epi by a strict monic is strict epi, and the push-out of a strict monic by a strict epi is strict monic.  }
   \end{enumerate}\end{defn} 

Axiom (2) allows to deal with {\it strict subquotients} without ambiguity.

\begin{exs}\label{exproto2} 1)  Any abelian category is proto-abelian. In fact, a proto-abelian category is abelian if and only if it has finite products and coproducts, and any epi-monic is an isomorphism\footnote{indeed, a category with kernels and cokernels is abelian if and only if it has finite products and coproducts and any morphism with zero kernel and cokernel is an isomorphism (this implies additivity), \cf \cite[p. 201]{ML}.}.

\smallskip\noindent 2) The category of (finite dimensional) hermitian vector spaces is proto-abelian: on a subquotient space, the quotient norm of the induced norm is the norm induced by the quotient norm.

\smallskip\noindent 3) The category of euclidean lattices is proto-abelian. 

\smallskip\noindent 4) On the other hand, the category of groups fails to be proto-abelian: it satisfies neither (1) nor (2).
\end{exs} 

\medskip Let $\sC$ be a proto-abelian category. 
  
\begin{lemma}\label{pabc} \begin{enumerate} 
\item  any \emph{pull-back} square
\begin{equation}\label{CD2}\begin{CD} M @>f'>>  N  \\ @Vg'VV @VVgV   \\    P  @>f>>   Q  \end{CD}\end{equation}
in which $f$ is strict monic and $g$ is strict epi is also a \emph{push-out square} in which $g'$ is strict epi and $f'$ is strict monic, and conversely. It extends to a commutative diagram with exact rows and columns

\begin{equation}\label{CD3}\begin{CD}
 0 @>>> 0 @>>>  0 @>>>  0 @>>>  0
\\ @VVV @VVV @VVV @VVV @VVV
 \\  0 @>>> \Ker g' @>f''>> \Ker g @>>> 0 @>>>  0 
 \\ @VVV @VVV @VVV @VVV @VVV 
 \\  0 @>>> M @>f'>>   N @>>>  N/M @>>>  0
 \\ @VVV @Vg'VV @VVgV @VVV @VVV
  \\  0 @>>> P @>f>> Q @>>> Q/P @>>>  0
  \\ @VVV @VVV @VVV @VVV @VVV
   \\  0 @>>> 0 @>>>  0 @>>>  0 @>>>  0. \end{CD}\end{equation}
 \smallskip\item  If $M\overset{f}{\to} N\overset{g}{\to} P$ are strict monic (\resp strict epi), so are $\,M\overset{gf}{\to} P$ and $N/M {\to} P/M$ (\resp and $\Ker (M\to P)\to \Ker (N\to P)$). In fact, one has  short exact sequences  
 \begin{equation}\label{SESmonic}\;0\to N/M {\to} P/M \to P/N\to 0\end{equation}
 \begin{equation}\label{SESepi}   \;0\to \Ker(M\to N) \to \Ker (M\to P)  {\to} \Ker(N\to P)\to 0 .\end{equation}
  
\smallskip\item If the pair of composable morphisms $M\overset{f}{\to} N\overset{g}{\to} P$ satisfies $gf=0$, there is a short exact sequence $$ 0\to \Ker g/f(M)  \to \Coker f  \to  \Coim g\to 0.
$$ 
 \smallskip \item  In the canonical factorization 
$$ f = \im f \circ \bar f \circ \coim f  $$ of any $f$, $\bar f$ is epi-monic.
   \end{enumerate}\end{lemma}

 \proof (1) follows from item (2) of \ref{mmmeee}, duality, and the second axiom of proto-abelian categories.  
      
\smallskip\noindent (2). Let $M\overset{f}{\to} N\overset{g}{\to} P$ be strict monic, and let us consider the push-out
 \begin{equation}\label{CD4}\begin{CD} N @>g>>  P  \\ @Vh\,=\,\coker fVV @VVh'V   \\    N/M  @>g' >>   Q.  \end{CD}\end{equation}
By item (1), $g$ induces an isomorphism $M= \Ker h \cong \Ker h'$, hence the sequence 
$$0\to M  \overset{gf}{\to} P \overset{h'}{\to} Q\to 0$$ is exact, which proves that $gf$ is strict, as well as $N/M \overset{g'}{\to} Q = P/M$. In fact, since \eqref{CD4} is a push-out square, $\Coker g \cong \Coker g'$, which gives the short exact sequence \eqref{SESmonic}.

 The other part of the assertion follows by duality.

 \smallskip\noindent (3). Since $\coker g \circ g \circ f$ is zero and factors through the epi $\im f$, $f(M)\to \Coker g $ is zero, hence $f(M)\to N$ factors through $g= \ker\coker g$. Item (3) then follows from item (2) applied to $f(M) \to \Ker g \to N$.

 \smallskip\noindent (4). By item (2),  $\coim \bar f \circ \coim f$ is strict epi, hence equal to $\coim f$ by universality of the canonical factorization. This implies $\ker \bar f=0$.  By duality, $\coker \bar f= 0$. By axiom (1), $\bar f$ is therefore epi-monic.
    \qed

\medskip
 \subsection{Flags.} Let $\sC$ be a proto-abelian category.
 
  \begin{defn}\label{FLAG} A \emph{flag}  of length $r$ on $N$ is a finite sequence  $$\sF:\; 0= N_0 \inj  N_1 \inj   \cdots \inj N_r = N$$ of strict subobjects of $N$, with $N_i\neq N_{i-1}$ for $1\leq i\leq r$.   \end{defn} 

Note that by lemmas \ref{mmmeee} and \ref{pabc}, it amounts to the same to say that the $N_i$ are strict subobjects of $N$, or to say that  $N_i\inj N_{i+1}$ are strict monic.  In particular, it makes sense to consider the graded pieces $N_i/N_{i+1}$. 

\smallskip The following lemma will allow us to make some induction arguments.

  \begin{lemma}\label{PB} Let $$0 \to M\to N\overset{e}{\to} P\to 0$$ be a short exact sequence with $M\neq 0$.
  \begin{enumerate}
  \item Let $$\sF:  0\inj   P_1 \inj   \cdots \inj P_r = P$$ be a flag of length $r$ on $P$. Then the (step-by step) pull-back $$ e^\ast\sF: 0\inj   N_1 \inj   \cdots \inj N_{r+1} =N$$ is a flag of length $r+1$ on $N$, and $N_1= M,\, P_{i}= N_{i+1}/N_1$.   
    
  \smallskip \item Conversely, let $$ \sG: 0\inj  N_1 \inj   \cdots \inj N_{r+1} =N$$ be a flag of length $r+1$ on $N$ with $N_1=M$. Then the (step-by step) push-out $$ e_\ast\sG: 0\inj   P_1 \inj   \cdots \inj P_{r+1} =P$$ is a flag of length $r $ on $P$, and $  P_{i}= N_{i+1}/N_1$.   
\end{enumerate} \end{lemma} 
 \proof This follows from lemma  \ref{pabc}. \qed

   \begin{lemma}\label{REFI}
   Any two flags (of equal or unequal length) on $N$ admit a common refinement.  \end{lemma} 
   
   \proof   A common refinement of  the flags $(N_i)$ and $(N_j)$  is given by the following non-decreasing sequence of strict subobjects of $N$ (with respect to the lexicographic order):
  $$N_{ij} := (N_i\cap N'_j) + N_{i-1}   $$  
 (defined by the push-out of $N_i\cap N'_j\to N\leftarrow N_{i-1}$,  $N_i\cap N'_j$ being the pull-back of $N_i\leftarrow N \rightarrow N'_j$). \qed

\medskip\subsection{Rank function.}

Assume that $\sC$ is essentially small, let  $\,\sk\, \sC$ be the set  of isomorphism classes of objects of $\sC$ (skeleton). Taking $\sk\, \sC$ as set of generators and short exact sequences as relations, one builds the Grothendieck group $K_0(\sC)$. Any element of $K_0(\sC)$ can be written as a difference $[M]-[N]$ of elements of $\sk \, \sC$.

\begin{defn}\label{RANK} A \emph{rank function} on $\sC$ is a function $$\,\rk\,: \, \sk\, \sC \to \N$$ that 
 is additive on short exact sequences and takes the value $0$ only on the $0$ object. 
 \end{defn} 

In particular, a rank function gives rise to a group homomorphism still denoted by
 $$\rk: K_0(\sC) \to \Z.$$ 
 
\begin{rem}\label{FINITELENGTH} The length of any flag on $N$ is bounded by $\rk N $. It follows that {\it any abelian category with a rank function is noetherian and artinian}, and that the Jordan-H\"older length (given, for any object, by the maximal length of a flag on this object) is a rank function. Any non-zero subobject of $N$ of minimal rank is simple. 
 \end{rem}
  
 \begin{ex} If $\sC$ is the category of finitely generated torsionfree modules over a domain $R$, the usual rank (\ie the dimension of the vector space obtained by tensoring with the fraction field of $R$) provides a rank function. 
  \end{ex}

\medskip\subsection{The additive situation: quasi-abelian categories.}\label{qac}

 Recall that an additive category with kernels and cokernels has all finite limits and colimits, in particular all pull-backs and push-outs \cite[p. 113]{ML}.  
 
 \begin{defn} An additive category with kernels and cokernels is \emph{quasi-abelian} if {\it every pull-back of a strict {epi} is strict {epi}}, and  {\it every push-out of a strict monic is strict monic}.
 \end{defn} 
 
If $\sC$ is essentially small, this amounts to requiring that the set $Ext(P,M)$ of isomorphism classes  of extensions of an object $P$ by an object $M$ is bifunctorial.  
 
\medskip It is immediate that  {\it any quasi-abelian category is proto-abelian\footnote{we do not know if, conversely, any additive proto-abelian category is quasi-abelian.  
 }.} In particular, in the canonical factorization of any morphism $f= \im f \circ \bar f\circ \coim f$, $\bar f$ is always epi-monic (item (4) of lemma \ref{pabc})\footnote{as we already mentioned, this property is not true in any additive category with kernels and cokernels \cite{Rum}}. 

 \begin{exs} 1) The category of torsionfree finitely generated modules over any domain $R$  is quasi-abelian.
 
 If $R$ is Dedekind (or more generally Pr\"{u}fer), this is the category of projective modules of finite rank. If $R$ is principal (or more generally B\'ezout), this is the category of free modules of finite rank. 
  
 \smallskip \noindent 2) The category of (finitely generated) reflexive modules over an integrally closed domain $R$ is quasi-abelian. Kernel and cokernels in this category are the double duals of kernels and cokernels taken in the category of $R$-modules.
 
 If $R$ is regular of dimension $2$, this is the category of projective modules of finite rank.
 
 \smallskip \noindent 3) The category of torsionfree coherent sheaves over a reduced irreducible analytic space or algebraic variety $X$ is quasi-abelian. If $X$ is a normal curve, this is the category of vector bundles (of finite rank).
 
  \smallskip \noindent 4) The category of reflexive coherent sheaves over a normal analytic space or algebraic variety $X$ is quasi-abelian. 
 
 \smallskip \noindent 5) The category of filtered modules over any ring is quasi-abelian.

   \smallskip  Beside these algebro-geometric examples, there are many examples from functional analysis:  various categories of topological vector spaces - Banach and 
Fr\'echet spaces, locally convex and nuclear spaces, bornological spaces of convex type - are 
quasi-abelian. \end{exs}

\medskip The notion of quasi-abelian category seems to go back to Yoneda's 1950 paper \cite{Y}, and has been rediscovered a number of times (with various names, viewpoints and languages\footnote{\cf  \cite{Ju}\cite{Rai}; the adjective quasi-abelian, with this meaning, seems to stem from \cite{SC}. We refer to \cite[\S 2]{Rum} for a short history of this notion}. From the definition, a quasi-abelian category is just an exact category with kernels and cokernels, in which short exact sequences are defined as above. 

A systematic exposition is due to J.-P. Schneiders \cite{Sc}. In \cite{BvB}, it is shown that this notion is equivalent to the notion of cotilting torsion pairs introduced by D. Happel, I. Reiten and S. Smalo in representation theory\cite{HRS}. The main result can be summarized as follows.

\begin{prop}\label{P1}\cite[\S 1.2.4]{Sc}\cite[5.4, app. B]{BvB} An additive category $\sC$ is quasi-abelian if and only if it can be fully embedded in an abelian category ${\sA}$ with the following properties:
\begin{enumerate}
\item any object of ${\sA}$ is a quotient of an object of $\sC$,
\item there is a strictly full\footnote{\ie full and closed by isomorphism.} subcategory $\sT\subset  {\sA}$ (the ``torsion subcategory")  such that
 \begin{itemize}
 \item any object $A$ of ${\sA}$ sits in a unique (up to unique isomorphism) short exact sequence 
$$0\to A_{tor}\to A\to M\to 0$$ where $\,A_{tor}\in Ob \, \sT$ and $M\in Ob \, \sC$,
\item  there is no non-zero morphism from objects of $\sT$ to objects of $\sC$.
\end{itemize}
\end{enumerate}
 Condition (2) implies that any subobject of $M$ in ${\sA}$ is (isomorphic to an object) in $\sC$.
 Together with (1), it characterizes the pair $({\sA}, \sT )$. 
 
 A short sequence in $\sC$ 
$$0\to M_1 \to M \to M_2\to 0$$ is exact in $\sC$ if and only if it is exact in ${\sA}$. A morphism in $\sC$ is monic (\resp strict {epi}) if and only if it is monic (\resp {epi}) in $\sA$. 
  \end{prop}

 In the sequel,  $\sA$ will be called the \emph{left abelian envelope} of $\sC$: in \cite{Sc}, ${\sA}$ appears as the heart of the derived category of $\sC$ with respect to the ``left $t$-structure", for which  $ D(\sC)^{\leq 0}$ is represented by complexes in degree $\leq 0$, and $D(\sC)^{\geq 0}$ by complexes in degree $\geq -1$, the morphism $d_{-1}$ being monic\footnote{there is also a right $t$-structure whose heart, the ``right abelian envelope" of $\sC$, is the ``tilting" of $\sA$, \cf \cite{Sc}\cite{BvB}.}. 
 
 The canonical embedding 
 $\,\sC \to {\sA}\,$ has a left adjoint and induces an equivalence 
 $\, D(\sC) \cong D({\sA})\,$ compatible with the $t$-structures (the left one on the left-hand side, the canonical one on the right-hand side), hence an equivalence of categories with bounded $t$-structures
 $\, D^b(\sC) \cong D^b({\sA})\,$ (actually, this construction already appears in \cite[1.3.22]{BBD}). Any object of ${\sA}$ is represented by the complex $[M\to N]$, in degrees $-1$ and $0$, associated to a monic in $\sC$. One deduces that there is a canonical isomorphism
$$K_0(\sC)\cong K_0({\sA}).$$
 
    \begin{rem} In particular, any rank function $\rk$ extends to a function $$ \sk\, \sA \to \Z,$$ and any object of $\sA$ of rank $0$ is in $\sT$. In \ref{qtrf}, we will deal with a class of quasi-abelian categories in which any object of the left abelian envelope has non-negative rank, hence objets of $\sT$ are precisely the objects of $\sA$ of rank $0$. 
 \end{rem}

\begin{exs} 1) If $\sC$ is the category of finitely generated torsionfree modules over a domain $R$, 
  ${\sA}$ is the abelian category of all finitely generated modules, and $\sT$ the full subcategory of torsion modules; a cokernel in $\sC$ is a cokernel in ${\sA}$ divided out by torsion.
 
 \smallskip \noindent 2) If $\sC$ is the category of reflexive modules over an integrally closed domain $R$, ${\sA}$ is the localization of the abelian category of all finitely generated modules with respect to pseudo-isomorphisms \cite[ch. 7]{B}, \ie morphisms which are isomorphisms outside a closed subset of codimension $2$ of $\Spec R$.
 \end{exs}

\begin{lemma}\label{TORS} An object of $\sA$ belongs to $\sT$ if and only if it is the cokernel in $\sA$ of a epi-monic in $\sC$.  
\end{lemma}

\proof Let $Q$ be an object of $\sT$. By \ref{P1}, there is an epi $N\to Q$ and an epi $M\to \Ker_\sA(N\to Q)$, with $N,M\in Ob\,\sC$. 
Replacing $M$ by its $\sC$-image in $N$, one may assume that $M\overset{f}{\to} N$ is monic. By construction $Q = \Coker_\sA f$. Let $N\overset{g}{\to}P$ be a morphism in $\sC$ such that $gf=0$. Then $g$ factors through $Q$, and since $P$ has no torsion, $g=0$. Hence $f$ is epi.

Conversely, let us write $Q= \Coker_\sA f$ as an extension of $P\in Ob \, \sC$ by $R\in Ob\, \sT$. The composed morphism $M\overset{f}{\to} N\to P $ is zero. Since $f$ is epi,  $N\to P$ is zero, \ie $\coker f$ factors through $R$. Hence $Q=R$. \qed

\begin{lemma} Any strict subobjects $N, P$ of an object $Q$ of $\sC$ are also strict subobjects of their sum $N+P$ (in the sense of \ref{pbpo}) and the natural morphism
 $$N/(N\cap P) \to (N+P)/P$$ is an isomorphism.  In fact, one has a commutative diagram with exact rows and columns
\begin{equation}\begin{CD}\label{CDqab} 0 @>>> 0 @>>>  0 @>>>  0 @>>>  0\\ @VVV @VVV @VVV @VVV @VVV \\  0 @>>> N\cap P @>>>  N @>>>  N/(N\cap P) @>>>  0 \\ @VVV @VVV @VVV @VVV @VVV \\  0 @>>> P @>>>  Q @>>>  Q/P @>>>  0\\ @VVV @VVV @VVV @VVV @VVV \\  0 @>>> P/(N\cap P) @>>>  Q/N@>>> Q/(N+P) @>>>  0\\ @VVV @VVV @VVV @VVV @VVV \\  0 @>>> 0 @>>>  0 @>>>  0 @>>>  0. \end{CD}\end{equation}
 \end{lemma}

\proof Indeed,  $N\cap P $ is the kernel of $M {\to} M/N \oplus M/P,\; N+P$ is the image of $N\oplus P  {\to} Q$, and the natural morphism $N/(N\cap P) \to (N+P)/P$ in $\sC$ is an isomorphism since it is so in $\sA$; \eqref{CDqab} has exact rows and columns in $\sA$, hence in $\sC$. \qed
 
 \begin{rem} This holds does not in a proto-abelian category in general (in the category of hermitian vector spaces, a counterexample is constructed by taking $N$ and $P$ to be non-orthogonal supplementary subspaces).\end{rem}

 \bigskip
   \section{Slope functions and the ``yoga" of semistability.}\label{s1}

\medskip In this section, we introduce the yoga of (semi)stability with respect to a slope function $\mu$ in a proto-abelian category.

   \bigskip In the sequel, {\it $\sC$ stands for an essentially small proto-abelian category equipped with a rank function $\rk$}.    

\medskip  In addition, we fix {\it a totally ordered, uniquely divisible, abelian group $\La$} (in practice, this will be $\Q$, or a $\Q$-subspace of $\R$).

 \medskip\subsection{Slope functions and degree functions. }     
 
  \begin{defn}\label{sf} A \emph{slope function} on $\sC$, with values in $\La$, is a map
  $$\mu:  \, \sk\, \sC \setminus \{0\} \to \La$$ that satisfies the following two conditions:
  \begin{enumerate}
   \item for any epi-monic $M\to N$, one has $\mu(M) \leq \mu(N)$,
 
 \smallskip \item the associated \emph{degree function} 
  $$ \deg = \mu \cdot \rk: \,\sk\, \sC \to \La $$ (taking value $0$ at the $0$ object) is {\it additive on short exact sequences}.  
     \end{enumerate} 
   \end{defn}
 Of course $\mu$ and $\deg$ determine each other, and the latter induces a group homomorphism  $$\deg: K_0(\sC) \to \La .$$  
    
 \begin{rems}\label{remsf} {1)} If $\sC$ is abelian, condition (1) is empty. If $\sC$ is abelian semisimple, to give a slope function amounts to attaching to every simple object (up to isomorphism) a label in $\La$. 
 
\smallskip\noindent {2)} If $\sC$ is abelian, one can define, up to isomorphism, the semi-simplification $N_{ssi}$ of any object $N$. Then $\mu(N)=\mu(N_{ssi})$ by additivity of the degree.
  
\smallskip\noindent {3)}  If $\mu$ is a slope function on $\sC$, then $-\mu$ is a slope function on $\sC^{op}$.  

\smallskip\noindent {4)} Degree functions form a convex cone in the $\Q$-vector space  $Hom(K_0(\sC), \La)$.
   \end{rems}

 \begin{exs}\label{exs1} 1) Let $\sC$ (\resp ${\sA}$) be the category of torsionfree coherent sheaves (\resp all coherent sheaves) on a polarized normal connected projective variety $(X, \sO(1))$. Take $\La= \Q$. Then the quotient 
 $\,\mu = \frac{\deg}{\rk},$ where  
 $\deg$ and $\rk$ have their usual meaning, is a slope function on $\sC$ (Mumford-Takemoto \cite{Mu}\cite{Ta}\footnote{introduced by Mumford  for $\dim X=1$,  Takemoto generalized to $\dim X >1$.}). 
 
 When $\dim X>1$, a more refined choice of slope function on $\sC$ is the following (Gieseker-Maruyama), which is more useful in moduli problems \cf \eg \cite{Gi2}. Take $\La = \Q[t]$, with the total order given by $P\geq Q$ if $P(t)\geq Q(t)$ for $t>>0$\footnote{this is the lexicographical order on the coefficients.}. The function which associates to any non-zero torsionfree coherent sheaf its Hilbert polynomial divided by its rank is a slope function on $\sC
$ (this follows from the fact that the Hilbert polynomial is additive on short exact sequences in ${\sA}$, being an Euler characteristic in the large, and has non-negative leading coefficient).

\smallskip 2) Let $\sC$ be the category of euclidean lattices (with additive maps of norm $\leq 1$ as morphisms), \cf \ref{exproto} 3). Take $\La= \R$. Then the quotient 
 $\,\mu = \frac{\deg}{\rk},\,$ where $\deg$ is minus the logarithm of the covolume, is a slope function on $\sC$ (Grayson-Stuhler \cite{Gra}\cite{St}).

     \end{exs}   
 
  \begin{lemma}\label{L1}   \begin{enumerate}
     \item For any short exact sequence $0\to M\to N\to P\to 0$ of non-zero objects, one has $$\min((\mu(M), \mu(P)) \leq \mu(N)\leq \max(\mu(M), \mu(P)),$$ both inequalities being strict unless  $\mu(M)=\mu(N)=\mu(P)$.
 
  \smallskip  \item More generally,
for any flag $0= M_0 \inj  M_1 \inj   \cdots \inj M_r  = M$ with non-zero quotients $M_i/M_{i-1}$, one has 
$$ \min(\mu(M_i/M_{i-1}))\leq \mu(M)\leq \max(\mu(M_i/M_{i-1})),$$ both inequalities being strict unless all the   $\mu(M_i/M_{i-1})$ are equal to $\mu(M)$. \end{enumerate}
 \end{lemma}

\proof Item (1) follows immediately form the additivity of $\deg$, and item (2) follows from the first by induction.\qed

 \begin{lemma}\label{L1.5} Let $\theta:\,\sC\to \sC'$ be a faithful exact\footnote{\ie which preserves short exact sequences, \cf \ref{kcok}.} functor between proto-abelian categories. Then any slope function  $\mu$ on $\sC'$ induces a slope function $\theta^\ast\mu$ on $\sC$.  
\end{lemma} 

\proof  $\theta^\ast\mu(M):= \mu(\theta(M))$ satisfies the two axioms of a slope function: (1) because $\theta$ preserves epi-monics by faithfulness, (2) because $\theta$ is exact. \qed

  \medskip\subsection{(Semi)stability.}\label{sst} Let $\sC, \rk, \La$ be as above, and let $\mu$ be a slope function on $\sC$ with values in $\La$.
 
    \begin{defn} A non-zero object $N$ of $\sC$ is called \emph{$\mu$-semistable} (\resp \emph{$\mu$-stable}) if for any non-zero  subobject $M \neq N$, $\mu(M)\leq \mu(N)$ (\resp $\mu(M)< \mu(N)$). 
  \end{defn} 
  
  If there is no ambiguity on $\mu$, one just says semistable (\resp stable).  
  
  \smallskip The next lemma deals with the behaviour of semistability with respect to monic and/or epi morphisms.  
  
    \begin{lemma}\label{PL2} Let $N$ be a non-zero object.  
   \begin{enumerate}
     \item $N$ is semistable if and only if for any non-zero {\rm strict} subobject $M$ of $N$, $\mu(M)\leq \mu(N)$.
      \smallskip  \item $N$ is semistable if and only if for any non-zero quotient (\resp {\rm strict} quotient) $P$ of $N$, $\mu(P)\geq \mu(N)$.
     
        \smallskip  \item If $N$ is semistable of slope $\la$, any non-zero subobject $M$ with $\mu(M)= \la$ (\resp quotient $P$ with $\mu(P)= \la$) is semistable of slope $\la$.
      
        \smallskip  \item  If $N$ is semistable of slope $\la$, any non-zero direct summand of $N$ is semistable of slope $\la$.
        
           \smallskip    \item Any non-zero strict subobject $M$ (\resp strict quotient $P$) of $N$ of minimal rank with $\mu(M)\geq \mu(N)$ (\resp $\mu(P)\leq \mu(N))$ is semistable. 
         
         In particular, any object of rank $1$ is semistable. 

  \smallskip  \item  If $N$ is semistable, any non-zero subobject (\resp quotient) with the same slope is semistable.
  
     \smallskip    \item Let $ 0\to M\to N\to P\to 0$ be a short exact sequence. If two of the three objects are semistable of slope $\la$, so is the third, unless it is zero. 
    
    In particular, in the additive case, a direct sum of semistable objects of slope $\la$ is semistable of slope $\la$. 
     \end{enumerate}
       \end{lemma} 
   
  \proof (1) Let $m: M\inj N$ be the given monic. Then $\bar m : M \to    \IM  \,m$  is monic-{epi}, hence $\mu(M)\leq \mu(  \IM  \,m)$, and the assertion follows. 
  
  \smallskip   \noindent   (2) Let  $e: N \to P$ be the given {epi}. Then  $\bar {e}:  \Coim\, {e}  \to P$  is monic-{epi}, hence $\mu( \Coim\, {e} )\leq \mu(P)$, and it suffices to show that $N$ is semistable if and only if $\mu(P)\geq \mu(N)$ for any non-zero strict quotient $P$.
  
Let us denote by $M$ the kernel (which we may assume to be non zero). Then lemma \ref{L1} shows that 
$\mu(P)\geq \mu(N)\Leftrightarrow \mu(M)\leq \mu(N)$ and $\mu(P)< \mu(N)\Leftrightarrow \mu(M)> \mu(N)$, from which the assertion follows (by item (1)).  
    
   \smallskip   \noindent  (3) Any subobject $M'$ of $M$ is a subobject of $N$, hence $\mu(M')\leq \la$; therefore $M$ is semistable of slope $\la$. Any quotient $P'$ of $P$ is a quotient of $N$, hence $\mu(P')\geq \la$; therefore $P$ is semistable of slope $\la$ by item (2). 

   \smallskip  \noindent   (4) and (6) follow immediately from the definition and (2).
   
    \smallskip  \noindent   (5) There is no proper strict subobject of $M$  (\resp proper strict quotient of $P$) of slope $\geq \mu(M)$ (\resp $\leq \mu(M)$) by minimality of the rank. The assertion then follows immediately from (1) and (2),  taking into account the fact the composition of strict monic (\resp epi) morphisms is strict.
 
    \smallskip  \noindent   (7) By lemma \ref{L1}, if two of the objects are of slope $\la$, so is the third. It follows that if $N$ and either $M$ or $P$ is semistable of slope $\la$, so is the third one, by (6). 
  
 Let us next assume that $M$ and $P$ are semistable of slope $\la$. If $N$ is not semistable, there is a semistable strict subobject $N'$ of $N$ with $\mu(N')>\la = \mu(N)$ (by items (1) and (5) above). By item (6), the induced morphism $N'\to P$ is zero, hence $N'\subset M$. Since $M$ is semistable, $\mu(N')\leq \mu(M)=\la$, a contradiction.  
   \qed
 
 \begin{lemma}\label{MSS} For any non-zero morphism $M\overset{f}{\to} N$ with $M$ and $N$ semistable, $\mu(M)\leq \mu(N)$.\end{lemma}
 
 \proof Let us consider the canonical factorization $M \to \Coim\, f  \inj N$. Then $\mu(M) \leq \mu(\Coim\, f  )\leq   \mu(N)$ by semistability of $M$ and $N$, taking into account item (2) of the previous lemma.
  \qed

 \begin{lemma}\label{SSQAB} \begin{enumerate} 
 \item
  Let $f: M\to N$ be a morphism between semistable objects of the same slope $\la$. Then $\Ker  f, \IM f, \Coker f, \Coim f$ are either zero or semistable of slope $\la$.
\item Let $\sC(\la)$ be the full subcategory of $\sC$ consisting of $0$ and the semistable objects of slope $\la$. Then $\sC(\la)$ is proto-abelian, and the notion of short exact sequence is compatible with the one in $\sC$. 
\end{enumerate} \end{lemma}
 
 \proof   Since $M$ and $N$ are semistable of slope $\la$, one has
  $\la \leq \mu(\Coim f) \leq \mu(\IM f) \leq \la$, whence equality of slopes, which implies $\mu(\Ker  f)= \mu(\Coker f)= \la$. The assertion then follows from items (3) and (7) of  \ref{PL2}. \qed

   \begin{rems}\label{remst}  {1)} Stable objects need not exist in general. If $M$ is simple (\ie irreducible), it is stable. The converse is not true, even if $\sC$ is abelian, \cf \ref{EXX}. See however \ref{SSS}. 
      
      \smallskip  \noindent  {2)} The analog of item (1) for stability (as opposed to semistability) is not true in general (\eg for objects of rank one). See however \ref{STABAB1} for a condition under which it holds. 
      
     \smallskip   \noindent {3)} See \ref{STABAB0} and \ref{STABAB1} for a condition under which $\sC(\la)$ is abelian.
    
    \smallskip \noindent {4)} Let $L\inj M\inj N$ be strict monic, with $L$ and $N$ semistable of the same slope $\la$. If $M$ is semistable, then $\mu(M)=\la$, and conversely. But it may happen that $\mu(M)<\la$. An example is given by the sequence of vector bundles 
  $$\sO_{{\mathbb P}^1}\overset{\iota_1}{\inj} \sO_{{\mathbb P}^1}\oplus  \sO_{{\mathbb P}^1}(-1)\overset{id\oplus  \iota}{\inj}  \sO_{{\mathbb P}^1}\oplus  \sO_{{\mathbb P}^1}\oplus  \sO_{{\mathbb P}^1} $$ on the projective line, with respect to the Harder-Narasimhan filtration.  
  
    \smallskip \noindent {5)} $\mu$-(semi)stable objects in $\sC$ are the same as $(-\mu)$-(semi)stable objects in $\sC^{op}$. \end{rems}

   \medskip\subsection{Universal destabilizing subobject.}\label{uds} Let $N$ be a non-zero object of $\sC$. 
   
   \begin{defn} A \emph{universal destabilizing subobject} of $N$ (with respect to $\mu$) is a non-zero strict subobject $M\inj N$ such that for any non-zero strict subobject $M'$ of $N$, 
  
  \noindent $i)$ $\mu(M')\leq \mu(M)$,
  
  \noindent $ii)$ if  $\mu(M')= \mu(M)$, then $M'\inj N$  factors through $M$. 
    \end{defn}
   
   In order to check these conditions, one may assume that $M'$ is semistable in virtue of  \ref{PL2} (3). 
   
\begin{lemma}\label{UDO} A universal destabilizing subobject exists and is unique. Moreover, it is semistable.
\end{lemma}
 
\proof Uniqueness follows from universality. Semistability follows from condition $i)$ (and \ref{PL2} (1)).

   Let us prove existence by induction on $\rk N$. 
  
  If $N$ is semistable, $M=N$ works. Otherwise, let $P$ be a strict quotient of $N$ of minimal rank with $\mu(P)\leq \mu(N)$, and set $ N' := \Ker  (N\to P)$. By \ref{PL2} (5), $P$ is semistable,  and by  \ref{L1},   one has $$ \mu(P)\leq  \mu(N)\leq \mu(N'),  \; \rk N'< \rk N.$$ Let $M$ be the universal destabilizing subobject for $N'$; in particular, $M$ is semistable of slope $\geq \mu(N)$. In case of equality, $N'$ would be semistable of slope $\mu(N)$, so would be $P$ by  \ref{L1} and $N$ by \ref{PL2} (7) contrary to assumption. Therefore $\mu(M) >\mu(N)$.
  
   Let $M'$ be a semistable strict subobject of $N$. 
   
   If the composed morphism $M'\to P$ is non-zero, we have $\mu(M')\leq \mu(P)$ by lemma \ref{MSS}, hence $\mu(M')\leq \mu(N) <\mu(M)$. 
   
   Otherwise $M'$ is a strict subobject of $N'$. Therefore $\mu(M')\leq \mu(M)$, with equality only if $M'$ is a subobject of $M$. In both cases, this shows that $M $ is the universal destabilizing subobject for $N$. \qed

 \bigskip
 \section{Slope filtrations and Newton polygons.}\label{sfnp} In this section, we introduce the concept of a slope filtration (a functorial filtration of objects of $\sC$ by strict subobjects, satisfying some conditions). 

We establish a one-to-one correspondence between slope filtrations and slope functions $\mu$, which synthesizes a lot of (more or less ad hoc) constructions of slope filtrations in the literature.

We then discuss Newton polygons, and examine in some detail the exactness properties of slope filtrations. 
 
\medskip
  \subsection{Filtrations by strict subobjects.}\label{FSS} As usual, we may consider the (totally) ordered set $\La$  as a small category.   

\medskip A \emph{decreasing functorial filtration on $\sC$ by strict subobjects, indexed by $\La$,} is a functor $$F^{\geq .}(.) :\;\La^{op}\times \sC \to \sC $$ which sends any object $(\la, M)$  to a strict subobject $F^{\geq {\la}}M$ of $M$. 

It is \emph{separated} (\resp \emph{exhaustive})  if for any $M$,
$$  \varprojlim \, F^{\geq {\la}}M = 0,\;\;    \varinjlim \, F^{\geq {\la}}M = M. $$ 
 
 It is \emph{left continuous} if for any object
$$F^{\geq {\la }}M = \varprojlim_{\la'<\la}\, F^{\geq {\la'}}M. $$ 

\smallskip Using the fact that the ranks bounds the length of any flag, it is easy to see that for any separated, exhaustive, left continuous filtration, and any object $M$, there is a partition of $\La$ by intervals 
$$I_r = ]-\infty, \la_r ],\ldots, \; I_2=  ]\la_2, \la_1 ],  I_1=  ]\la_1, +\infty [$$ such that $F^{\geq\la} M$ is {\it constant} on each of these intervals, and a flag\footnote{sometimes called the {\it Harder-Narasimhan flag} of $M$, with example \ref{exs1} in mind. } of length $r$ 
$$\sF(M):\;0\inj M_1= F^{\geq \la_1}M \inj \cdots \inj  M_r = F^{\geq \la_r}M = M.$$

\begin{defn}\label{break} The elements $\la_1 > \la_2 > \ldots > \la_r$ are called the  {\it breaks}\footnote{we refrain from calling them the slopes of $M$, in order to prevent confusion with $\mu(M)$.} of $M$ (with respect to the filtration $F^{\geq .} $).   \end{defn}

\medskip From $F^{\geq .}$, one also gets another {descending functorial filtration $F^{> .}$ by strict subobjects, indexed by $\La$,}  by setting 
$$ F^{> \la}M := \varinjlim_{\la' >\la}\, F^{\geq {\la'}}M.$$ The colimit exits:  indeed, $F^{> \la}M = F^{\geq {\la'}}M$ for $\la'>\la$ close enough to $\la$, hence is a strict subobject of $F^{\geq \la}M$, so that there is a short exact sequence in $\sC$:
$$0\to  F^{> \la}M \to F^{\geq \la}M \to {\gr}^\la M\to 0.$$ 
We set ${\gr}\, M := \oplus \,  {\gr}^\la\, M = {\gr}^{\la_1}\, M\oplus \ldots \oplus {\gr}^{\la_r}\,M.$ This is functorial in $M$.

\begin{defn}\label{mul} The {\it multiplicity} of the break $\la$ (in $M$) is the rank of $\gr^\la M$.
 \end{defn}

 \begin{rem} 1) The data of the filtration $F^{\geq .}$ is equivalent to the data of the filtration $F^{>.}$, thanks to the formula
 $$F^{\geq {\la }}M = \varprojlim_{\la'>\la}\, F^{> {\la'}}M. $$ 
The filtration $F^{> \la}M$ is {\it right continuous}:
$ F^{> {\la }}M = \varinjlim_{\la'>\la}\, F^{> {\la'}}M. $ 
   \end{rem}

      \begin{rem}\label{asc} In the literature, one also encounters {\it increasing filtrations by strict subobjects} $F_{\leq .}$. They are defined in the same way as decreasing filtrations by strict objects, except that $\La^{op}$ is replaced by $\La$, and left continuity by right continuity. 

In practice, the distinction between descending and ascending filtrations is not essential: one passes from one to the other by changing $\la$ into $-\la$, more precisely, by setting $$F_{\leq {\la }}M = F^{\geq {-\la}}M,$$ and by reversing all inequalities in the definition of slope functions and (semi)stability.  It is therefore just a matter of convention on the signs of slopes. For uniformity, we concentrate on descending filtrations in the sequel. 
 \end{rem}

  \medskip\subsection{Slope filtrations.}\label{slf}    
Given a descending filtration by strict subobjects as before, one defines functions
 $$\deg: \sk\, \sC \to \La, \;\;  \deg M =   \sum_\la \, \la \cdot \rk \gr^\la M,$$
 and 
 $$\mu: \sk\, \sC \setminus\{0\} \to \La,\;\; \mu(M)= \frac{\deg M}{\rk M}.$$
 
   \begin{defn} A \emph{descending slope filtration $F^{\geq .}$ on $\sC$} (indexed by $\La$) is a separated, exhaustive, left continuous decreasing functorial filtration on $\sC$ by strict subobjects,  that satisfies  
\begin{enumerate}
{\it \item for any $\la$,  the filtration of $F^{\geq \la}N$ (\resp $N/F^{\geq \la}N$) is induced by the  filtration of $N$,

\smallskip  \item  the associated function $\mu $  is a slope function} in the sense of \ref{sf}. 
 \end{enumerate}
 \end{defn}

\noindent The trivial slope filtration is the one attached to the $0$ slope function: $F^0M=M, \, F^{>0}M=0.$
 
    \begin{prop}\label{P3}  Let $F^{\geq .}$ be a slope filtration on $\sC$ and let $N$ be a non-zero object of $\sC$. \begin{enumerate}
  \item
 The flag $\sF(N)$ attached to $F^{\geq.}N$ is the unique flag on $N$ (up to unique isomorphism) whose graded pieces are semistable of slopes arranged in decreasing order: 
 $$\mu(N_1)>\mu(N_2)>\ldots >\mu(N_r).$$  
  
   In particular, $N$ is semistable if and only if it has a unique break (which is then $\mu(N)$). 

\smallskip \item In fact, $N_1$ is the universal destabilizing subobject of $N$. More generally, $N_i$ is the pull-back by $N\to N/N_{i-1}$ of the universal destabilizing subobject of $N/N_{i-1}$.

\smallskip \item  Under axiom $(2)$ of slope filtrations, axiom $(1)$ is equivalent to:

 \smallskip $(1)' $ for any $\la $, $\gr^\la \circ \gr^\la = \gr^\la$.
\end{enumerate} \end{prop} 
 
 \proof (1) and (2). Functoriality of $F^{\geq .}$ implies that $N$ is semistable if $N= \gr^{\la}\,N$. 
 
 Assume either axiom (1) or (1)'. Then, by the inequality of \ref{L1}, the converse holds: if $N$ is semistable, $N= \gr^{\la}\,N$. It follows that the graded pieces of $\sF(N)$ are semistable of slopes equal to the breaks of $N$, in decreasing order (taking into account the condition $\gr \circ \gr = \gr$). 
  
 Let $$\sF':  0= N'_0 \inj  N'_1 \inj   \cdots \inj N'_s = N$$ be a flag on $N$ with $N'_i/N'_{i+1}$ semistable of slope $\la'_i$ and $\la'_i>\la'_{i+1}$. 
 
 We prove at the same time equality $\sF'= \sF$ and assertion (2) by showing that $N'_i$ is the pull-back by$N'\to N'/N'_{i-1}$ of the universal destabilizing subobject of $N'/N'_{i-1}$. By induction on the rank, it is enough to deal with $i=1$. 
 
Let $M$ be the universal destabilizing subobject of $N$. One has $\la'_1\leq \mu(M)$, and equality only if $N'_1 \subset M$. 

  Let $j\geq 1$ be the first index for which $M\inj  N$ factors through $N'_{j }$. The composition $N_1 \to  N'_j/N'_{j-1}$ is a non-zero morphism between semistable objects, hence $\mu(M)\leq \la'_j$ by lemma \ref{MSS}. 
  
  One concludes that $j=1,\,\mu(M)= \la'_1 $, and $M= N'_1$.

\smallskip \noindent (3)  Assume (1)'.  Then the graded pieces of the subflag $\sF(N)\cap F^{\geq \la}N$  (\resp quotient flag $\sF(N)/ F^{\geq \la}N$) of $\sF(N)$ are semistable of slopes arranged in decreasing order, hence $\sF(N)\cap F^{\geq \la}N = \sF(F^{\geq \la}N)$  (\resp   $\sF(N)/ F^{\geq \la}N= \sF(N / F^{\geq \la}N)$) by item (1) of the proposition. 

Conversely, let us assume axiom (1). Then the filtration of the subquotient $\gr^\la N$ of $N$ is induced by the filtration of $N$, hence $ \gr^\la N= F^{\geq \la} \gr^\la N= \gr^\la \gr^\la N$.   \qed

   \bigskip
   The following theorem digests most of the avatars found in the literature of existence theorems for filtrations of Harder-Narasimhan type and their functoriality.

   \medskip
  \begin{thm}\label{T1} The rule $\,F^{\geq .} \mapsto \mu\,$ induces a \emph{bijection} between slope filtrations (up to unique isomorphism) and slope functions on $\sC$.
  \end{thm}

\smallskip  \proof Let us fix a slope function $\mu$. By item
   (2) of the previous proposition, we know the right candidate for the first notch $N_1$ of $\sF(N)$ (the notch of maximal slope): it is the universal destabilizing subobject $N_1$ (which depends only on $N$ and $\mu$). 
    
  \medskip Let us show the existence of a flag $\sF(N)$ with the property that $N_i/ N_{i-1}$ is semistable of slope $\la_i := \mu(N_i/ N_{i-1})$, with $\la_{i} > \la_{i+1}  $ ($i= 1,\ldots, r$). We proceed by induction on $\rk N$. We consider such a flag $\sF(N/N_1)$ for $N/N_1$. The pull-back of $\sF(N/N_1)$ by $N\to N/N_1$  is a flag $\sF(N)$ (lemma \ref{PB}) and the corresponding morphisms $N_i \to (N/N_1)_i$ are strict {epi}; moreover $N_i/N_{i-i}\cong  (N/N_1)_i/ (N/N_1)_{i-1}$ is semistable for $i\geq 1$. It is then clear that $\sF(N)$ has the desired properties. 
  
 \smallskip Let us set 
 
 \centerline { $ F^{\geq\la} N = N$ if $\la\leq \la_r , \; F^{\geq\la} N = N_i$ if $\la \in   ]\la_{i+1 }, \la_{i  }  ],\;\; F^{\geq\la} N = 0$ if $\la > \la_1$.  } 
      
  \smallskip It is clear that this is a separated, exhaustive, left continuous decreasing filtration on $\sC$ by strict subobjects, indexed by $\La$, which satisfies $\gr^\la \circ \gr^\la = \gr^\la$. Moreover, the associated degree (\resp slope) function is the original one. 
   
   Let us finally prove functoriality, \ie that any morphism $f: M\to N$ sends $F^{\geq\la} (M)$ to $F^{\geq\la} (N)$. By descending induction, we may assume that  $F^{>\la} (M)\to N$ factors through $ F^{>\la} (N)$, and we have to prove that $\gr^\la \,M \to N/ F^{>\la} (N) $ factors through $\gr^\la \,N $. Since $\gr^\la \,M $ is semistable of slope $\la$, its image $P$ in $N/ F^{>\la} (N)$ has $\mu(P)\geq \la$, hence is contained in $\gr^\la \,N$ by construction of the filtration.

  \smallskip This proves the surjectivity of $F^{\geq .} \mapsto \mu$.  
  
  Injectivity follows from item (1) of the previous proposition.  \qed

  \begin{cor}\label{CO1} Let $\sC$ be a full subcategory of a proto-abelian category $\sC'$, such that any strict subquotient in $\sC'$ of an object of $\sC$ is an object of $\sC$ (so that $\sC$ is proto-abelian, and the embedding $\sC\inj \sC'$ is exact). Let $\rk$ be a rank function on $\sC'$.

Let $\mu$ be a slope function on $\sC'$, and $\mu_{\mid \sC}$ be the slope function on $\sC$ induced by $\mu$ (\cf \ref{L1.5}). Then the slope filtration attached to $\mu_{\mid \sC}$ is the restriction to $\sC$ of the slope filtration attached to $\mu$.
\end{cor}    
 
 Indeed, these two slope filtrations on $\sC$ have slope function $\mu$, hence coincide. \qed

\begin{rems}\label{remsslo} 1) Up to now, the additivity of $\deg$ on short exact sequences $0\to M\to N\to P\to 0$  has been used only via the inequalities of lemma \ref{L1} (item (1)). 
 One could thus weaken the definition of slope functions and slope filtrations, retaining these inequalities instead of $\deg N= \deg M +\deg P$. 
 
 In \cite{Rud}, a formalism of stability is developped where these inequalities are taken as an axiom, but only in the context of abelian categories. 
 
  \smallskip \noindent 2) The slope filtration $\check F^{\geq .}$ on $\sC^{op}$ corresponding to the slope function $-\mu$ is given by
  $$ \check F^{\geq \la} M = M / F^{>-\la} M.$$
 \end{rems}

\medskip
Let $\sC$ be quasi-abelian, with left abelian envelope $\sA$.  Recall that 
 $K_0(\sC)=K_0(\sA)$, which contains the subgroup generated by the torsion classes $[Q], \, Q\in Ob \, \sT$.
 
  \begin{cor}\label{STABAB0} If $\sC$ is quasi-abelian, the rule $\,F^{\geq .} \mapsto \deg\,$ induces a \emph{bijection} between slope filtrations on $\sC$ (up to unique isomorphism) and homomorphisms
  $$K_0(\sA)\to \La$$ that are \emph{non-negative on torsion classes $[Q], \, Q\in Ob \, \sT$}.   
  
  \smallskip If moreover $\deg$ is positive on non-zero torsion classes, and if all torsion classes have rank $0$, then the full subcategory $\sC(\la)$ of $\sC$ consisting of $0$ and the semistable objects of slope $\la$ is {\rm abelian} (hence artinian and noetherian by \ref{FINITELENGTH}). \end{cor}
 
 \proof For the first assertion, it only remains to see that an additive map $\deg: K_0(\sA)\to \La$ is a degree function, \ie satisfies 

\smallskip \centerline{ $\forall M\to N$ epi-monic in $\sC$, $\;\deg M / \rk M \leq \deg N/\rk N$} if and only if $$\forall Q\in Ob\, \sT, \,\deg Q \geq 0.$$
 This follows immediately from lemma \ref{TORS}. 
 
 For the second assertion, one has to see that for any epi-monic $M\overset{f}{\to} N$ in $\sC$ with $M$ and $N$ semistable of the same degree is an isomorphism. By \ref{TORS} again, one has a short exact sequence in $\sA$
 $$0\to M\overset{f}{\to} N\to \Coker f \to 0$$ with $\Coker f\in \sT$, and $\deg \Coker f = \deg N- \deg M=0$. Hence $\Coker f=0$ by assumption. \qed

   \begin{exs}\label{Ff} 1) In the case of example \ref{exs1} (torsionfree coherent sheaves on a polarized normal connected variety $(X, \sO(1))$), the filtration attached to $\mu $ (in either the Mumford-Takemoto or the Gieseker-Maruyama version) is the {\it Harder-Nararasimhan filtration} which is generally used in the study of moduli spaces (at least when $X$ is smooth), which extends, as is well-known, to the abelian category $\sA$ of coherent sheaves. The assumptions of \ref{STABAB0} are satisfied.
      
\smallskip   Note that, in the construction of the slope filtration, we have not used the fact that the slope function on subsheaves of a given torsionfree coherent sheaf is bounded from above. Rather, this fact appears as an immediate corollary of the construction.
   
   \smallskip \noindent   2) Let $\sC$ be the quasi-abelian category of finite flat commutative group schemes over a $p$-adic field $(K, v)$, of $p$-primary order. The height $\, {\rm ht}\,$ provides a rank function on $\sC$. In \cite{Fa}, L. Fargues considers the following degree function $\,\deg$: if the conormal sheaf $\omega_G$ decomposes as $\oplus \sO_K/a_i\sO_K$,  $\deg   G = \sum v(a_i)$. He shows that $\mu = \frac{\deg}{\rm ht}$ satisfies the axioms of a slope function (and takes values in $[0,1]$), and studies the associated slope filtration on $\sC$.  
    Moreover, by \cite[prop. 2]{Fa} and \ref{TORS}, the assumptions of \ref{STABAB0} are satisfied.
    \end{exs}

 \begin{cor}\label{OBABATS} If $\sC$ is abelian, any slope filtration on its socle $\sC_{ssi}$ (\ie the full subcategory of $\sC$ consisting of semisimple objects) comes from a unique slope filtration on $\sC$.  \end{cor}
 
 \proof The degree function on $\sC_{ssi}$ to a degree function on $\sC$ since $K_0(\sC)= K_0(\sC_{ssi})$. The corresponding slope filtration on $\sC$ then extends that on $\sC_{ssi}$ by \ref{CO1}. \qed

 \bigskip
\subsection{Highest break function.}\label{hbf}  The highest break of the slope filtration defines a function  $$\rho: \sk\, \sC \setminus\{0\} \to \La .$$ Of course, $\mu\leq \rho$.

\begin{rem} In the case of the Turrittin-Levelt filtration, the highest break is called the {\it Poincar\'e-Katz rank} and can be interpreted as a spectral radius; it is commonly denoted by $\rho$, whence the choice of this symbol; another common notation is $\mu_{max}$. 
 \end{rem}

 \begin{prop} \begin{enumerate} 
 \item An object $N$ is semi-stable of slope $\la$ if and only if for any non-zero strict quotient $P$, $\rho(N)\leq \rho(P)$.
 \item A slope filtration is determined by its highest break function.
 \end{enumerate} \end{prop}
 
 \proof  (1)  If $N$ is semistable, one has $  \rho(N)= \mu(N) \leq \mu(P) \leq \rho(P)$ by item 2 of  \ref{PL2}. Conversely, let $\nu$ be the lowest break of $P$. Then $Q= \gr^\nu P$ is a strict quotient of $N$ which is semistable of slope $\rho(Q) =\nu$. By assumption, $\rho(N)\leq \rho(Q)$, hence $\mu(N)  \leq  \nu \leq \mu(P) $.

\smallskip\noindent (2) By item (1), for two slope filtrations with the same highest break function, an object $M$ is semi-stable of slope $\la$ for one filtration if and only if it is so for the other filtration. By the characterization \ref{P3} of the canonical flags, the filtrations coincide.  \qed

 \bigskip
\subsection{Newton polygons.}\label{newp} To fix ideas, we assume in this subsection that $$\La \subset \R.$$
   Let $F^{\geq .}$ be a descending slope filtration with values in $\La$, and let $\deg$ and $\mu$ be as before the associated degree and slope functions. 
  
To any object $M$, we attach its Newton polygon in $\R^2$ whose slopes are the breaks of $F^{\geq .}M$. Since we are dealing with a descending filtration, the natural convention is to arrange the slopes in decreasing order (from left to right), thus giving rise to {\it concave} piecewise affine functions\footnote{this is the usual convention in the context of stability for vector bundles.}.
More precisely, let us introduce the following definition.

\begin{defn} \begin{enumerate}\label{NP}
\item The \emph{polygon\footnote{{\it stricto sensu}, this is not a polygon, since this convex set is unbounded from below; but the terminology is traditional.} $P(\sF)$ of a flag}   $$  \sF: \, 0= M_0 \inj  M_1 \inj   \cdots \inj M_r = M $$ is the convex hull of the points  with coordinates $$\,(x= \rk M_i, \; y \leq \deg M_i).$$
 \item The \emph{Newton polygon} of  $M$ is the polygon of the flag $\sF(M)$ attached to  $F^{\geq .}M$ 
 $$ NP(M) := P(\sF(M)).$$
 \end{enumerate}
\end{defn} 

\begin{lemma} The end-points of $ NP(M)$ are $(0,0)$ and $(\rk M, \, \deg M)$. The slope of the segment linking these points is the slope $\mu(M)$ of $M$. 

The points $\,(x= \rk M_i, \; y \leq \deg M_i) $ are extremal points of $ NP(M)$. The slopes of the edges of $ NP(M)$ are the breaks $\la$ of $M$, and the horizontal length of such an edge is $\rk \, \gr^\la M$ .
 \end{lemma}

This is immediate.\qed

 \begin{lemma} If $\sC$ is additive, and $N= M\oplus P$, then the breaks of $N$ are the breaks of $M$ and of $N$, counted with multiplicities. A fortiori $NP(N)= NP(M)+NP(P)$ (in the sense of the Minkovsky sum of convex sets).
\end{lemma}

\proof Indeed, for any $\la$, the additive functor $\gr^\la$ preserves $\oplus$, hence $\gr^\la N = \gr^\la M \oplus \gr^\la P$. \qed

 \begin{prop}\label{P2} The polygon of any flag $\sF $ (of any length) ending with $M$  lies below $NP(M)$, with the same end-points. \end{prop} 

 \proof Let $\,\sF' \,$ and $\,\sF''\,$ be other flags on $M$. If $\sF''$ is a refinement of $\sF'$, then $P(\sF')$ lies below $P(\sF'')$. On the other hand, for any refinement $\sF''$ of $\sF(M)$, $P( \sF'') = NP(M)$. Indeed, if $0\subset M_{i-1} \subset N\subset M_i$ is a flag, then the point $(\rk N, \deg N)$ lies below the segment joining the points $(\rk M_{i-1}, \deg M_{i-1})$ and $(\rk M_{i }, \deg M_{i })$ since $M_i/M_{i-1}$ is semistable. 
 
One concludes by using  a common refinement of $\sF(M)$ and $\sF'$, \cf \ref{REFI}.  \qed
 
   \medskip
  
In the context of vector bundles, this characterization of $NP(M)$ was given by Shatz \cite{Sh}.
 
 \begin{rem}\label{asc2} If one deals with ascending slopes filtrations, it is then natural to define $P(\sF)$  as the convex hull of the points  with coordinates $\,(x= \rk M_i, \; y \geq \deg M_i)$, which gives rise to a {\it convex} piecewise affine function\footnote{this is the usual convention in the context of $p$-adic Frobenius slopes.}. The end-points are again $(0,0)$ and  $(\rk M, \, \deg M)$.
Passing to the associated descending filtration ($F^{\geq \la}=F_{\leq -\la}$) results in changing the polygon of a flag on $M$ by a symmetry with respect to the horizontal axis (and changing the sign of $\deg M$).  

\smallskip On the other hand, if one insists on dealing with {\it convex} piecewise affine functions  in the presence of a descending slope filtration\footnote{this is the usual convention in the context of ramification theory and asymptotic analysis of differential equations.}, one may consider the  sequence of strict {epi}s\footnote{that is nothing but the flag on $M$ with respect to the dual slope filtration $\check F^{\geq .}$, \cf \ref{remsslo} 2).} 
$$  M'_r = M \surj  M'_{r-1} = M /M_{1} \surj \cdots \surj M'_1= M /M_{r-1} \surj M'_0= 0.$$
associated to the flag $\sF$ and redefine the polygon of $\sF$ to be the convex hull of the  points with coordinates $(x= \rk M'_i, y \geq \deg M'_i)$. The end-points are $(0, -\deg M)$ and $(\rk M, 0)$.
The relation with the polygon defined in \ref{NP} is a symmetry through the point $(\frac{\rk M}{2}, 0)$.  \end{rem}

   \bigskip\subsection{The topological space of all slope filtrations.} Let us endow $\La$ with the canonical topology generated by the open intervals. 
   
   By theorem \ref{T1}, slope filtrations are in bijection with degree functions 
   $$\deg: K_0(\sC)\to \La.$$ Endowing $Hom(K_0(\sC), \La)$ with its natural linear (weak) topology, the space of degree functions (which is defined by the linear inequalities $$\deg ([N]-[M])\geq 0$$ if there is a monic-{epi} from $M$ to $N$) is a convex cone\footnote{\ie is stable under linear combinations with non-negative coefficients.}, whose apex corresponds to the trivial slope filtration. It is in fact a {\it closed convex cone}. 
      
   Given a non-zero object $N$, the condition that $N$ is $\mu$-semistable (\ie the set of linear inequalities 
   $$\rk M . \deg [N] - \rk N. \deg [M] \geq 0$$  if there is a monic from $M$ to $N$) defines a {\it closed convex subcone}.   \medskip

 \begin{ex}\label{EHN} Let us compute the space of all slope filtrations on the quasi-abelian category of vector bundles on a smooth connected projective curve $X$ of genus $g$. 
 
 Since $\rk $ is additive in short exact sequences and $\rk M\leq \rk N$ whenever there is a monic $M\inj N$, any constant function on $Sk\, \sC\setminus \{0\}$ is a slope function. 
 By addition of a constant, we may consider only slope functions $\mu$ with $\mu(\sO_X)=0$.
  Among them, there is the canonical (Mumford) slope function $\mu_{can}$.
  
 One has a group isomorphism 
 $$K_0(\sC)\cong K_0(X)\overset{(\det, \rk)}{\to} Pic(X)\oplus \Z,$$
 and the obvious mapping $Pic(X) \to K_0(X)$ is a set-theoretic section of the projection $K_0(\sC) \overset{\det }{\to}Pic(X)$. Therefore $ \mu\cdot \rk $ factors through this projection, and is determined by its value on $Pic(X) \subset K_0(X)$. Moreover, for any $L, L'\in Pic(X)$, one has, $\mu(L\otimes L')= \mu(\det (L\oplus L'))= (\mu\cdot \rk)(L\oplus L') = \mu(L)+\mu(L')$.

 On the other hand, there is an exact sequence
 $$ 0 \to Pic^0(X) \to Pic(X) \overset{\mu_{can}}{\to} \Z.$$   
 Let us show that $\mu$ vanishes on $Pic^0(X)$. Indeed, for any $L, L'\in Pic(X)$ of degree $0$ and $d\geq g$ respectively, and for any integer $n$, there is, by Riemann-Roch, a  monic $L^{\otimes n} \inj L'$. This implies $\mu(L^{\otimes n})= n\mu(L) \leq \mu(L')$, whence $\mu(L)=0, \mu(L')\geq 0$. Therefore $\mu$ factories through a non-negative multiple of $\mu_{can}$ on $Pic(X)$.
 
 In conclusion, {\it any slope function on $\sC$ is of the form 
 $$ \lambda\cdot  \mu_{can}+ \lambda'  ,\; \lambda \in \La_{\geq 0}, \lambda'\in \La .$$}
   \end{ex}

        \bigskip\subsection{Rees deformation from ${\rm gr}\, M$ to $M$.}\label{Rees} Assume that  $\sC$ consists of modules (or sheaves of modules) over some domain $R$, with some extra structure. For a given object $M$, let $\La_M\subset \La$ be a finitely generated sub-semigroup such that the associated group is free and equal to the subgroup of $\La$ generated by the breaks of $M$ (for instance, if $\La=\Q$, one may choose ${\La_M}=   \frac{1}{d}\N$ to be the semigroup generated by the inverse of the common denominator $d$ of the breaks). Without loss of generality, one may assume that $\Z[[x^{\La_M}]]$ is a regular algebra  and that the least non-zero element of ${\La_M}$ is less or equal to the positive differences between breaks. 
  
 Then one can form the following variant of the Rees module over $R[[x^{\La_M}]]$:
  $$R(M) = R_{\La_M}(M)=  \sum_\la \; F^{\geq\la}M \cdot x^{-\la}R[[x^{\La_M}]]$$ (as a submodule of $x^{-\nu}M\otimes_R R[[x^{\La_M}]]$). The generic fiber is isomorphic to $M$, whereas the special fiber is isomorphic to ${\rm gr}\, M$ (variant: one could work with $R[x^{\La_M}]$ instead of $R[[x^{\La_M}]]$).
 
 This construction is functorial: any $f\in \sC(M, N)$ gives rise, for suitable $\La_{M,N}$, to a morphism $R(M)\to R(N)$ whose special fiber is $\gr f$.

       \bigskip
 \section{Exactness properties.}

\subsection{Exact filtrations.}\label{exa}   
  Let $F^{\geq .}$ be a separated, exhaustive, left continuous decreasing filtration by strict subobjects on the proto-abelian category $\sC$, as in subsection \ref{FSS}. 
 
 \begin{defn} A morphism $f: M\to N$ is \emph{strictly compatible} with $F^{\geq .}$ if for any $\la$, the canonical (strict) monic 
 $$ f(F^{\geq \la} M)\inj f(M) \cap F^{\geq \la} N $$ is an isomorphism.
 \end{defn}
 
 This is equivalent to saying that in the canonical factorization of $f = m\circ e$ as a strict {epi} $e: M \surj \Coim f$ followed by a monic $m$, both $e$ and $m$ are strictly compatible with $F^{\geq .}\,$ (the composed monic 
 $$ f(F^{\geq \la} M)\inj m(F^{\geq \la} eM) \inj m(eM) \cap F^{\geq \la} N = f(M)\cap  F^{\geq \la} N$$ being an isomorphism if anf only if so are the two intermediate monics).
 
 \smallskip Caution: the composition of two morphisms that are strictly compatible  with $F^{\geq .}$ is not necessarily strictly compatible with $F^{\geq .}$.

 \begin{defn}  $F^{\geq .}$ is \emph{exact} (\resp \emph{strongly exact}) if any \emph{strict} morphism\footnote{\ie a composition of a strict epi followed by a strict monic.  Note, on the other hand, that in the additive case, {\it any} morphism can be factored, in the opposite order, as $e\circ m$ where $m$ is a strict monic and $e$ a strict {epi}; take $m=(id, f): M \to M\oplus N$ and $e= pr_2$ and note that there is a short exact sequence  $ 0\to   M \overset{ (id, f)}{ \to} M\oplus N \overset{f\circ pr_1 - pr_2}{\to} N \to 0   $.}  (\resp  \emph{any} morphism) is strictly compatible with $F^{\geq .}$. 
   \end{defn}
 
  If $\sC$ is abelian, there is no difference between these two notions, of course. 
  
  \begin{exs}\label{exex}  1) Among our five basic examples, it turns out that the Turritin-Levelt, Hasse-Arf and Dieudonn\'e-Manin slope filtrations are exact, but  the Harder-Narasimhan slope filtration is not, as the consideration of the standard short exact sequence of vector bundles
 \begin{equation}\label{exseqHN} 0\to \sO_{\bP^1}(-1)\to \sO_{\bP^1}^2\to \sO_{\bP^1}(1)\to 0\end{equation} shows. The Grayson-Stuhler filtration is also non-exact, as the as the consideration of the standard short exact sequence  of euclidean lattices
  \begin{equation}\label{exseqGS} 0\to (1,1)\cdot\Z \to \Z^2 \to (\frac{1}{2},-\frac{1}{2})\cdot \Z\to 0 \end{equation} shows.
  
  \smallskip \noindent 2) On the proto-abelian category of vector bundles of rank $\leq 1$ on a smooth connected projective curve, the standard slope filtration is exact but not strongly exact. However, we do not know any example of an exact, but not strongly exact, slope filtration on a {quasi-abelian} category.  \end{exs} 
 
    \medskip
 
In the sequel, we assume that $\sC$ is quasi-abelian.

 \smallskip 
 \begin{defn}(\cite[1.1.18]{Sc}). A functor $\theta$
 between quasi-abelian categories is \emph{exact} (\resp \emph{strongly exact}) if it preserves short exact sequences  (\resp if it preserves kernels and cokernels).
 \end{defn}
 
 It follows from \cite[1.1.15, 1.1.16]{Sc} that $\theta$ is strongly exact if and only if it is exact and preserves epi-monics. This is the characterization that we shall use.
 
 \medskip

\begin{lemma}\label{eqex}  The following properties are equivalent: 
  \begin{enumerate} 
  \item  $F^{\geq .}$ is exact (\resp strongly exact),
 \smallskip\item  for every $\la$, $F^{\geq \la}$ is an exact functor (\resp strongly exact functor),
 \smallskip \item for every $\la$, $\gr^\la\,$ is an exact functor (\resp strongly exact functor),
 \smallskip \item the ``dual filtration" $\check F^{\geq .}$ on $\sC^{op}$ (given by $\check F^{\geq \la}N= N/F^{>-\la}N$) is exact (\resp strongly exact).
 \end{enumerate}
  \end{lemma} 
 
  \proof  $(1)\Leftrightarrow (2)$ is straightforward: the strict compatibility of any strict epi and any strict monic (\resp and also any epi-monic) with $F^{\geq .}$ implies that for every $\la$, $F^{\geq \la}$ is an exact functor (\resp and preserves epi-monics), and conversely.
  
     \smallskip  $(2)\Leftrightarrow (3)$:  let $0\to M \to N \to P\to 0$ be a short exact sequence. Let us consider the following commutative diagram in $ \sC$ with exact columns: 
    \[\begin{CD} 0 @>>> 0 @>>>  0 @>>>  0 @>>>  0\\ @VVV @VVV @VVV @VVV @VVV \\  0 @>>> F^{> \la}M @>>>  F^{> \la}N @>>>  F^{> \la}P @>>>  0 \\ @VVV @VVV @VVV @VVV @VVV \\  0 @>>> F^{\geq \la}M @>>>  F^{\geq \la}N @>>>  F^{\geq \la}P @>>>  0\\ @VVV @VVV @VVV @VVV @VVV \\  0 @>>> {\rm gr}^\la M @>>>  {\rm gr}^\la N @>>>  {\rm gr}^\la P @>>>  0\\ @VVV @VVV @VVV @VVV @VVV \\  0 @>>> 0 @>>>  0 @>>>  0 @>>>  0. \end{CD}\]  
If $F^{\geq .}$ is exact, the second and third rows are exact. By the snake lemma in the left abelian envelope ${\sA}$, \cf \ref{P1}, it follows that the fourth row is also exact in ${\sA}$, hence in $\sC$. This shows that $gr$ is exact. 

 For the converse, we argue by descending induction on $\la$: we assume that the second row of the diagram is exact. If $gr$ is exact, the fourth line is exact, and it follows that the third is also exact.
 
 \smallskip It remains to prove that $F^{\geq \la}$ preserves epi-monics for any $\la$ if and only if so does $\gr^\la$ for any $\la$.  Let us note that $F^{\geq \la}$ and $F^{> \la}$ always preserve monics (independently of exactness).

Let $M\overset{f}{\to} N$ be epi-monic in $\sC$ and let us consider the following commutative diagram in $ \sA$ with exact columns: 

    \[\begin{CD} 0 @>>> 0 @>>>  0 @>>>  0 @>>>  0\\ @VVV @VVV @VVV @VVV @VVV \\  0 @>>> F^{> \la}M @>>>  F^{> \la}N @>>>  T = F^{> \la}N/F^{> \la}M @>>>  0 \\ @VVV @VVV @VVV @VVV @VVV \\  0 @>>> F^{\geq \la}M @>>>  F^{\geq \la}N @>>>  T' = F^{\geq \la}N /F^{\geq \la}M @>>>  0\\ @VVV @VVV @VVV @VVV @VVV \\  0 @>>> {\rm gr}^\la M @>>>  {\rm gr}^\la N @>>>  T'/T@>>>  0\\ @VVV @VVV @VVV @VVV @VVV \\  0 @>>> 0 @>>>  0 @>>>  0 @>>>  0. \end{CD}\] 
     
If  $ F^{\geq \la}f$ and $F^{> \la}f$ are epi-monic in $\sC$, $\gr^\la f$ is epi in $\sC$ and the second and third rows are exact in $\sA$. By the snake lemma, it follows that the fourth row is also exact in ${\sA}$,  hence $\gr^\la f$ is epi-monic in $\sC$.  

For the converse, we argue by descending induction on $\la$: we assume that $F^{>\la} f$ and $\gr^\la f$ are epi, then the composition $F^{>\la} N \to F^{\geq\la} N \to \Coker F^{\geq\la}  f$ factors through a morphism $\gr^\la N \to \Coker F^{\geq\la}  f$ whose composition with the epi $\gr^\la f$ is $0$, hence is itself  $0$.    Therefore $F^{\geq\la}  f$ is epi.
  
   \smallskip  $(1)\Leftrightarrow (4)$ follows from $(1)\Leftrightarrow (3)$ since $\gr^\la_F = \gr^{-\la}_{\check F}$.
   \qed
  
  \begin{cor}\label{efsf}  Any \emph{strongly exact} functorial decreasing separated exhaustive left-continuous filtration by strict subobjects (indexed by $\La$) is a slope filtration. 
  
  A fortiori, if $\sC$ is abelian, any \emph{exact} functorial decreasing separated exhaustive left-continuous filtration (indexed by $\La$) on objects of $\sC$ is a slope filtration. 
  \end{cor}
  
  \proof  Indeed, exactness implies that for any $\la$,  the filtration of $F^{\geq \la}N$ (\resp $N/F^{\geq \la}N$) is induced by the  filtration of $N$. On the other hand,
   strong exactness implies, via item (3) of the previous lemma, that the function 
   $ M \mapsto \deg M = \sum_\la \la \cdot \rk \gr^\la M$  is additive with respect to short exact sequences. \qed
   
 \begin{ex}\label{weight} In the abelian category of rational mixed Hodge structures (Deligne), the decreasing filtration attached to the (increasing) {\it weight filtration} (which is exact \cite{De1})
   $$F^{\geq \la} M= W_{[-\la]}\, M $$ is an exact slope filtration. It also induces a strongly exact slope filtration on the quasi-abelian category of (torsionfree) integral mixed Hodge structures, hence also on the full quasi-abelian subcategory of $1$-motives over $\C$, \cf \cite{De2}. The left abelian envelope of the latter category was considered in \cite{BRS}.
 \end{ex}  
  
   \begin{rems} 1) The Newton polygons associated to an exact slope filtration are additive in short exact sequences. 
   
 \smallskip\noindent 2) For any morphism $f: M\to N$, it follows from lemma \ref{SSQAB} that $\Ker  \gr f= \gr \Ker  \gr f$, whence a canonical morphism
  $$   \gr \Ker  f\to \Ker  \gr f,$$  which is neither injective nor surjective in general (as one can see in the  short exact sequence \eqref{exseqHN}). It is an isomorphism for strongly exact filtrations.
   \end{rems}

\medskip \subsection{Characterization of (strongly) exact slope filtrations.}

 We now assume that $F^{\geq .}$ is a slope filtration.
    
 \begin{thm}\label{T2}  The following properties are equivalent: 
  \begin{enumerate} 
  \item  $F^{\geq .}$ is exact (\resp strongly exact),
    \smallskip \item any non-zero strict subobject (\resp any non-zero subobject) of a semistable object has the same slope,
 \smallskip \item any non-zero strict quotient (\resp any non-zero quotient) of a semistable object has the same slope,
 \smallskip \item there is no non-zero strict morphism (\resp non-zero morphism) between semistable objects of different slopes. \end{enumerate}
  \end{thm}

  \smallskip  Items (2) and (3) of the theorem justify the common use of the terminology \emph{isoclinic} or \emph{pure} (of slope $\la$) instead of semistable, in the case of a strongly exact filtration. 
   
    \smallskip  \begin{cor}\label{SSS} If $F^{\geq .}$ is a strongly exact slope filtration, then the \emph{stable} objects are the \emph{simple} objects of $\sC$. 
\end{cor} 

This follows from item (2) of \ref{T2}.  \qed 
 
 \proof  $(2)\Leftrightarrow (3)$ follows from the fact that $F^{\geq .}$ and its dual $\check F^{\geq .}$ are simultaneously exact (\resp strong exact) or not. Note that in these items, the subobject (\resp quotient) is necessarily semistable (by item (6) of \ref{PL2}).
 
 \smallskip $(2)+(3) \Rightarrow (4)$: Let $f:L\to M$ be a strict morphism (\resp a morphism) between semistable objects of slopes $\la$ and $\nu$ respectively. Let $L\overset{e}{\surj} M\inj N$ be its canonical factorization, with $e $ strict epi. Then $(2)+(3)$ imply $\mu(M)= \la= \nu$. 

  \smallskip $(4) \Rightarrow (2)$: Let $N$ be a semistable object of slope $\nu$, $M $ be a non-zero strict subobject (\resp subobject) of $N$, and $L$ be the universal destabilizing subobject of $M$. Then $\mu(L)=\mu(\IM(L)) = \nu$ by (4). In particular, $\IM L$ and $N/\IM L$ are semistable of slope $\nu$ (or zero). Again, the universal destabilizing subobject of $M/  L $, which is a subobject $N/\IM L$,  has slope $\nu$ if it is non-zero, but this contradicts the definition of $L$. Thus $M=L$ is semistable of slope $\nu$.

    \smallskip
 $(1)\Leftrightarrow (2)$  It suffices to prove that the following two assertions (in their respective avatars).
  
 \smallskip $i)$ If any monic (\resp strict monic) is strictly compatible with $ F^{\geq .}$, then any non-zero subobject (\resp strict subobject) of a semistable object is semistable of the same slope. 
 
 \smallskip $ii)$ If any non-zero subobject (\resp strict subobject, \resp strict quotient) of a semistable object is semistable of the same slope, then any monic (\resp strict monic, \resp strict {epi}) is strictly compatible with $ F^{\geq .}$.

\smallskip Proof of $(i)$. If $f: M\inj N$ is monic (\resp strict monic), with $N$ semistable of slope $\la$, the functoriality of $F^{\geq .}$ implies that $F^{>\la}M=0$, and the strict compatibility of $f$ with  $F^{\geq .}$ implies that  $f(F^{>\la}M) = f(M)\cap F^{>\la}N= M$. Hence $M= \gr^\la\, M$ is zero or semistable of slope $\la$.
  
\smallskip Proof of $(ii)$. Let $f: M\inj N$ be a non-zero monic. Arguing by descending induction on $\la$, we assume that $f(F^{>\la}M) = f(M)\cap F^{>\la }N$ and have to show that $f(F^{\geq\la}M) = f(M)\cap F^{\geq\la }N$. Then $\gr^\la\, M\to \gr^\la\, N$ is monic (\resp strict monic), being a push-out of  $F^{\geq \la}M\to F^{\geq \la}N$ by a strict epi, and that the natural morphism 
$$  \frac{f(M)\cap F^{\geq\la }N}{f(F^{\geq\la}M) } \to \frac{F^{\geq\la }N}{f(F^{\geq\la}M)+F^{>\la }N} = \frac{\gr^\la\, N}{\IM (\gr^\la\, M\to \gr^\la\, N)}$$  is  monic (\resp strict monic). By assumption, this implies that $\,\frac{f(M)\cap F^{\geq\la }N}{f(F^{\geq\la}M) }\,$ is zero or semistable of slope $\la$. Since $$\,\frac{f^{-1}( F^{\geq\la }N)}{  F^{\geq\la}M  }\to \frac{f(M)\cap F^{\geq\la }N}{f(F^{\geq\la}M) }\,$$ is {epi-monic} and  $\,\frac{f^{-1}( F^{\geq\la }N)}{  F^{\geq\la}M  }\,$ has slopes $\leq \la$, we conclude that 
$f(F^{\geq\la}M) = f(M)\cap F^{\geq\la }N$.

\smallskip Let $g: N\surj P$ is a non-zero strict {epi}, with kernel denoted by $f: M\inj N$. Taking into account  
the previous step, we know that $f$ is strictly compatible with $F^{\geq .}$ 
Arguing by ascending induction on $\la$, we assume that $ F^{\geq \la}P  = g(  F^{\geq \la }N)$ and have to show that $ F^{> \la}P  = g(  F^{> \la }N)$. One has a commutative diagram with exact rows and columns

 \[\begin{CD} 0 @>>> 0 @>>>  0 @>>>  0 @>>>  0\\ @VVV @VVV @VVV @VVV @VVV \\  0 @>>> F^{> \la}M @>>>  F^{> \la}N @>>>  \frac{F^{> \la}N}{F^{> \la}M} @>>>  0 \\ @VVV @VVV @VVV @VVV @VVV \\  0 @>>> F^{\geq \la}M @>>>  F^{\geq \la}N @>>>  F^{\geq \la}P @>>>  0\\ @VVV @VVV @VVV @VVV @VVV \\  0 @>>> {\rm gr}^\la M @>>>  {\rm gr}^\la N @>>> \frac{\gr^{  \la}\,N}{\gr^{  \la}\,M} @>>>  0\\ @VVV @VVV @VVV @VVV @VVV \\  0 @>>> 0 @>>>  0 @>>>  0 @>>>  0. \end{CD}\] 
 
We have to show that the natural strict monic $  \frac{F^{> \la}N}{F^{> \la}M} \to  {F^{> \la}P}$  is an isomorphism, or equivalently, that the natural strict {epi}
$ \frac{\gr^{  \la}\,N}{\gr^{  \la}\,M} \to  {\gr^{  \la}\,P}$  is an isomorphism.
 By assumption, $\frac{\gr^{  \la}\,N}{\gr^{  \la}\,M} $ is zero or semistable of slope $\la$. Thus the morphism $ F^{\geq\la}P\to \frac{\gr^{  \la}\,N}{\gr^{  \la}\,M} $ factors through  $\gr^\la \, P$. 
\qed

\medskip \subsection{Split slope filtrations.}\label{ssf}  

\begin{defn} A slope filtration $F^{\geq .}$ is \emph{split} if $\gr \cong id$ (as a functor). 
\end{defn}

In other words, the canonical flag $\sF(M)$ splits, functorially in $M$. 

\begin{exs}\label{exsp} 1) Among our four basic ``additive examples", the Turrittin-Levelt and Hasse-Arf filtrations are split (see \ref{bsf} below for an explanation of this fact), as well as the Dieudonn\'e-Manin filtration if $\phi$ is bijective. 

\smallskip\noindent  {2)} Any vector bundle on a smooth connected projective curve of genus $\leq 1$ is a direct sum of semistable bundles, \ie $\gr$ is the (isomorphic to the) identity on objects; however, it is not identity on morphisms, and, as we have seen, the Harder-Narasimhan filtration is not exact.

\smallskip\noindent  {3)}  Exact slope filtrations may be non-split, even in the abelian case, \cf \ref{weight} or \ref{EXX}.
  \end{exs}

\begin{lemma}\label{ssex} \begin{enumerate}
\item Any split slope filtration is strongly exact.
\item  In the presence of a split slope filtration, the additive groups of morphisms $\sC(M,N)$ are naturally graded.  
\end{enumerate}
 \end{lemma}

 \proof  (1) follows from item (3) of lemma \ref{eqex}. (2) is immediate.
 \qed

 \bigskip 
 
 \section{Slope filtrations and triangulated categories.}
 
 \medskip\subsection{Extension of a slope filtration from $\sC$ to its left abelian envelope.} We set
 $$ \bar \La =  \La \cup \{ {\infty}\}$$ (totally ordered set with maximum ${\infty}$). 
 
 \smallskip Let $\sC$ be a quasi-abelian category, with left abelian envelope $\sA$.   According to corollary
  \ref{STABAB0}, to give a slope filtration indexed by $\La$ on $\sC$ (with respect to a fixed rank function $\rk$) is equivalent to giving a homomorphism 
  $$\deg: \,K_0(\sA)\to \La$$ that is \emph{non-negative on torsion classes $[Q], \, Q\in Ob \, \sT$}.    
  
  Let us assume that {\it $\sT$ consists precisely of objects of rank $0$}. 
   One can extend the slope function $\mu = \deg/ \rk$ to a function 
    $$\mu :  Sk\, \sA \to \bar \La  $$
 which is $  \infty$ exactly on $ Sk \, \sT$.
 
 The slope filtration on $\sC$ then {\it extends} to a unique decreasing separated exhaustive functorial left-continuous filtration {\it on $\sA$  indexed by $\bar \La $}: with the notation of \ref{P1}, for any $A \in Ob\, \sA$, $F^{\geq\la} A$ is the pull-back of $F^{\geq\la} (A/ A_{tor})$ by $A\to A/A_{tor}$, and 
  $ F^{\geq  \infty} A = A_{tor}$.

\medskip \subsection{Stability structures on a triangulated category.}\label{sstc} Let $\sD$ be an essentially small triangulated category, and let
 $$\rk : \; K_0(\sD) \to \Z $$ be a group homomorphism. 
 The following definition is a slight reformulation of Bridgeland's notion of ``stability condition". 

\smallskip
\begin{defn} A \emph{stability structure} (or \emph{$s$-structure}) on $\sD$ consists of a group homomorphism  
$$\deg: \; K_0(\sD) \to  \La $$  called the \emph{degree function}, and full additive subcategories $\sC(\la)$ for each $\la \in \bar \La$, such that 

\begin{enumerate} 
\item the values of $\,\rk \,$ on $\,\sC(\la)\setminus 0\,$ are positive if $\la \in \La$ (\resp $\rk = 0$ on $\sC( {\infty})$),

\smallskip \item one has $\,\deg = \la \cdot  \rk  \,$ on $\,\sC(\la)\setminus 0\,$ if $\,\la \in \La$ (\resp $\,\deg \in   \La_{>0}\,$ on $\,\sC({\infty})$),

\smallskip \item for any $E\in \sC(\la) $ and  $ E'\in \sC(\la')$, one has 
$$\sD(E[n],E'[n'])= 0\;\;\; \text{if}\;\; (n,\la)> (n',\la')\,$$ (with respect to the lexicographic order in  $\,\Z \times \bar \La$),

\smallskip \item for any non-zero object $E$ of $\sD$, there is a finite sequence
$$(n_1, \la_1)> \ldots > (n_r, \la_r)\;\;\;\; \text{in}\;\;\,\Z \times \bar \La$$ and a collection of triangles (Postnikov tower)
$$\begin{matrix} 0 = E_0 & \longrightarrow & E_1 &  \longrightarrow & E_2 &  \to  \cdots     \to  & E_{r-1} & \longrightarrow & E_r= E\\
& \overset{+1}{\nwarrow}    F_1  \swarrow && \overset{+1}{\nwarrow}   F_2    \swarrow  &   &&    & \overset{+1}{\nwarrow}    F_r    \swarrow & \end{matrix}$$
\smallskip with $F_j\in \sC(\la_j)[n_j]$.

\end{enumerate} 

\end{defn}

  \smallskip
 \begin{rem} For $\La= \R$, this corresponds to Bridgeland's notion of  ``stability condition" modulo the following dictionary. Bridgeland's ``central charge" is
 $$Z(E) = -\deg E + \sqrt {-1} . \rk E\; \in \C.$$  The categories $\sP(\phi), \, \phi\in \R,$ from \cite[def. 1.1]{Br} are the shifts $\sC(\la)[n]$, according to the rule
 $$(n, \la) \mapsto \phi = n + \frac{1}{\pi}Arctg (-\frac{1}{\la}),$$  which induces an increasing bijection $\Z \times \bar \La\cong \R\,$ (here, $Arctg$ takes its values in $]0,\pi]$). By working directly with the totally ordered set $\Z \times \bar \La$ instead of $\R$, all arguments of \cite{Br} apply, mutatis mutandis, without assuming $\La= \R$.
   \end{rem}
 
\medskip For any interval $I\subset \bar \La$, we denote by $\sC(I)$ the smallest strictly full extension-closed subcategory of $\sD$ containing the objects of $\sC(\la),\; \la \in I$. This is nothing but the full subcategory of $\sD$ consisting of objects that admit a Postnikov tower as above with $n_j=0, \la_j \in I$. 
 
 Similarly, for any interval $J\subset \Z \times \bar \La$, we denote by $\sD(J)$ the smallest strictly full extension-closed subcategory of $\sD$ containing the objects of $\sC(\la)[n],\; (n,\la) \in J$.

\begin{lemma} \begin{enumerate}
\item A Postnikov tower as in $(4)$ is unique (up to unique isomorphism). 

\smallskip \item The subcategories $\sC(\la)$ are abelian, and all the subcategories $\sC(I)$ are quasi-abelian. The short exact sequences in $\sC(I)$ are the triangles in $\sD$ whose vertices belong to $\sC(I)$.
\end{enumerate}
 \end{lemma}

\proof (1): \cf \cite[4.1]{GKR}.

\noindent (2): \cf \cite[5.2, 4.3]{Br}.
\qed

\begin{thm}[Bridgeland] To give an $s$-structure on $\sD$ is equivalent to giving 
\begin{itemize} \item a bounded $t$-structure on $\sD$, 

\item a quasi-abelian full subcategory $\sC$ of the heart $\sA$ of this $t$-structure, such that $\sA$ is the left abelian envelope of $\sC$ 
and the associated torsion subcategory $ \sT$ (\cf \ref{P1}) consists of the objects of $\sA$ of rank $0$,  and
 \item  a slope filtration on $\sC$ whose degree function $$\deg: \, K_0(\sD)=K_0(\sA)= K_0(\sC)\to \La$$ is positive on non-zero torsion classes $[Q],\, Q\in Ob\, \sT$.
  \end{itemize}
 In fact, the $t$-structure attached to a given $s$-structure is 
 $$ \sD^{\leq 0} = \sD(\N \times \bar \La) ,\;\; \sD^{\geq 0} =\sD((-\N) \times \bar \La) ,$$
 one has $$\, \sA= \sC(\bar \La),\;\; \sC = \sC(\La) $$ and 
 $\sC(\la)$ is the full subcategory of $\sC$ consisting of $0$ and the semistable objects of slope $\la$.\end{thm}

\proof \cf \cite[5.3]{Br}.   
 \qed

 \begin{exs}\label{extri} The Harder-Narasimhan filtration on vector bundles (Mumford or Gieseker version, \cf \ref{exs1}) on a polarized smooth normal connected projective variety $X$ satisfies the assumption of the corollary and induces a canonical $s$-structure on the bounded derived category $D^b(\sO_X)$. We refer to \cite{GKR} for a discussion of this $s$-structure and more exotic ones.   \end{exs} 
 
 \begin{rem} Actually, Bridgeland \cite{Br} allows rank functions 
 $\,K_0(\sD) \to \R\,$ with real values as well; the set of such $s$-structures then acquires a $GL^+_2(\R)$-action, coming from the homographic action of $GL^+_2(\R)$ on pairs $(\rk, \deg)$.
 \end{rem}

 \newpage

\addtocontents{toc}{{\bf II. Behaviour of slope filtrations with respect to a tensor product.} \hfill\thepage}
\
\bigskip
\begin{center}
\large\bf II. Behaviour of slope filtrations with respect to a tensor product.
\end{center}
\bigskip

 In our first four basic examples, the underlying quasi-abelian categories are endowed with a natural tensor product $\otimes$, and the slope filtrations exhibit rather different behaviours with respect to $\otimes$ and duality: in the Turrittin-Levelt and Hasse-Arf cases, the slopes are non-negative and invariant under duality; in the Dieudonn\'e-Manin and Harder-Narasimhan cases, the slopes are changed to their opposite by duality.
 
In these two types of slope filtrations, the breaks remain bounded or grow linearly, respectively, when one takes arbitrarily large tensor (or symmetric) powers.
 
 The aim of this chapter is to analyze these two types (which we call $\otimes$-bounded and $\otimes$-multiplicative respectively) in the general context of quasi-tannakian categories, that are quasi-abelian generalizations of tannakian categories.

 \section{Quasi-tannakian categories.}\label{qtc}

\subsection{Quasi-tannakian categories and rank function.}\label{qtrf} Let $F$ be a field of characteristic $0$. 

\begin{defn} An $F$-linear symmetric monoidal category $(\sC, \otimes)$ is {\emph{quasi-tannakian over $F$}} if
\begin{enumerate}{\it
\item it is quasi-abelian,
\item it is rigid (\ie any object has a (strong) dual, \cf \cite{Saa}),
\item $End\,\un = F\,$ (where $\un$ denotes the unit object),
\item there is an exact faithful (symmetric) rigid monoidal functor ${\omega}$ from $\sC$ to the monoidal category $\,Vec_{F'}\,$ of finite dimensional vector spaces over some fixed extension $F'/F$.}
 \end{enumerate}
 \end{defn}

\begin{rems}  If one replaces (1) by the stronger condition
 
\smallskip  $\;$ (1)' {\it it is abelian},
  
\smallskip  one recovers the definition of a tannakian category over $F$.
  
 \medskip On the other hand, it is well-known that $\omega$ is automatically compatible with duality (\cf \cite[I.5.2.2]{Saa}).  Note that the functor $()^\vee: \sC^{op}\to \sC$ being an equivalence, it respects monics and {epi}s, kernels and cokernels, images and coimages. Note also that ${\omega}$ can be used to detect when a morphism in $\sC$ is non-zero, \resp monic,   \resp {epi}. \end{rems}

\subsection{Quasi-tannakian rank.}\label{qtr} By rigidity, there is a notion of trace of any endomorphism, and of rank 
$$\rk \, M := tr\,id_M .$$ One has 
$\rk \,M= \dim_{F'} \omega(M)$, which is a natural integer (here, the fact that ${\rm char}\, F = 0$ is essential). This shows that $\rk$ takes the value $0$ only on the zero object, and is additive on short exact sequences (since $\omega$ is exact). 
 Thus $\rk$ is a rank function in the sense of \ref{RANK}.
 
\smallskip Tensor product and duality make  $K_0(\sC)$ into a commutative ring with involution, and $\rk$ defines a {\it ring homomorphism} 
$$\,\rk:\,K_0(\sC)\to \Z $$  with $\rk M= \rk M^\vee$. 

  \bigskip In the sequel, $\sC$ will be an {\it essentially small quasi-tannakian category over $F$, equipped with its canonical rank function $\rk$}.

\subsection{The semisimple tannakian quotient category.}
  Being quasi-abelian, $\sC$ is pseudo-abelian\footnote{\ie idempotent endomorphisms have kernels.}, and since ${\rm char}\, F = 0$, it is possible to define symmetric and exterior powers of an object as direct summands of its tensor powers.
 
 Condition (4) implies 
 
 $\,$ (4)'  {\it For any object $M$, $\displaystyle\bigwedge^{\rk M +1}\, M=0.$}
    
\smallskip Essentially small pseudo-abelian $F$-linear symmetric monoidal categories satisfying (2), (3) and (4)' have been studied in \cite{AK} and by P. O'Sullivan (independently). 

\begin{prop}\label{pi}\cite[\S 9]{AK}\cite{OS2} Assume $ \sC $ satisfies (2), (3), (4'). Then
the maximal $\otimes$-ideal $\sN$ of $\sC$ is locally nilpotent,  $\bar \sC = \sC/\sN$ is a semisimple tannakian category (with the same objects as $\sC$). The canonical $\otimes$-functor $\sC\to \bar\sC$ is conservative (\ie any morphism $f$ in $\sC$ is an isomorphism if $\omega(f)$ is an isomorphism) and full, and $\sk \sC = \sk \bar\sC$. 
\end{prop}

\subsection{O'Sullivan's description.}\label{OSD} In fact, O'Sullivan went further and elucidated the  structure of $\otimes$-categories satisfying (2), (3) and (4)'.  Although we will make little use of it, we briefly survey this enlightening viewpoint (\cf \cite{OS2} and \cite[3.7]{A4} for more detail).

\smallskip The functor $\sC\to \bar \sC$ actually admits a $\otimes$-section $\sigma$.

  Let us first assume that $F'=F$. Then $\omega\circ \sigma$ induces an equivalence $$\bar\sC \cong Rep \, H,$$ where $Rep\, H$ denotes the tannakian category of finite dimensional representations of the  proreductive group $H = Aut^\otimes (\omega\circ \sigma)$ over $F$. 
Moreover, there is an integral affine scheme $X = \Spec A$ with $H$-action, with $A^H= F$,  a $F$-point $x\in X$ fixed under  $H$, and an equivalence 
$$ \sC \cong Vec(H, X)$$ between $\sC$ and the category of $H$-equivariant vector bundles on $X$, such that the projection $\sC\to \bar\sC$ corresponds to the functor ``fiber at $x$":  $Vec(H, X)\to Rep \, H$ it turns out that, in this situation, any object of $Vec(H, X)$ is of the form $V\otimes \sO_X$ for some object $V\in Rep\, H$.

\smallskip In general, the result is similar:  $A$ becomes an integral\footnote{in the idealistic sense.} algebra in $Ind\, \bar\sC$, $x$ an augmentation $A\to \un$, $Vec(H, X)$ is replaced by the category $Proj_A$ whose objects are those of $\sC$ and whose morphisms are given by  $ \, Hom^{A-linear}_{Ind\, \bar\sC}(M\otimes A, N\otimes  A).$

    \begin{exs} 1) The quasi-tannakian category of finite-dimensional $F$-vector spaces with a (separated exhaustive) $\Z$-filtration is $\otimes$-equivalent to $Vec(\bG_m, \bA^1)$ (for the natural action of $\bG_m $ on $\bA^1$ by homotheties).
 
 \smallskip\noindent 2) The quasi-tannakian category of vector bundles over $\bP^1$ is $\otimes$-equivalent to $Vec(\bG_m, \bA^2)$ (this is a reformulation of Grothendieck's theorem).
  
 \smallskip\noindent 3) The tannakian category of finite-dimensional $F$-vector spaces with a nilpotent endomorphism is $\otimes$-equivalent to $Vec(SL_2, \bA^2)$ (this is a reformulation of the Jacobson-Morozov theorem). 
 
The latter embeds as a full subcategory in the abelian category $Mod(SL_2, \bA^2)$ of $SL_2$-equivariant coherent sheaves on $\bA^2$. This subcategory is unstable under taking subobjects or quotients, but monics and {epi}s in $Vec(SL_2, \bA^2)$ remain so in $Mod(SL_2, \bA^2)$, respectively. An object of $Mod(SL_2, \bA^2)$ lies in $Vec(SL_2, \bA^2)$ if and only if it is reflexive (\ie isomorphic to its bidual).  
   \end{exs}

 \begin{lemma}\label{propqt1} Let $\sC$ be a quasi-tannakian category. \begin{enumerate}  
\item $\otimes$ is exact in both arguments and $()^\vee$ is exact. 
\item a morphism $f\in \sC(M,N)$ is monic (\resp epi) if and only $\omega(f)$ is injective (\resp surjective).

\noindent In that case, one has $\rk \, M \leq \rk \, N$ (\resp $\rk \, M \geq \rk \, N$).  \end{enumerate}
\end{lemma} 

\proof (1). The functor $()^\vee: \sC^{op}\to \sC$ preserves kernels and cokernels, hence is exact. 

Let $MOD_A$ be the abelian monoidal category of $A$-modules in $Ind\, \bar\sC$ and $MOD_{ \omega(A)}$ be the abelian monoidal category of $\omega(A)$-modules in $ Ind\, Vec_{F'}$. Notice that $\omega$ extends to a faithful exact monoidal functor $ MOD_A \to  MOD_{ \omega(A)}$. Then, in view of O'Sullivan's monoidal equivalence $$\sC\cong Proj_A,\;\; M\mapsto M\otimes A,$$ the bi-exactness of $\otimes$ in $\sC$ follows from the exactness of the endofunctors $$- \otimes_{\omega(A)}\, (\omega(M)\otimes_{F'} \omega(A)),\;\; (\omega(M)\otimes_{F'} \omega(A)) \otimes_{\omega(A)} - $$  in $MOD_{ \omega(A)}$, free $ \omega(A)$-modules being flat.

\smallskip\noindent (2) Since $\omega$ is faithful, $\omega(f)$ injective (\resp surjective) implies $f$ monic (\resp epi). For the converse, by duality, it suffices to treat the case of a monic $f: M\inj N$ in $\sC \cong Proj_A$.  Let $W$ be the kernel of $f$ in $MOD_A$. If $W\neq 0$,  there is a non-zero morphism $P\to W$ in $Ind\,\bar\sC$ with $P$ in $\sC$. Whence a non-zero morphism $P\otimes A\to W$ in $MOD_A$, and by composition a non-zero morphism $N\otimes A \to M\otimes A$ in $Proj_A$ such that the composed morphism $P\otimes A \to N\otimes A$ is zero. This is a contradiction, thus  $f$ remains monic in $MOD_A$. Therefore $\omega(f\otimes 1_A)$ is also monic in $MOD_{ \omega(A)}$, hence $\omega(f)$ is injective. 

The last assertion is immediate. \qed

\begin{rems}\label{remqt} 1) \noindent  In concrete situations,  the assertions of the proposition can be checked direcly, without reference to O'Sullivan's theory.

If $\sC$ is abelian, the proposition is standard and may be obtained directly using $\omega$.
  
\smallskip \noindent 2) Items (2) and (3) imply that any constant function $\mu$ with values in $\La$ defines a slope filtration on $\sC$. 

\smallskip \noindent 3) The description of quasi-tannakian in terms of equivariant vector bundles (or of objects of $Proj_A$ in the non-neutral case) allows to extend the Rees deformation of \ref{Rees} to this setting. Applying $\omega$, on gets a finitely generated $F'[[x^{\La_M}]]$-module which is a deformation from $\omega(\gr M)$ to $\omega(M)$.

\smallskip \noindent 4) There is a natural surjective ring homomorphism  $K_0(\bar \sC) \to K_0(\sC)$\footnote{in fact, $K_0(\bar \sC)$ is the Grothendieck group of $\sC$ with respect to split short exact sequences.}, that can be identified with the standard morphism $R(H)\to K_H(X)$ between the representation ring and the equivariant $K$-theory ring, when $\sC\cong Vec(H, X)$.

 \smallskip \noindent 5) It is an open problem to determine which categories of type $Vec(H, X)$ are quasi-abelian. 
 
 \smallskip \noindent 6) Any slope function on $\sC$ induces a slope function on $\bar \sC$. The corresponding slope filtrations are compatible if and only if any object $M$ of $\sC$ is a direct sum of semistable objects, \ie $M = \gr M$. \end{rems}

\bigskip\section{Invertible objects and determinantal slope filtrations.}

\subsection{Determinants.}\label{deter} The invertible objects with respect to $\otimes$ are the rank one objects (the inverse being the dual).  Any non-zero morphism between invertible objects is monic-{epi}.

 We denote by ${\rm Pic}\, \sC$ the {\it Picard group} of $\sC$, \ie the group of isomorphism classes of rank one objects, with respect to $\otimes$. Since $\sk \sC = \sk \bar\sC$, one has  ${\rm Pic}\, \sC= {\rm Pic}\, \bar\sC$ (which is identified with the group of $F$-characters of $H$ in case $\sC \cong Vec(H,X)$). 

\smallskip For any non-zero object $M$, its {\it determinant} 
$$\det M = \bigwedge^{\rk M}\, M\,$$  is an invertible object. One has a canonical isomorphism 
\begin{equation}\label{dualtens} M^\vee \cong (\bigwedge^{\rk M -1}\, M) \otimes (\det M)^\vee.\end{equation}
There are two ways to see this. One can use the fact that there is a natural $\otimes$-functor $Rep_F GL(\rk M)\to \sC$ sending the standard representation to $M$ (\cf \eg \cite[3.21]{A4}), and that such an isomorphism is already available in $Rep_F GL(\rk M)$. Or one can use the fact the functor $\sC\to \bar\sC$ being conservative and full, it is essentially injective (\ie two objects of $\sC$ are isomorphic if and only if their images in $\bar\sC$ are isomorphic), and that such an isomorphism is actually known to be avalaible in any semisimple tannakian category.

\begin{lemma} The rule $\displaystyle M \mapsto  {\det}\, M   $ induces a surjective group homomorphism
$$ \, K_0(\sC)\overset{\det}{\to} {\rm Pic}\, \sC $$ with $\det ([M ]^\vee)=  \det [M ]^{-1}$. The natural set-theoretic map ${\rm Pic}\, \sC\to K_0(\sC)$ is a section of $\det$.
\end{lemma} 
 
 \proof In order to establish the existence, one has to see that ${\det}$ is multiplicative on short exact sequences $0\to M\to N\to P\to 0$. This is seen as usual by introducing the Koszul filtration by strict subobjects
 $$ K^i(\bigwedge^j\,N) = \IM(\bigwedge^i\,M  \otimes  \bigwedge^{j-i}\,N\to \bigwedge^j\,N) \;\;\; (i\leq j)$$
 with $$K^i(\bigwedge^j\,N)/K^{i+1}(\bigwedge^j\,N) \cong \bigwedge^i\,M  \otimes  \bigwedge^{j-i}\,P.$$
  The surjectivity and the other assertions are straightforward .
   \qed

   \begin{prop}\label{propqt2} \begin{enumerate}
 \item For any objects $M, N$, there is an isomorphism of invertible objects
$$\det (M\otimes N) \cong \det M^{\otimes \rk N}\otimes  \det N^{\otimes \rk M}.$$   For any positive integer $m$, and any positive integer $n$ less than $\rk N$, there are isomorphisms of invertible objects 
$$\det (S^m M ) \cong \det M^{\otimes r} ,\;\; \det (\bigwedge^n M ) \cong   \det N^{\otimes s},$$  with 
$$r=  \frac {(m+\rk M - 1)!}{(m-1)!\rk M !},\;\; s= \frac{(\rk N-1)!}{n!(\rk N -n)!}.$$

  \item For any morphism $f: M\to N$ that is monic-{epi}, $\displaystyle\,\det f = \bigwedge^{\rk M} \, f$ is a epi-monic morphism of invertible objects.

\smallskip \item A morphism $ M\overset{f}{\to} N$ is an isomorphism if and only if $\,\rk M = \rk N\,$ and $\displaystyle\,\det f  $ is an isomorphism.  \end{enumerate} \end{prop}
 
  \proof  
 (1) is be proven in the same way as \eqref{dualtens}.

\smallskip\noindent (2) By \ref{propqt1}, $\omega(f)$ is a bijective linear map, hence $\det \omega(f) = \omega (\det f)$ is non-zero. Therefore, $\det f $ is a non-zero morphism between invertible, hence monic-{epi}.

 \smallskip\noindent  (3): if $f$ is an isomorphism, so is $\det f$ and $\rk M = \rk N$. The converse follows from the fact that $\displaystyle (\bigwedge^{\rk M-1} f) \otimes ({\det}f)^{-1}$ is then left inverse to $f$.
\qed

\begin{cor}\label{Cdet1}  Let $\sC'$ be a quasi-tannakian category over an extension $F'/F$, and let $\theta: \sC \to \sC'$ be an $F$-linear $\otimes$-functor. 

Then $\theta$ is conservative if and only if any non-zero monic $m: L\inj \un$ such that $\theta(m)$ is an isomorphism is an isomorphism.
\end{cor} 

\proof The ``only if" part is obvious. To prove the ``if" part, let us consider a morphism $ f: M\to N$ be a morphism such that $\theta(f)$ is an isomorphism. Then $\rk M= \rk N$ and $\theta(\det f)= \det \theta(f)$ is an isomorphism. Tensoring $\theta(\det f$ with $1_{\det N^\vee}$, one gets a non-zero monic $m: L\inj \un$. By assumption, this is an isomorphism, hence $\det f$ is an isomorphism, and so is $f$ by the last proposition.
\qed

\begin{cor}\label{Cdet2}
 \item The following conditions are equivalent:
   \begin{enumerate} 
  \item $\sC$ is abelian (hence tannakian),
  \item the unit $\un$ is simple (\ie irreducible),
  \item $\omega$ is conservative.
\end{enumerate}  \end{cor}
  
 \proof Note that $\un$ is simple if and only if any object of rank one $L$ is simple (using $\otimes L^\vee$). 
 
 $(1) \Rightarrow (2)$: in any tannakian category, any object of rank $1$ is simple. 
 
 $(2)\Rightarrow (3)$: follows from the previous corollary.
 
 $(3) \Rightarrow (1)$: Let $f: M\to N$ be monic-{epi}. By \ref{propqt1},  the linear map $\omega(f)$ is bijective. Since $\omega$ is conservative, $f$ is an isomorphism.    \qed

  \begin{prop} For any slope filtration on the quasi-tannakian category $\sC$, the endofunctor $\gr$ of $\sC$  is conservative\footnote{but not essentially injective in general, of course.}.
 \end{prop}
   
  \proof By corollary \ref{Cdet1}, this reduces to the fact that any non-zero monic $m: L\inj \un$ such that $\gr(m)$ is an isomorphism is an isomorphism. If $\gr(m)$ is an isomorphism, then $L$ is semistable of slope $0$, so that $\gr(m)=m$, and $L\cong \un$.  
  \qed

 \subsection{Determinantal slope filtrations.}\label{dsf}

\begin{defn} A slope filtration on $\sC$ is \emph{determinantal} if for any object $M$, 
$$ \deg(M)= \deg(\det M).$$
\end{defn}
In other words, one requires that its degree function $$\deg= \rk \cdot \mu : \,K_0(\sC) \to \La$$ factors through $$ \det :\, K_0(\sC) \to Pic(\sC).$$ 

It is immediate that the set of determinantal slope filtrations is a convex subcone of the cone of all slope filtrations.

    \begin{exs}\label{EXF} 1) Let $\sC$ be the (quasi-tannakian) category of finite-dimensional $F$-vector spaces with a separated exhaustive filtration indexed by $\Z$. Then $K_0(\sC)\cong \Z[t, t^{-1}]$, with $\rk (\sum a_n t^n)= \sum a_n$. Moreover $\det : K_0(\sC) \to Pic(\sC)\cong \Z$ is given by $\det(\sum a_n t^n)= \sum na_n,$ and the determinantal slope function with degree function $\det$ is the original filtration.

 \smallskip\noindent 2) Let $\sC$ be the (quasi-tannakian) category of  vector bundles over a smooth geometrically connected projective curve $X$ over $F$. One has $Pic(\sC)=Pic(X)$. We have seen in \ref{EHN} that any slope function $\mu$ such that $\mu(\sO_X)=0$ is determinantal: in fact, it is a non-negative multiple of the standard (Harder-Narasimhan) slope function.
 \end{exs}

Let us introduce a partial order on  $Pic\, \sC$ as follows:
 
\smallskip {\it $[L]\leq [L']   \Leftrightarrow  $ there is a non-zero morphism $L\to L'$} (clearly, this does not depend on the choice of representatives $L, L'$).

\smallskip\noindent This makes ${\rm Pic }\, \sC$ into an ordered abelian group.

\begin{thm}\label{DETFI} Let  $$\delta:  {\rm Pic }\, \sC \to \La$$ be a {\it non-decreasing homomorphism} (\ie $\,\delta([L])\geq 0 $ if there is a non-zero morphism $\un \to L$). Then
the function $$M\mapsto  \mu(M)= \frac{\delta(\det M)}{\rk M}$$ is the slope function attached to a (unique) determinantal slope filtration. 

Any determinantal slope filtration arises in this way.  
 \end{thm} 
 
\proof  Let $f: M\to N$ be a non-zero monic-{epi}. Then according to item (2) of \ref{propqt2},  $\det f \neq 0$, hence $[\det M] \leq [\det N]$. Therefore $\delta([\det M]) \leq \delta([\det N])$ and $\mu(M)\leq \mu(N)$. 

On the other hand, $\deg(M)= \delta(\det M)$  factors through the composed homomorphism $K_0(\sC)\to {\rm Pic}\, \sC \to \La$, hence is additive on short exact sequences. Therefore $\delta$ gives rise to a (unique) determinantal slope filtration via \ref{T1}.
 The converse is immediate.
   \qed 

\begin{prop}\label{DUAL}  Let $F^{\geq .}$ be a determinantal slope filtration.  \begin{enumerate}  
\item For any non-zero objects $M,N$, \begin{equation}\label{eqI}  \mu(M \otimes N)=  \mu(M )+ \mu(N).\end{equation}
\noindent In particular, if $\La$ is a commutative ring, the rule $$[M]\mapsto \rk M + \epsilon \deg M$$ induces a {\rm ring homomorphism} $$K_0(\sC) \to \La[\epsilon ]/(\epsilon^2).$$

\item For any non-zero object $M$,  
\begin{equation}\label{eqII}\mu(M)= -\mu(M^\vee). \end{equation}
 For any  $\la$, one has
$$F^{\geq \la}(M^\vee) = (F^{>-\la}M)^\perp ,\;\;  (F^{\geq \la} M)^\vee = M^\vee/ F^{>-\la}M, $$ and  a canonical functorial isomorphism
\begin{equation}\label{eqIII} \gr\, (M^\vee) \overset{\cong}{\to} (\gr\, M)^\vee.\end{equation}  Thus the breaks of $M$ are the opposite of the breaks of $M^\vee$. In particular, $M$ is semistable of slope $\la$ if and only if $M^\vee$ is semistable of slope $-\la$.
 
 \item For any positive integer $m$ and any positive integer $n$  less than $\rk M$, one has 
\begin{equation} \mu(S^m M)= m\, \mu(M ), \;\; 
  \mu(\bigwedge^n N)= n\, \mu(N).\end{equation}
  \end{enumerate}\end{prop}

\proof
 (1) comes from the isomorphism $ \det (M\otimes N) \cong \det M^{\otimes \rk N}\otimes  \det N^{\otimes \rk M} $ of item (1) of \ref{propqt2}.   
 
  (2). One has $\deg M^\vee= \deg \det M^\vee= -\deg \det M = -\deg M$, whence \eqref{eqII}.
 
 Let us set $\,\tilde F^{\geq \la}(M^\vee) = (F^{>-\la}M)^\perp \,$ (that is by definition the  kernel of the strict {epi} 
$M^\vee\surj (  F^{> -\la}  M)^\vee $ dual to $F^{> -\la}  M\inj  M$). 
 This defines a separated, exhaustive, left continuous decreasing filtration $\bar F^{\geq .}$ on $\sC$ by strict subobjects, and one has a canonical functorial isomorphism
$$ \tilde \gr^\la M \cong  (\gr^{-\la} M^\vee)^\vee.$$ 
In particular, $\mu( \tilde \gr^\la M)= \la$ by \eqref{eqII}. It follows that $\bar F^{\geq .}$ is a slope filtration with slope function $\mu$, hence $\bar F^{\geq .}=  F^{\geq .}$ (\cf \ref{T1}), and a canonical isomorphism \eqref{eqIII}.

(3) comes from the other isomorphisms  of item (1) of \ref{propqt2}.  
\qed

    \medskip
\begin{rems} 1) The formula $(F^{\geq \la} M)^\vee = M^\vee/ F^{>-\la}M$ means that the $\otimes$-equivalence 
$\sC\to \sC^{op}$ given by duality is compatible with the slope filtration (for the slope function $-\mu$ on  $\sC^{op}$, \cf \ref{remsslo} 2)). 

\smallskip \noindent 2) We do not know whether, conversely, the formulas \eqref{eqI} (or even \eqref{eqI} + \eqref{eqII}) imply that $F^{\geq .}$ is determinantal.
\end{rems} 

   \medskip \begin{prop}\label{STABAB1} Let $F^{\geq .}$ be a determinantal slope filtration. The following properties are equivalent:
  \begin{enumerate}
  \item
  $\delta$ is (strictly) increasing,

\smallskip  \item $\un$ is stable,

\smallskip  \item for any $\la$ and any object $N$, $N$ is stable of slope $\la$ if and only if any {\rm strict} subobject (\resp subquotient) has slope $<\la$ (\resp $>\la$).

\smallskip \item for any $\la$, the full subcategory $\sC(\la)$ of $\sC$ consisting of $0$ and of the semistable objects of slope $\la$ is {\rm abelian}, hence artinian and noetherian (by \ref{FINITELENGTH}); its simple objects are the stable objects of slope $\la$.
    \end{enumerate}
  Under these conditions, the simple objects of $\sC(\la)$ are the stable objects of $\sC$ of slope $\la$. The natural functor from the socle (\ie the full subcategory consisting of semisimple objects) $\sC(\la)_{ssi}$ of $\sC(\la)$ to $\bar\sC$ is fully faithful. 
    \end{prop}
  
 \proof $(1)\Leftrightarrow (2)$ follows from the fact that any non-zero morphism $L\to L'$ between objects of rank $1$ is monic, and gives rise to a monic $L\otimes (L')^\vee \inj \un$. 
  
\smallskip \noindent $(1)\Rightarrow (3)$. Let $L\overset{g}{\inj} N$ be any non-zero subobject of the stable object $N$. In order to show that $\mu(L)<\la$, let us consider the canonical factorization $L\overset{f}{\inj} M\overset{h}{\inj} N$, where $f$ is monic-{epi}, and $h$ is a strict {epi}. If $M\neq  N$, one has $\mu(L)\leq \mu(M)<\mu(N)$ by assumption. If $M=N$, $\det f = \det g$ is a non-zero morphism  $\det L \to \det N$ (\cf item (2) of \ref{propqt2}), whence $\mu(\det L) < \mu(\det N)$ by (1), and $\mu(L)<\mu(N)$ after division by $\rk M=\rk N$.

\smallskip \noindent $(3)\Rightarrow (2)$ is immediate. 

\smallskip \noindent $(1)\Rightarrow (4)$\footnote{another, independent, proof is provided by \ref{STABAB0}.}.  We already know that $\sC(\la)$ is quasi-abelian (\cf \ref{SSQAB}). Let $f: M\to N$ be a non-zero morphism in $\sC(\la)$, and let us consider its canonical factorization $M\surj \Coim f \overset{\bar f}{\to} \IM f \inj N$ in $\sC$. We have $\la = \mu(M)\leq \mu(\Coim f)\leq \mu(\IM f) \leq \mu(N)=\la$, whence equality.  This implies $\delta(\det \Coim f)= \delta(\det \IM f)$. Assuming that  $\delta $ is decreasing, we get that $\det \bar f$ is an isomorphism, an so is $\bar f$ itself. On the other hand additivity of the degree implies that $\Ker  f$ (\resp $\Coker f$) is either $0$ or is of slope $ \la$. Using \ref{SSQAB}, one concludes that $\Ker  f$ and $\Coker f$ belong to $\sC(\la)$. This shows that $\sC(\la)$ is abelian.   
It is immediate that the simple objects in $\sC(\la)$ are the stable objects in $\sC(\la)$. 

  The natural functor $\sC(\la)_{ssi}\to \bar\sC $ is a full, conservative, additive functor between semisimple categories, hence faithful.  
 
\smallskip \noindent $(4)\Rightarrow (2)$. By (4), any non-zero monic $L\inj \un$, with $L$ of slope $0$, is an isomorphism, hence $\un$ is stable.  \qed

    \medskip 
 \begin{exs}\label{EXX} Let $G\subset GL_2$ be the algebraic group over $F$ consisting of matrices of the form $\begin{pmatrix} x & y\\ 0&1\end{pmatrix}$, and let $\sC = Rep\, G$ be the tannakian category of its finite dimensional representations. 
 Then any object of $\sC$ is the restriction of a representation of $GL_2$  and in fact $\bar\sC \cong Rep\, GL_2$, \cf \cite[C5]{OS2}. The group ${\rm Pic}\, \sC= {\rm Pic}\, \bar\sC$ is freely generated by the determinant $Det: G \inj GL_2\to \bG_m$. The standard representation $V$ of $G$ sits in a short exact sequence 
 $$0\to Det \to V \to \un\to 0.$$ 
 If one takes $\delta(Det)= 1$, one gets an exact non-split determinantal slope filtration for which the only stable objects are the powers of $Det$. Thus $\sC(\la)$ consists of direct sums of copies of $Det^{\otimes \la}$ if $\la$ is an integer, and is $\{0\}$ otherwise.
  
\noindent If one takes $\delta(Det)= -1$, one gets a non-exact determinantal slope filtration for which $V$ is stable of slope $-1/2$. One can check that $\sC(\la)$ is $\{0\}$ if $\la$ is not half an integer, whereas if $\la$ is half an integer, $\sC(\la)$ consists of direct sums of copies of objects of the form $Det^{\otimes m} \otimes S^{n} V $ with $\la = -m - n/2, \, m\in \Z, n\in \N$. \end{exs}

 \bigskip
  \subsection{Integrality.}\label{Int}  
 
 \begin{defn} A slope filtration on a quasi-abelian category $\sC$ is said to be \emph{integral} if its \emph{degree} function takes values in $\Z$.
 \end{defn}
 
 This amounts to saying that the \emph{vertices of the Newton polygon of any object $M$ belong to $\Z^2$}. One may then assume that $\La = \Q$.

\smallskip Most slope filtrations of the literature have this property. This is for instance the case in our first four basic examples, for which this property is actually trivial, except for the Hasse-Arf filtration (Hasse-Arf theorem).
 
 The significance of this property is illustrated by the following 
 
   \medskip \begin{prop}\label{INT} Let $F^{\geq .}$ be an integral slope filtration. Assume either that $\sC$ is abelian, {\rm or} that $\sC$ is quasi-tannakian, $F^{\geq .}$ is determinantal and $\un$ is stable.
  
 Then any {\rm semistable} object $N$ such that $\deg N$ and $\rk N$ are relatively prime is {\rm stable}.
   \end{prop} 
 
  \proof  Assume, on the contrary, that $N$ is semistable but not stable. If $\sC$ is abelian,  {\rm or}  $F^{\geq .}$ is determinantal and $\un$ is stable, there is a \emph{strict} non-zero subobject   $M\neq N$ with $\mu(M)=\mu(N)$, \ie $$\rk N\cdot \deg M= \rk M\cdot \deg N.$$ By assumption, $\rk N$ is prime to $\deg N$, hence divides $\rk M$. But  $\rk M< \rk N$, a contradiction. 
  \qed

 \begin{rem}\label{REMF} In many situations, $\sC$ belongs to a family of quasi-tannakian categories $\sC_Y$ indexed by objects $Y$ of a certain small category $\sY$. Morphisms in $\sY$ have a degree, that is a natural integer which is multiplicative with respect to composition of morphisms. 
 
 Any morphism $\phi: Y\to Y'$ of degree $d$ in $\sY$ gives rise to $F$-linear functors 
 $$\phi^\ast: \sC_{Y'}\to \sC_{Y},\; \phi_\ast:  \sC_{Y}\to \sC_{Y'},$$
 $\phi^\ast$ being a $\otimes$-functor, while $\phi_\ast(id \otimes \phi^\ast)= \phi_\ast \otimes id$ and $\rk \phi_\ast M = d\cdot \rk M $.
  
 Moreover, the categories $\sC_Y$ are endowed with slope filtrations $F^{\geq .}_Y$ as in the proposition, which are related to each others by the conditions
 \begin{itemize}
 \item $M' \in \sC_{Y'}$ is semistable of slope $\la'$ if and only if $\phi^\ast M'$ is semistable of slope $d\la'$,
 \item $M  \in \sC_{Y }$ is semistable of slope $\la$ if and only if $\phi_\ast M$ is semistable of slope $ \la/d$.
 \end{itemize} 
 
 One can then use this last condition in order to create objects of non-integral slopes to which the above proposition applies. See for instance \cite{A2} for an application of this technique (to the $p$-adic local monodromy theorem conjectured by Crew). See also \ref{REMFF} below. \end{rem}
 
   \medskip \begin{exs}\label{extr} This setting occurs in the context of our first four basic examples. In the Harder-Narasimhan case, $\sY$ is the category of finite etale coverings $Y$ of the curve $X$ (recall that ${\rm char}\, F= 0$). In the other example, $\sY$ is the category of finite unramified extensions of the ground complete discretely valued field. 
 \end{exs}
 
 \medskip
 
{\it In the next two sections, $\sC$ is a quasi-tannakian category over a field $F$ of characteristic zero, and $F^{\geq .}$ is a slope filtration on $\sC$ indexed by the totally ordered divisible abelian group $\La$.}

\medskip \section{$\otimes$-multiplicative slope filtrations.}

 \subsection{Definition and characterization.}\label{msf}
 
 \begin{defn} $F^{\geq .}$ is \emph{$\otimes$-multiplicative} if it satisfies the following condition:

\smallskip {\it if $M_1$ is semistable of slope $\la_1$ and $M_2$ is semistable of slope $\la_2$, then
 $M_1\otimes M_2$ is semistable of slope $\la_1+\la_2$.} 
 \end{defn}
 
 \begin{exs}\label{EXFF} 1) The Dieudonn\'e-Manin filtration $\otimes$-multiplicative. The Harder-Narasimhan filtration is $\otimes$-multiplicative (a purely algebraic proof, based on geometric invariant theory, of the semistability of the tensor product of two semistable vector bundles appeared in \cite{RR}. Other proofs, relating semistability to numerical effectivity and ampleness, appeared in \cite{Maru} and \cite{Mi}, \cf also \cite[6.4.14]{Laz} and \cite{A6}).  
 
 The examples \ref{EXX} are also $\otimes$-multiplicative.
 
 \smallskip \noindent 2) Let $\La$-$Fil_{F}$ (\resp $\La$-$Bifil_{F}$) be the category of finite-dimensional $F$-vector spaces endowed with a (\resp two) separated, exhaustive, left continuous decreasing filtration (\resp filtrations) indexed by $\La$. This is a quasi-tannakian category over $F$. 
 
 There is a tautological slope filtration on $\La$-$Fil_{F}$, which is obviously $\otimes$-multiplicative (this generalizes example \ref{EXF} 1).
 
On the other hand, the formula 
 $$( F_1^{\geq .}V, F_2^{\geq .}V)  \mapsto  \frac{1}{\dim V}\sum \la (\dim gr^\la_1V + \dim gr^\la_2V)$$ defines a slope function on $\La$-$Bifil_{F}$. 
 
 It is known that the associated slope filtration is $\otimes$-multiplicative if $\La = \Q$, \cf \cite[p. 650]{F2}\footnote{which gives a purely algebraic proof, inspired by arguments of Laffaille. See also \cite{To1}.}. One thus gets a $\otimes$-functor
 
 \medskip\centerline{ $ \La$-$Bifil_{F} \to  \La$-$Fil_{F},\;\;( F_1^{\geq .}V, F_2^{\geq .}V) \mapsto  F^{\geq .}V $.}
    \end{exs}
 
 \medskip
 
  \begin{thm}\label{T3}\begin{enumerate}
   \item There is equivalence between
  \begin{enumerate}
  \smallskip \item $F^{\geq .}$ is $\otimes$-multiplicative,
  \smallskip  \item for any $\la$ and any pair  $(M_1,M_2)$, there is a canonical functorial isomorphism
 \begin{equation}\label{MULT} \displaystyle F^{\geq \la} (M_1\otimes M_2 )\cong \sum_{\la_1+\la_2= \la}\, F^{\geq \la_1}  M_1  \otimes F^{\geq \la_2}  M_2,\end{equation}
 \smallskip \item for any pair  $(M_1,M_2)$, the breaks of $M_1\otimes M_2$ are the sums of a break of $M_1$ and a break of $M_2$.
 \end{enumerate}
 
\medskip  \item Any  $\otimes$-multiplicative slope filtration $F^{\geq .}$ also satisfies:  
 \begin{enumerate} 
  \smallskip  \item $\gr$ is a $\otimes$-functor.
   \smallskip  \item $F^{\geq .}$ \emph{determinantal}. In particular (\ref{DUAL}), there is a canonical functorial isomorphism
 \begin{equation}\label{DUA} F^{\geq \la}(M^\vee) \cong (F^{>-\la}M)^\perp  .\end{equation}
  \item If moreover $\un$ is stable, the category $\sC(0)$ of semistable objects of slope $0$ is   \emph{tannakian}. The tensor product of stable objects of slope $0$ is a direct sum of stable objects.
 \end{enumerate}  \end{enumerate} 
 \end{thm}
 
\smallskip Here, $\displaystyle\sum_{\la_1+\la_2= \la}\, F^{\geq \la_1}  M_1  \otimes F^{\geq \la_2}  M_2 $ denotes the image of the natural morphism $\,\bigoplus_{\la_1+\la_2= \la}\, F^{\geq \la_1}  M_1  \otimes F^{\geq \la_2}  M_2  \to  M_1\otimes M_2 \,$  (note that each $F^{\geq \la_1}  M_1  \otimes F^{\geq \la_2}   M_2$ can be considered as a strict subobject of $M_1\otimes M_2$ since $\otimes$ is bi-exact, \cf item (1) of \ref{propqt1}). 
 
\smallskip   The terminology ``$\otimes$-multiplicative" comes from item (1b): the filtration on a tensor product is the product filtration.

\smallskip \proof (2b).  It suffices to show that $\deg M = \deg \det M$. Since this property is ``stable by extension", it suffices to prove it for $M$ semistable. In that case, $M^{\otimes \rk M}$ is semistable of slope $\mu(M). \rk M= \deg M$. Hence $\det M$, which is a direct summand of $M^{\otimes \rk M}$ is also semistable of slope $\deg M$.
 
\smallskip $(2b) \Rightarrow (2a)$.  Let us first show that for any $M_1, M_2$, the image of the strict monic 
$$f: F^{\geq \la_1}  M_1  \otimes F^{\geq \la_2}  M_2  \to  M_1\otimes M_2  $$ is contained in $F^{\geq \la} (M_1\otimes M_2 )$.  This is done by descending induction on $\la = \la_1+ \la_2$: we assume that the image of the morphism 
$$(F^{> \la_1}  M_1  \otimes F^{\geq \la_2}  M_2) \oplus  (F^{\geq \la_1}  M_1  \otimes  F^{> \la_2}  M_2)  \to   M_1\otimes M_2  $$ is contained in $F^{> \la} (M_1\otimes M_2 )$. 
 It follows that the composed morphism
$$\bar f: F^{\geq \la_1}  M_1  \otimes F^{\geq \la_2}  M_2  \to  M_1\otimes M_2   \to M_1\otimes M_2/ F^{\geq \la} (M_1\otimes M_2 )$$ factors through the object $\gr^{\la_1}
 M_1 \otimes \gr^{\la_2}
 M_2$. Since the latter is semistable of slope $\la $ by assumption, and since $F^{\geq \la}(M_1\otimes M_2/ F^{\geq \la} (M_1\otimes M_2 ))=0$ (\cf \ref{P3} (3)), one has $\bar f=0$. Therefore , for every $\la_1, \la_2$, one has
 \begin{equation}\label{eq3.1}F^{\geq \la_1}  M_1  \otimes F^{\geq \la_2}  M_2 \subset F^{\geq \la_1+\la_2}  (M_1\otimes M_2) ,\end{equation}
and similarly for $\geq$ replaced by $>$.
 Whence a canonical morphism 
\begin{equation}\label{eqIV}\tilde\gr_{M_1,M_2}:\, \gr M_1 \otimes \gr M_2 \to \gr (M_1\otimes M_2). 
\end{equation}
This makes $(\gr, \tilde\gr, \gr_\un = 1_\un)$ into a pseudo-$\otimes$-functor from $\sC$ to $\sC$ (\ie it satisfies all the axioms of a $\otimes$-functor, except that \eqref{eqIV} may not be an isomorphism a priori, \cf appendix), and \eqref{eqIII} 
$$ \gr\, (M^\vee) \overset{\cong}{\to} (\gr\, M)^\vee$$
(which holds by item (2)) is the canonical morphism $\hat\gr_M$ corresponding to 
$$ ev_{\gr M}\circ  \tilde\gr_{M^\vee, M}:\; \gr M \otimes (\gr M)^\vee \to \gr(M\otimes M^\vee)\to \un.$$
To check these assertions directly may be tedious, but they become clear if one considers the Rees deformation from $\gr M$ to $M$ (\cf \ref{Rees} and \ref{remqt}).

According to the corollary in the appendix, the fact that $\hat\gr_M$ is an isomorphism for any $M$  implies that $\tilde \gr_{M_1,M_2}$ is an isomorphism for any $(M_1,M_2)$, \ie {\it $\gr$ is a $\otimes$-functor}.

\smallskip  $(1a)+ (2a)\Rightarrow (1b)$. 
  One has a morphism of (horizontal) short exact sequences (which is functorial in $M_1, M_2$)
$$\begin{matrix}   \oplus  F^{> \la_1}  M_1  \otimes F^{\geq \la_2}  M_2 &\to & \oplus F^{\geq \la_1}  M_1  \otimes F^{\geq \la_2}  M_2 &\to & \oplus \gr^{  \la_1}  M_1  \otimes \gr^{ \la_2}  M_2  \\   \downarrow && \downarrow &&\downarrow   \\     F^{> \la} (M_1\otimes M_2 ) &\to &  F^{\geq \la} (M_1\otimes M_2 ) &\to & \gr^\la (M_1\otimes M_2)    \end{matrix}$$
 in which the third vertical morphism is an isomorphism. By (ascending or descending) induction, one gets \eqref{MULT}. 
 
 \smallskip
  $(1b)\Rightarrow (1c) \Rightarrow (1a)$ are immediate.

\smallskip  
 (2c)  By \ref{STABAB1}, we know that $\sC(0)$ is abelian. On the other hand it is stable under $\otimes$ and $()^\vee$. Therefore it is tannakian and the restriction of $\omega$ to $\sC(0)$ is a fiber functor. 
 
Since ${\rm char}\, F=0$, the socle $\sC(0)_{ssi}$ is a (semisimple) tannakian subcategory of $\sC(0)$. Therefore,   if $M_1$ and $M_2$ are stable of slope $0$, \ie simple objects of $\sC(0) $, then  $M_1\otimes M_2$ is a semisimple object of $\sC(0)$, \ie  a direct sum of stable objects of slope $0$.  
 \qed
 
 \begin{rem} When $\La= \Q$,  one can shorten the proof of $(1a)\Rightarrow (1b)$  by avoiding the devissage via $\gr$, on replacing $\gr$ by the Rees deformation functor
 $$R: \sC \cong Proj_A \to \{\text{\rm flat families of projective $A$-modules over ${\mathbb A}^1$}\}$$
 and applying the last corollary of the appendix to this pseudo-monoidal functor.
 \end{rem}

  \begin{prop} Assume that $\sC$ is abelian. For any exact $\otimes$-multiplicative slope filtration  on $\sC$, ${\omega}\circ \gr$ is a fiber functor.  If $\sC$  admits a $\otimes$-generator and $\La= \Q$, then ${\omega}\circ \gr \cong {\omega}$ (as fiber functors). 
   \end{prop}
 
  \proof For an exact slope filtration, $gr$ is exact, whence the first assertion. For the second, see \cite[2.2.5, 2.4]{Saa} (in \loccit  only filtrations indexed by $\Z$ are considered; in the case where $\sC$ is algebraic  and $\La= \Q$, the abelian group generated by all breaks is of the form $\frac{1}{N}\Z$, so that one may reduce to the case of filtrations indexed by $\Z$). \qed

  \medskip\subsection{Is any determinantal slope filtration $\otimes$-multiplicative?}
 
  Any $\otimes$-multiplicative slope filtration is determinantal (item (2) of \ref{T3}). The converse is an  interesting {\it open problem} (for $\La = \Q$, say).  
 
\smallskip In fact, there is many an instance in the literature, where the proof of $\otimes$-multiplicativity 
of a concrete determinantal slope filtration is either difficult or ad hoc. 
 It would therefore be desirable to know whether this is a general fact. 

\smallskip We propose a partial result in this direction, 
 assuming that $F$ is algebraically closed, and that $\La=\Q$.
 
Let $\sC'$ be the full subcategory of  $\sC$ consisting of direct sums of semistable objects.  Let us introduce a tensor product on $\sC'$ by setting
 $$M\hat\otimes N = \gr(   M\otimes  N).$$
 
 \begin{prop} The determinantal slope filtration $F^{\geq .}$  is $\otimes$-multiplicative if and only if $\hat\otimes$ is associative in the sense that for any three objects, $(M_1\hat\otimes M_2)  \hat\otimes M_3\cong  M_1\hat\otimes ( M_2   \hat\otimes M_3)$.
 \end{prop} 
 
 \proof  If $F^{\geq .}$  is $\otimes$-multiplicative, then $M\hat\otimes N =   M\otimes  N $ on $\sC'$ and the associativity follows.  
 
For the converse, let $M_1,M_2$ be semistable objects of slopes $\la_1$ and $\la_2$ respectively. We have to show that $M_1\otimes M_2$ is semistable of slope $\la_1+\la_2$; equivalently, that $M_1\hat\otimes M_2 = \gr^{\la_1+\la_2}M_1\hat\otimes M_2$. 
  Let $\sD$ be the smallest strictly full subcategory of $\bar \sC$ containing $M_1 $ and $M_2$, and stable under $\hat\otimes$, duality, sums and direct summands. This is a semisimple abelian subcategory of $\bar\sC$, and every object of $\sD$ is a direct sum of semistable objects.
 
Replacing $\sD$ by an equivalent small category $\sD_{sk}$ with $Sk\, \sD$ as set of objects,  we may assume, using Schur's lemma, that the associativity property of $\hat\otimes$ gives rise to a (functorial) associativity constraint on $  \sD_{sk}$. Then $ \sD_{sk}$ is a semisimple tannakian category generated by the classes of $M_1$ and $M_2$ (with respect to the tensor product $\hat\otimes$ and duality ${}^\vee$). 

Let $G$ be the associated reductive tannakian group over $F$. Then the subgroup $Pic(\sD)=Pic(\sD_{sk})$ of $Pic(\sC)=Pic(\bar\sC)$ is the character group $X(G)$ of $G$.

Let us consider the $\Q$-vector space $X(G)_\Q  \subset X(G^0)_\Q = X(Z(G^0))_\Q $ (where the superscript ${}^0$ stands for the connected component of identity, and $Z$ stands for the center). 
Note that the finite group $\pi_0(G)$ acts on $Z(G^0)$ on one hand, and on $X(G^0)$ on the other hand, and that 
$$X(G)_\Q = X(G^0)_\Q^{\pi_0(G)}  = X(Z(G^0))_\Q^{\pi_0(G)} = X(Z(G))_\Q. $$ 
In particular, there is a central cocharacter $y: {{\mathbb G}_m} \to G$ and an element $r\in \Q^\times$ such that for any $L \in  Pic(\sD)$, identified with a character $\chi_L$ of $G$, one has
$\,y \circ \chi_L = r \cdot \mu(L) \in Hom({\mathbb G}_m, {\mathbb G}_m)_\Q= \Q.\,$
 Up to scaling the slopes by the factor $r$, the decomposition of objects of $\sD_{sk}$ (viewed as representations of $G$) according to the action of $y$ amounts to the decomposition into semistable direct factors according to the slopes. 
 In particular, $y$ acts diagonally on $M_1$ and $M_2$ with respective weights $r\la_1, r\la_2$. Therefore it acts diagonally on $M_1\hat\otimes M_2$ with weight $r(\la_1+\la_2)$, \ie $M_1\otimes M_2$ is semistable of slope $\la_1+\la_2$.
 \qed
 
\begin{rem} In particular, if $\gr$ is identity on objects (a case which occurs for  the Harder-Narasimhan filtration of vector bundles on smooth projective curves of genus $\leq 1$), the slope filtration is $\otimes$-multiplicative.
\end{rem}

 \bigskip
 \section{$\otimes$-bounded slope filtrations.}  
 
  \subsection{Definition and characterization.}\label{bsf}
  
   \begin{defn} $F^{\geq .}$ is \emph{$\otimes$-bounded} if it satisfies the following conditions:
 \begin{enumerate}
\item $\mu(\un)=0$,
\item If $M_1$ and $M_2$ are semistable of slopes $\leq  \la$, the breaks of $M_1\otimes M_2$ are $\leq \la$,
 \item  if $M$ is semistable of slope $\la$, so is its dual $M^\vee$.
 \end{enumerate}
   \end{defn}
 
 \begin{exs} The Turrittin-Levelt and Hasse-Arf filtrations are $\otimes$-bounded. 
    \end{exs}
 
    \begin{thm}\label{T4}  \begin{enumerate}
  \item  There is equivalence between
 \begin{enumerate}  
 \smallskip \item $F^{\geq .}$ is $\otimes$-bounded,
 
 \smallskip\item $F^{\geq .}$  satisfies \begin{itemize} \item $\mu(\un)=0$,
\item for any non-zero $M_1, M_2$, the breaks of $M_1\otimes M_2$ are bounded from above by the maximum of the breaks of $M_1$ and $M_2$,
 \item  for any rank one object $L$, $\mu(L)=\mu(L^\vee)$,\end{itemize}

 \smallskip\item for any $\la $, the full subcategory $\,\sC(<\la)\,$ of $\,\sC\,$ consisting of objects $N$ with $F^{\geq \la} M=0$ is stable under $\otimes$ and ${}^\vee$, and contains $\un$ if and only if $\la>0$,
  
 \smallskip\item for any $\la $, the full subcategory $\,\sC(\leq\la)\,$ of $\,\sC\,$ consisting of objects $N$ with $F^{> \la} M=0$ is stable under $\otimes$ and ${}^\vee$,  and contains $\un$ if and only if $\la \geq 0$,
 
   \smallskip \item   \begin{itemize}
 \smallskip\item $F^{\geq 0}\un =\un, \,  F^{>0} {\bf 1} = 0$, and  for any $\la  $,
  \smallskip\item $\displaystyle F^{> \la} M_1 = F^{> \la} M_2 = 0 \Rightarrow F^{> \la} (M_1\otimes M_2 )=0$,
 \smallskip\item $\displaystyle F^{> \la} M  = 0\Rightarrow F^{> \la} M^\vee  = 0$.
 \end{itemize}  \end{enumerate}

\medskip\item Any $\otimes$-bounded slope filtration $F^{\geq .}$ also satisfies: 
   \begin{enumerate}  
 \smallskip \item $F^{\geq .}$ \emph{splits canonically}, \ie there is a canonical isomorphism of functors $\gr \cong id$. In particular, $F^{\geq .}$ is strongly exact (\ref{ssex}),
 \smallskip\item for any non-zero $M$, the breaks of $M$ are the breaks of $M^\vee$ and $\mu(M)=  \mu(M^\vee)$,
 \smallskip\item all breaks are \emph{non-negative},
  \smallskip\item the \emph{stable} objects are the \emph{simple} objects,
    \smallskip\item any subobject (\resp quotient) in $\sC$ of an object of $\sC(<\la)$ is in $\sC(<\la)$; same for $\sC(\leq \la)$, 
       \smallskip\item if $\sC$ is tannakian, so are $\sC(<\la)$ for $\la> 0$ and  $\sC(\leq \la)$ for $\la\geq 0$.
      \end{enumerate} \end{enumerate}
  \end{thm} 
   
   \begin{rem} In the case where $\sC$ is abelian, it follows from \ref{efsf} and items (1e) and (2a) that $\otimes$-bounded slope filtrations are exactly the slope filtrations discussed in \cite{A3}; indeed, the latter were defined to be exact filtrations that satisfy (1e).  
  The Hasse-Arf slope filtrations of \loccit are exactly the integral $\otimes$-bounded slope filtrations of the present paper.   
\end{rem} 
  
  We need a lemma.
  
  \begin{lemma} Assume that if $M$ is semistable of slope $\la$, so is $\,M^\vee$. Then
  \begin{enumerate}
\item for any non-zero object $N$, $\mu(N)= \mu(N^\vee)$,
\item the breaks of any non-zero $N$ are the breaks of $N^\vee$,
 \item $F^{\geq .}$ is strongly exact, 
\item $F^{\geq .}$ splits canonically. \end{enumerate}
  \end{lemma}

  \proof (1) is equivalent to $\deg N= \deg N^\vee$ and follows, by induction on the rank, from addivity of $\deg$ applied to the short exact sequence $0\to M\to N\to P\to 0$, where $M$ is the universal destabilizing subobject, and to its dual.
  
  \smallskip (2). The equivalence of categories $\sC \overset{^\vee}{\to} \sC^{op}$ sends the slope function $\mu$ to itself (by the previous item). It is then clear that the image by $\vee$ of $F^{\geq .}$ is the (unique) slope filtration $\hat F^{\geq .}$ on $\sC^{op}$ attached to $\mu$, and is given by 
  $$\hat F^{\geq \la}M= (F^{\geq \la } M^\vee)^\vee.$$ In particular,  $\hat gr^{ \la}M= (gr^{  \la } M^\vee)^\vee $, and it follows that $NP(M)=NP(M^\vee)$.    
  
 \smallskip (3). Let us show that for any subobject $M$ of a semistable object $N$, $\mu(M)= \mu(N)$ (\cf \ref{T2}). Indeed,  $\mu(M)\leq \mu(N)$ (by semistability), and since $M^\vee$ is a quotient of $N^\vee$, which is semistable of slope $\mu(N)$, $\mu(M^\vee)\geq \mu(N)$. Thus $\mu(M)= \mu(M^\vee)=\mu(N)$, which shows that $F^{\geq .}$ is strongly exact.

\smallskip (4). It suffices to construct, for any $\la$, a canonical right inverse $\iota^\la$ to the natural transformation $F^{\geq \la}\to \gr^\la$: indeed, $\iota^\la$'s composed with the natural transformations $F^{\geq \la} \to id$ will sum up to an isomorphism $\gr \to id$. 

We first construct $\iota^\la$ on objects. By descending induction, on may assume that $\la$ is the highest break of $M$, hence also of $M^\vee$ by the previous item. Since the filtration is exact (by the previous item) and $(\gr^\la M^\vee)^\vee $ is semistable of slope $\la$ (by assumption), the strict epi $$ M = M^{\vee \vee} \to  (\gr^\la M^\vee)^\vee $$ gives rise to a strict epi $$\gr^\la M  \to (\gr^\la M^\vee)^\vee .$$  In particular, $\rk   \gr^\la\, (M^\vee)\leq \rk \gr^\la\, M$, and in fact  $\rk   \gr^\la\, (M^\vee)= \rk \gr^\la\, M$ by exchanging $M$ and $M^\vee$. It follows that 
$ \gr^\la M \to  (\gr^\la M^\vee)^\vee  $ is actually an isomorphism. Composing $M   \to  (\gr^\la M^\vee)^\vee $ with the inverse of this isomorphism gives $\iota^\la$. 

The functoriality of $\iota^\la$ with respect to morphisms $M\overset{f}{\to} N$ is also established by descending induction on $\la$, the case when $\la$ is larger or equal to the breaks of $M$ and $N$ being clear.  
 \qed

\medskip \noindent {\it Proof of \ref{T4}}.  $(1a)\Rightarrow (1b)$ is immediate.

 \smallskip $(1b)\Rightarrow (1e)$. Note that
 $$\mu(\un)=0\Leftrightarrow F^{\geq 0}\un =\un, \,  F^{>0} {\bf 1} = 0.$$ Note also that the condition 
\begin{itemize} \item  the breaks of $M_1\otimes M_2$ are bounded from above by the maximum of the breaks of $M_1$ and $M_2$ 
\end{itemize} is equivalent to 
\begin{itemize} \item $\forall \la, \displaystyle F^{> \la} M_1 = F^{> \la} M_2 = 0 \Rightarrow F^{> \la} (M_1\otimes M_2 )=0,$ 
\end{itemize}
and the condition 
\begin{itemize} \item  the breaks of $M^\vee$ are bounded from above by the maximum of the breaks of $M$ 
\end{itemize} is equivalent to 
\begin{itemize} \item $\forall \la,\displaystyle F^{> \la} M   = 0 \Rightarrow F^{> \la} M^\vee =0$.
\end{itemize}
 It thus suffices to see that the latter condition follows from the special case of rank one objects and from the former condition. This follows from the fact that  $M^\vee $ is a direct summand of $\, M^{\otimes {(\rk M-1})} \otimes (\det M)^\vee\,$ (\cf  \eqref{dualtens}).

 \smallskip $(1e)\Rightarrow (1d)$ is immediate.

\smallskip $(1d)\Leftrightarrow (1c)$ since $$\,\sC(<\la)= \bigcup_{\la'<\la} \, \sC(\leq\la'),\,\sC(\leq\la)= \bigcap_{\la'>\la} \, \sC(<\la').$$

\smallskip  $(1d)\Rightarrow (1a)$. The condition that $\sC(\leq\la)$ contains $\un$ if and only if $\la \geq 0$ means that $F^{\geq 0}\un =\un,\; F^{>0} \un = 0$ (by left continuity of $F^{\geq .}$).

Let us choose for $\la $ the maximum of the breaks of $M_1$ and $M_2$. Then $M_1,M_2\in \sC(\leq \la)$, hence $M_1\otimes M_2\in \sC(\leq \la)$, which means that the breaks of  $M_1\otimes M_2$ are $\leq \la$. 

Similarly, the highest break $\rho(M^\vee)$ of $M^\vee$ is bounded by, hence equal to by symmetry, the highest break $\rho(M)$ of $M$. Assume that $M$ is semistable of slope $\la$ (\ie $\rho(M)=\mu(M)=\la$), and let $L$ be the universal destabilizing object of $M^\vee$. Then 
$ \mu(L^\vee)= \rho(M^\vee) = \rho(M)= \la$. On the other hand, $L^\vee$ is a quotient of the semistable object $M$, hence is zero or semistable of slope $\la$. It follows that $\Ker  (M\to L^\vee)= (M^\vee/L)^\vee$ is zero or semistable of slope $\la$. Since $\rho(M^\vee/L)<\mu(M^\vee) $ by definition of $L$, we have $L=M^\vee$, \ie $M^\vee$ is semistable of slope $\la$. 

 \medskip (2a) and (2b) follow from the lemma.

\smallskip  (2c).   Let $M$ be semistable of slope $\la$. Then the breaks of $M\otimes M^\vee$ are bounded from above by $\la$. On the other hand,  $\un $ is a direct summand of $M\otimes M^\vee$ (the coevaluation divided by $\rk M$ is a section of the evaluation morphism $M\otimes M^\vee \to \un$). Therefore $\la\geq 0$.

\smallskip  (2d) follows from item (1e) and \ref{SSS}.

\smallskip (2e) follows from the fact that $F^{\geq .}$ is strongly exact and split, and (1c) + (1d).

\smallskip (2f) follows from (2e) and (1c) + (1d).
    \qed

 \begin{prop}\label{SSIM} The highest break function attached to a $\otimes$-bounded slope filtration satisfies
 \begin{itemize}
 \item $\rho(\un)=0$
 \item $\rho(M_1\otimes M_2)\leq\rho(M_1\oplus M_2)= \max (\rho(M_1),\rho(  M_2))$
  \item for any rank one object $L$, $\rho(L)=\rho(L^\vee)$.
 \end{itemize}
 Conversely, if $\sC$ is abelian semisimple, any function $\rho: \, Sk \, \sC \setminus \{0\} \to \La$ which satisfies these conditions is the highest break function of a unique slope filtration  on $\sC$, which is $\otimes$-bounded.
      \end{prop}

\proof The conditions on $\rho$ are clear from item (1b) of the theorem. Conversely, if $\sC$ is abelian semisimple, it is clear that the split slope filtration defined by $\displaystyle{ gr^\la M = \oplus_{\rho(M_i)=\la}\, M_i}$, where $M= \oplus M_i$ is the isotypical decomposition, is the unique slope filtration  on $\sC$ with highest break function $\rho$, and that it is $\otimes$-bounded. \qed

\begin{prop}\label{extbou} Assume that $\sC$ is abelian, and let $\sC_{ssi}$ be its socle (the full subcategory of semisimple objects). Let  $F_{ssi}^{\geq .}$ be a $\otimes$-bounded slope filtration on $\sC_{ssi}$. Then its unique extension $F^{\geq .}$ to $\sC$ (\cf \ref{OBABATS}) is $\otimes$-bounded if and only if $M= \gr M$ for any $M\in \sC$.

 In that case, one has $\rho(M)= \rho(M_{ssi})$ for any object $M$ of $\sC$ and  its semisimplication $M_{ssi}$ in $Sk \,\sC_{ssi}$).
  \end{prop}

\proof   Indeed, this condition is necessary since any $\otimes$-bounded filtration is split. To prove sufficiency, it suffices (using the fact that $(M_1\otimes M_2)_{ssi}= (M_1)_{ssi}\otimes (M_2)_{ssi}$) to prove that $\rho(M)= \rho(M_{ssi})$ if $M=\gr M$. Actually $M=\oplus  \gr_\la M$ implies that the breaks of $M$ are the breaks of $M_{ssi} = \oplus (\gr_\la M)_{ssi}$, since $\mu((\gr_\la M)_{ssi})= \mu( \gr_\la M )=\la$. \qed

 \begin{prop}\label{MAX} The cone of $\otimes$-bounded slope filtrations is stable under the operation $(\mu_1, \mu_2) \mapsto \mu = \max(\mu_1, \mu_2)$ of slope functions.
   \end{prop}
   
   \proof  Since $\otimes$-bounded slope filtrations are split, hence strongly exact, any object $M$ has a canonical decomposition 
   $$M = \oplus \gr_{12}^{\la_1\la_2} M, \;\; \gr_{12}^{\la_1\la_2} M : = \gr_1^{\la_1}\gr_2^{\la_2}M= \gr_2^{\la_2}\gr_1^{\la_1}M.$$ Each summand $\gr_{12}^{\la_1\la_2} M$ is $\mu$-semistable of slope $\max(\la_1, \la_2)$. The statement then follows from characterization (1b) of $\otimes$-bounded slope filtrations.
    \qed

 \bigskip
 \subsection{The tannakian case.}  In this subsection, $\sC$ is an essentially small tannakian category over $F$, with a fiber functor ${\omega}: \sC \to Vec_{F'}$. Let $$  G = Aut^\otimes {\omega}$$ be the corresponding (tannakian) affine group scheme over $F'$.  Any slope filtration induces a separated, exhaustive, left continuous decreasing filtration of ${\omega}$ by $F$-linear subfunctors. 

\begin{prop}\label{p5}   Assume $F= F'$.  A  $\otimes$-bounded slope filtration on $\sC$ is equivalent to the data of a separated left continuous decreasing filtration $(G^{(\la)})_{\la\in \La_{>0}}$ of $G$ by closed normal subgroups satisfying the following condition:
 
\smallskip for any $M$ and any $\la>0$, the trivial subrepresentation $\omega(M)^{G^{(\la)} }$ is a direct summand of $\omega(M) $ (as representations of $G^{(\la)}$). 
 
\smallskip \noindent The correspondence is determined by the formula
$$\omega(F^{\geq \la}M) = \Ker  (\omega(M)\to \omega(M)_{G^{(\la)}}).$$
The quotient $G/G^{(\la)}$ is the tannakian group of the tannakian subcategory $\sC(<\la)$ of $\sC$.
  \end{prop}  
 
 \proof \cf \cite[1.2.3]{A2} (in \loccit  only the case $\La \subset \R$ is considered, but this restriction is unnecessary).\qed   

  See \loccit for a detailed study of integral $\otimes$-bounded slope filtrations.

  \begin{rem} The operation $(\mu_1, \mu_2)\mapsto \mu = \max(\mu_1, \mu_2)$ of \ref{MAX} corresponds, at the level of tannakian groups, to  $(G_1^{(\la)},  G_2^{(\la)})\mapsto G^{(\la)}= G_1^{(\la)}\cdot G_2^{(\la)}$ (the closed normal subgroup of $G$ generated by $G_1^{(\la)}$ and $ G_2^{(\la)}$).
\end{rem}

 \subsection{}\label{REMFF} We end this section with a special case of \cite[5.3.1]{A2}, in a setting reminiscent of remark \ref{REMF}.

We consider the poset of positive integers $n$  with respect to divisibility (which we also identify with the poset of open subgroups of $\hat \Z$).

We consider the following data:

\begin{itemize}
\item  for any $n$, a \emph{tannakian} category $\sC_n$ over an algebraically closed field $F$ (of characteristic zero) and a fiber functor $$\omega_n: \sC_n  \to Vec_F.$$  Let $G_n$ be the tannakian group of $\sC_n$;

  \item  a group-scheme epimorphism $$G_n \to n\hat\Z$$ such that for any multiple $n'$ of $n$, $G_{n'}$ is the inverse image of $n'\hat\Z\subset n\hat\Z$ in $G_n$. To $G_{n'} \inj G_n$ corresponds a faithful exact $\otimes$-functor $$\iota^\ast_{n,n'}: \sC_n \to \sC_{n'} $$ such that\footnote{the condition that $G_{n'}$ is the inverse image of $n'\hat\Z\subset n\hat\Z$ in $G_n$ amounts to saying that $\iota^\ast_{n,n'}Rep_F n'\hat\Z =  Rep_F n\hat\Z$ and any object of $\sC_{n'}$ is a subquotient of an object in $\iota^\ast_{n,n'}\sC_n$.} $$\omega_{n'} \circ \iota^\ast_{n,n'} = \omega_n ;$$

 \item an {\it integral $\otimes$-bounded} slope filtration on $\sC_n$, with the compatibility condition
 $$F^{\geq n'\la} \circ \iota^\ast_{n,n'} = \iota^\ast_{n,n'} \circ F^{\geq n\la};$$

 \item a $\otimes$-equivalence
 $\; \sC_n \to \sC_{n'} \,$ compatible with the slope filtrations.
 \end{itemize}

\begin{prop}\cite{A2} Let us assume moreover that for any $n$,
\begin{enumerate}
\item  characters of finite order of $G_n$ come from characters of $n \hat \Z$, and give rise to invertible objects of slope $0$,
 \item simple objects of $\sC_n$ of slope $0$ are invertible.
  \end{enumerate}
Then $G_n$ is an extension of $n\hat \Z $ by a connected prosolvable group. \qed
  \end{prop}

\begin{cor} For any object $M$ of $\sC_1$, there is a positive integer $n$ such that each graded direct summand $ \gr^\la \iota^\ast_{1,n}  M$ of $\iota^\ast_{1,n}  M$  is an \emph{iterated extension of invertible objects of slope $\la$}.

Moreover, if for any pair of non-isomorphic invertible objects $L, L'$ of slope $\la$, $L'\otimes L^\vee $ is of slope $\la$, then $ \gr^\la \iota^\ast_{1,n}  M$ is the tensor product of an invertible object of slope $\la$ by an iterated extension of $\un$ by itself.
\end{cor}

 \proof Since the image of $G_1$ in the representation $\omega( \gr^\la \iota^\ast_{1,n}  M)$ is connected solvable, this representation is triangulable by Kolchin's theorem. This justifies the first assertion.

 For the second assertion, notice that the assumption implies that there is no non-trivial extension between non-isomorphic invertible objects $L,L'$ of slope $\la$ (tensoring by $L^\vee$ and using the fact that the filtration is split).   \qed

  \begin{ex}  Let $\sC_n$ be the $\C$-tannakian category  of differential modules over $\C((x^{1/n}))$, together with its Turrittin-Levelt filtration relative to the variable $x^{1/n}$  (which is integral and $\otimes$-bounded). A fiber functor $\omega_n$ (with values in $Vec_\C$) may be constructed using Katz's canonical extensions \cite{Ka2}.  
   An obvious $\otimes$-equivalence
 $\; \sC_n \to \sC_{n'} \,$ compatible with the slope filtrations is given by substituting $x^{1/n}$ to $x^{1/n'}$. The $\otimes$-functor $\iota^\ast_{n,n'}$ corresponds to the pull-back
 $\Spec \C((x^{1/n'}))\to \Spec \C((x^{1/n}))$. 
 
 All the above conditions are satisfied (\cf \cite[5.3.3]{A2}).  The statement of the corollary, in this special case, is nothing but the Turrittin-Levelt theorem. 
  \end{ex}

     \newpage

 \addtocontents{toc}{{\bf III. A catalogue of determinantal slope filtrations.}\hfill\thepage}
\
\bigskip
\begin{center}
\large\bf III. A catalogue of determinantal slope filtrations.
\end{center}

\bigskip According to theorem \ref{DETFI},  given a quasi-tannakian category $\sC$, one can associate to any a homomorphism  $$\delta:  {\rm Pic }\, \sC \to \La  $$ which satisfies $\,\delta([L])\geq 0 $ whenever there is a non-zero morphism $\un \to L$, a unique (descending) slope filtration $\, F_\mu^{\geq .} \,$ on $\sC$ with slope function $$\mu(M)= \frac{\delta(\det M)}{\rk M}.$$

\medskip In this chapter, we review some examples of $(\sC, \delta)$.   Our point is that there is no need to provide an existence proof of the slope filtration in each case: all cases are covered at once by theorem \ref{DETFI}. Similarly, the fact that the subcategory $\sC(\la)$ of semistable objects of slope $\la$ is abelian if $\,\delta([L])> 0 $ whenever there is a non-zero non-iso morphism $\un \to L$ follows from the general result \ref{STABAB1}.   

We also discuss the structure of the semistable objects, and the relations between a few of these examples.

\bigskip  \section{Vector bundles and filtrations of Harder-Narasimhan type.}\label{HNf}   
  
   \subsection{Vector bundles on curves.}\label{vbc}  
  Let  $Vec_X$ be the quasi-abelian category of vector bundles over a smooth geometrically connected projective curve $X$ (defined over some field $F$).   The function 
   $$\delta ([L])  = \deg L\, \in \Z$$ on  $Pic(Vec_X) = Pic \,X$ gives rise to the classical {\it Harder-Narasimhan  filtration} on $Coh^{tf}_X$, indexed by $\La= \Q$.  It is integral and non-exact.
   
   If ${\rm car}\, F= 0$, $Vec_X$ is quasi-tannakian, and the Harder-Narasimhan filtration is $\otimes$-multiplicative \cite{NS}\cite{RR}\cite{Maru}\cite{Mi} (the shortest proof is in \cite{A6}).

  \subsection{Vector bundles on higher dimensional polarized varieties.}\label{vbhd} Let $Coh^{tf}_X$ be the quasi-abelian category of torsionfree coherent sheaves on a normal geometrically connected projective variety $X$ of dimension $d\geq 1$ defined over a field $F$.  Let $\sO(1)$ be an ample line bundle on $X$.

 The function  $$ \deg_{\sO(1)}  M : = (c_1(M)\cdot c_1(\sO(1))^{d-1})\, \in \Z$$ is a degree function on  $Coh^{tf}_X$. By \ref{T1}, it gives rise to a (unique) slope filtration on $Coh^{tf}_X$, indexed by $\La= \Q$, the {\it Harder-Narasimhan filtration} (for an analytic viewpoint on this filtration  in terms of Hermite-Einstein metrics, when $F=\C$, see \cite{Bra}).
 
  This filtration induces a slope filtration on the full subcategory $Coh^{refl}_X$ of reflexive coherent sheaves, which is also quasi-abelian. In particular, on a smooth surface, one has $Coh^{refl}_X= Vec_X$ and the Harder-Narasimhan filtration of a vector bundle is a filtration by sub-bundles. 
 
  However, $Coh^{tf}_X$ (\resp $Coh^{refl}_X$) is not quasi-tannakian over $F$ with respect to its natural $\otimes$ if $d>1$ (\resp $d>2$), since it contains non-locally free sheaves (which do not have duals). 
   To remedy this, one may consider, as in \cite{Sh}, the localized category $\underline{Coh}^{tf}_X$ obtained from $Coh^{tf}_X$  by inverting morphisms which are isomorphisms outside a closed subset of $X$ of codimension $\geq 2$; this is a quasi-tannakian category over $F$, if ${\rm car}\, F= 0$.  
   
  One may also work with the localized category $\underline{Coh}^{refl}_X$ obtained from $Coh^{refl}_X$  by inverting morphisms which are isomorphisms outside a closed subset of $X$ of codimension $> 2$, which is also a quasi-tannakian category over $F$, if ${\rm car}\, F= 0$. 
     
    \smallskip The function  $\,\delta ([L])  = \deg_{\sO(1)}  L \,$  
   on  $Pic(\underline{Coh}^{tf}_X) = Pic(\underline{Coh}^{refl}_X) = Pic \,X$ satisfies the positivity condition, hence gives rise (by \ref{DETFI}) to a determinantal slope filtration on $\underline{Coh}^{tf}_X$ (\resp $\underline{Coh}^{refl}_X)$)  indexed by $\La= \Q$, which is nothing but the filtration induced by the Harder-Narasimhan filtration on $Coh^{tf}_X$.
   
     One way to prove its $\otimes$-multiplicativity is to reduce to the one-dimensional case by taking linear sections of $X$, \cf \cite{MR1}.

   \subsection{Vector bundles on compact analytic varieties.}\label{bru} 
  Let $Coh^{tf}_X$ be the quasi-abelian category of torsionfree coherent sheaves on a compact complex manifold $X$ of dimension $d\geq 1$. 
   
   In this context, the (missing) polarization is replaced by a Gauduchon metric $g$ on $X$, \ie a hermitian metric whose associated K\"ahler form $\omega_g$ satisfies $\partial \bar\partial  \omega_g^{d-1}= 0.$\footnote{Gauduchon proves that any hermitian metric is conformally equivalent to a unique (up to homothety) Gauduchon metric.}.  

The function  $$ \deg_g  M   : =\int_X c_1(L,h)\cdot  \omega_g^{d-1}  \in \R,$$ where $L$ denotes the double dual of $\displaystyle\bigwedge^{\rk M} M$ and $h$ is an auxiliary hermitian metric on this line bundle (the integral does not depend on its choice),  is a degree function on  $Coh^{tf}_X$ (\cf \cite{Bru}).
By \ref{T1},  it gives rise to a slope filtration on $Coh^{tf}_X$ indexed by $\Lambda= \R$. This filtration was introduced by Bruasse \cite{Bru}, building on results of Kobayashi \cite{Kob}\footnote{following a widespread belief, according to which the first step in constructing a filtration of Harder-Narasimhan should consist in proving that the set of degrees of subsheaves of $M$ has a maximum, Bruasse establishes this fact in $Coh^{tf}_X$ using deep compacity arguments. However, as we have seen in \ref{T1}, there is no need to prove this statement a priori: it is a formal consequence of the properties of a degree function.}
(whereas the corresponding notion of stability was introduced earlier by Toma \cite{Tom}). 

The {\it Bruasse filtration} induces a slope filtration on the full subcategory $Coh^{refl}_X$ of reflexive coherent sheaves, which is also quasi-abelian. In particular, on a compact complex surface, one has $Coh^{refl}_X= Vec_X$ and the {Bruasse slope filtration} of a vector bundle is a filtration by sub-bundles.

As in the algebraic case, one can pass to the (quasi-tannakian) localization $\underline{Coh}^{tf}_X$ in order to get a determinantal slope filtration, attached (via \ref{DETFI})  to the function
  $\,\delta ([L])  = \deg_g L \,$ on  $Pic(\underline{Coh}^{tf}_X) = Pic\, X$.   Its $\otimes$-multiplicativity is an open problem.

     \subsection{Higgs bundles.}\label{hig}  
Let $(X, \sO(1))$ be a polarized smooth geometrically connected projective variety over a field $F$ of characteristic $0$. According to Hitchin \cite{Hi} and Simpson \cite{Si}, a Higgs sheaf is a coherent sheaf $M$ together with a morphism $\theta: M \to M\otimes \Omega^1_X$ such that $\theta\wedge \theta=0$. Torsion-free Higgs sheaves form a quasi-abelian category $Higgs^{tf}_X$. The degree of the underlying coherent sheaf (with respect to the polarization) induces a degree function, hence a slope filtration, on $Higg^{tf}_X$ (and on the full quasi-abelian subcategory $Higgs^{refl}_X$ of reflexive objects). This filtration was studied in detail in \cite{F1} when $X$ is a curve.

Passing to the localization $\underline{Higgs}^{tf}_X$ (\resp $ \underline{Higgs}^{refl}_X$) with respect to morphisms which are isomorphisms outside a closed subset of $X$ of codimension $\geq 2$ (\resp $>2$), one gets a determinantal slope filtration, which is $\otimes$-multiplicative (\cf \cite[Cor. 3.8]{Si}).  

\medskip   {The Hitchin-Simpson correspondence} (\cf \cite{Si}) is a one-to-one correspondence between stable Higgs bundles with vanishing Chern classes and irreducible representations of $\pi_1(X(\C))$. Combined with the Riemann-Hilbert correspondence, this can be reformulated, if $X$ is a curve, as a $\otimes$-equivalence 
\begin{equation} ({Higgs}_X(0))_{ssi} \cong (DMod_X)_{ssi} \end{equation}
 between  direct sums of stable Higgs bundles of slope $0$ and semisimple vector bundles with connection, which generalizes the Narasimhan-Seshadri correspondence: ordinary vector bundles ($\theta=0$) corresponding to unitary connections. 

The fact that ${Higgs}_X(0)_{ssi}$ is a semisimple abelian category follows formally from \ref{STABAB1}.

   \bigskip
 \section{Arithmetic vector bundles and filtrations of Grayson-Stuhler type.}\label{AGSf}
 
     \subsection{Hermitian lattices.}\label{GSf} Let $Vec_{\sO_K}$ be the quasi-abelian category of projective modules of finite rank over the ring of integers of a number field $K$. A hermitian lattice $\bar M$  is an object $M$ of $Vec_{\sO_K}$ together with a hermitian norm $\vert\;\vert_v$ on $M\otimes_{\sO_K, v}\C$ with respect to each archimedean place $v$ of $K$;  for $K=\Q$, this is the same as an euclidean lattice, as considered in \ref{gsf}, \ref{exproto}. 
     
     Morphisms of hermitian lattices are $\sO_K$-linear maps of norm $\leq 1$ with respect to each $\vert\;\vert_v$. 
       Hermitian lattices form a proto-abelian category $Herm_{\sO_K}$ (\cf \ref{exproto2}). 
     
  \smallskip  The function 
     $$ \widehat{deg} \bar M =(\log \sharp (M/\sum_1^{\rk M} \sO_K s_i)- \frac{d_v}{2}\sum_v \log \det (\langle s_i, s_j\rangle_v))$$ where $s_i\, (i= 1,\ldots {\rk M})$ are elements of $M$ which form a basis of $M_K$, and $d_v$ is $1$ or $2$ according to whether $v$ is real or complex (the above expression is independent of this choice),  is a degree function on  $Herm_{\sO_K}$. By \ref{T1},  it gives rise to slope filtration on $Herm_{\sO_K}$, indexed by $\La= \Q$, the {\it Grayson-Stuhler filtration}.

    \smallskip   Although $Herm_{\sO_K}$ is non-additive, it is a rigid monoidal category, and it is possible to define the determinant of any object. With proper normalization, $ \widehat{deg} \bar M$ depends only on $\det \bar M$. In this sense, the Grayson-Stuhler slope filtration looks like a determinantal filtration. 
     
     Whether it has the $\otimes$-muliplicativity property is an open problem, already for $K=\Q$ (it was conjectured by J.-B. Bost \cite{Bos}; \cf \cite{DSP}, \cite{Che}, \cite{BK} and \cite{A6} for partial results).

       \subsection{Arithmetic vector bundles in Arakelov geometry.}\label{Mf} 
     A.  Moriwaki has generalized this filtration to the case of hermitian torsion-free sheaves $\bar M$ on a polarized normal arithmetic projective variety $X$ of any dimension $d$. In the case of an arithmetic surface, endowed with a nef and big hermitian line bundle $\bar H$, this is the filtration on the proto-abelian category of hermitian torsion-free sheaves $\bar M$ on $X$, attached to the degree function\footnote{the main effort in Moriwaki's paper is devoted to proving that the set of degrees of subsheaves of $\bar M$ has a maximum. Again, as we have seen in \ref{T1}, it should not be necessary to prove this statement a priori: it is a formal consequence of the properties of a degree function. 
       
       On the other hand, Chen \cite{Ch} gives another axiomatic viewpoint on these filtrations; however, the existence of the universal destabilizing subobject and a version of our lemme \ref{MSS} are taken by him as axioms.} given by $$\widehat{\deg}\, \bar M := \widehat{\deg} (\hat c_1(\bar M)\cdot \hat c_1(\bar H)^{d-1})\in \R.$$

  \subsection{Variants.}\label{HJSf} 
  In order to strengthen the analogy between the Harder-Narasimhan filtration and the Grayson-Stuhler filtration, Hoffman, Jahnel and Stuhler have extended the Harder-Narasimhan filtration to quasi-abelian category of adelic vector bundles on smooth algebraic curves \cite{HJS}\footnote{also considered by Gaudron \cite{Gau}, in characteristic $p$.}. Again, existence and unicity follow directly from \ref{DETFI}.
  

  \bigskip
  
  \section{$\phi$-modules and filtrations of Dieudonn\'e-Manin type.}
  
    \subsection{$\phi$-modules.}\label{phim} Let $R$ be field or a B\'ezout ring, and let $\phi$ be an injective endomorphism of $R$ such that the invariant ring $F= R^\phi$ is a field of characteristic zero. 
    
    A {\it $\phi$-module} is a free $R$-module of finite rank $M$  together with an isomorphism 
    $$\Phi: \, M\otimes_{R, \phi}R\to M.$$
    The category $\phi$-$Mod_R$ of $\phi$-modules, with its natural $\otimes$, is quasi-tannakian over $F$, and even tannakian if $R$ is a field.    
    
\smallskip    One has $Pic( \phi$-$Mod_R)=  R^\times/ \{b/\phi(b),\; b\in R^\times\}$.
If $L$ is represented by $c\in R^\times$, the existence of a non-zero morphism $\un \to L$ translates into the existence of $ a\in R$ (possibly non-invertible) such that $c= a/\phi(a)$. 

Thus, by \ref{DETFI}, any homomorphism $$\delta:   \, R^\times/ \{b/\phi(b),\; b\in R^\times\} \to \Z  $$ which satisfies 

\smallskip \centerline{$\,\delta([a/\phi(a)])\geq 0$ for any $a\in R$ such that $a/\phi(a)\in R^\times$}

\smallskip\noindent gives rise to a unique integral (descending) slope filtration with slope function $\,  \frac{\delta(\det M)}{\rk M}$.

    \medskip For instance, if $R$ is endowed with a valuation $v$ with values in $\Z$ and $\phi$ is an isometry with respect to $v$, then both choices
    
    \smallskip \centerline{$\,\delta([b])= + \,v(b)\;\;\forall b\in R^\times$}
    \centerline{ $\;\delta([b])= - \,v(b)\;\;\forall b\in R^\times$}
    
    \smallskip\noindent satisfy the assumption. 
      It turns out that the interesting examples occur with the  $-$ sign.    
      
      \medskip  Slightly differently, assume that $R^\times \cup \{0\}$ is a subring of $R$, and that  $v$ is a valuation on this subring with values in $\Z$, and that $\phi$ is an isometry with respect to $v$. Then  
      
        \smallskip \centerline{$\,\delta([b])= - \,v(b)\;\;\forall b\in R^\times$}
        
      \smallskip\noindent satisfies the assumption (it may occur that $\delta =+v$ does not, \cf \ref{Kf} below).
    
  \subsection{Description of the Newton polygon when $R=K$ is a complete valued field.}\label{phiNP}   
      
   Let $M$ be a cyclic $\phi$-module over a complete valued field $(K,v)$ (of characteristic $0$). Since the twisted polynomial ring $K\langle \phi\rangle$ is left principal, $M$ is of  the form $K\langle \phi\rangle/ K\langle \phi\rangle P$, with $P$ monic.

   Let us define the Newton polygon $NP(P)$ of $P= \sum a_i \phi^i$ to be the convex envelope of the lines $x= i, y\leq  - v(a_{n-i})$ (the origin is the left-end point). Then  
 it is known that $P $ admits a unique factorization 
   $$ P = P_{\la_r}\cdots P_{\la_1}$$ where 
   $\la_1> \cdots >\la_r$, $P_i$ is monic and 
$NP(P_{\la_i})$ has just one slope $\la_i$ (\cf \cite[14.2.5]{Ke4}). 
   
   From this, one derives that $NP(P)= NP(M)$ (for $\delta = -v$), the factorization of $P$ corresponding to the slope filtration of $M$ (\loccit 14.4.15). The filtration is $\otimes$-multiplicative (\loccit 14.4.9). Moreover, it is split if $\phi$ is invertible (\loccit 14.4.13).  
       
      \subsection{Frobenius modules.}\label{DMf} Let $R= K$ be a complete valued field of characteristic $0$, with residue field $k$ of characteristic $p>0$. Let $\phi$ be a lifting of some fixed positive power of the Frobenius endomorphism of $k$, so that $\phi$ is an isometric endomorphism of $K$.

In this context,  $\phi$-modules are also called $\rm F$-isocrystals (over the point), after Grothendieck.  The determinantal slope function attached to the slope function
$\,\mu(M) := -\frac{v(\Phi_{\det M})}{\rk M}\,$ 
is the (descending version\footnote{see \ref{asc} and \ref{asc2}.} of) the classical {\it Dieudonn\'e-Manin filtration}.

This slope filtration is $\otimes$-multiplicative. It is exact. Moreover, it is split if $\phi$ is invertible, \ie if $k$ is perfect.

         \subsection{$q$-difference modules.}\label{ASVf}
       Let $K=R$ be either the field $\C(\{x\})$ of germs of meromorphic functions at the origin, or its $x$-adic completion $\C((x))$, endowed the the $x$-adic valuation $v$. Let $q$ be a non-zero complex number, not a root of unity, and 
       let $\phi$ be the isometric continuous $\C$-automorphism of $K$ given by 
       $$\phi(x)=qx $$ 
       (here $F=\C$ since $q$ is not a root of unity). In this context,  $\phi$-modules are called $q$-difference modules\footnote{these objects occur in the context of $q$-calculus, which has a long history (Euler, Gauss, Jacobi, Heine, 
Ramanujan,
\dots), and is
based on the replacement of ordinary integers $n$ by their
$q$-analogs
  $$ [n]_q= 1+q+ q^2+ \dots  + q^{n-1}.$$
  The usual derivation $d/dx$ is then replaced by the $q$-derivation 
$$ d_q : \;\;\;   f(x) \mapsto
{f(x)-f(qx)\over (1-q)x}$$  which sends $x^n$ to $[n]_qx^{n-1}$.
 Differential equations are thus replaced by $q$-difference
equations, which are nothing but functional equations
$$\,y(q^nx )+ a_{n-1}y(q^{n-1}x)+ \ldots + a_0 y(x)= 0.$$ The ``confluence" of $q$-difference equations to differential equations occurs when $q$ tends to $1$.
   The analytic theory of $q$-difference equations is well-developed when $\vert q\vert \neq 1$, \cf \eg \cite{dVRSZ} for a survey; when $\vert q\vert = 1$, one encounters phenomena of small divisors which make the study more delicate.}.

       \smallskip  
       
 The slope function
$\,\mu(M) := -\frac{v(\Phi_{\det M})}{\rk M}\, $
 gives rise to an integral determinantal slope filtration on $\phi$-$Mod_{\C(\{x\})}$ and to an other one on $\phi$-$Mod_{\C((x))}$. 
 
\smallskip  The filtration on $\phi$-$Mod_{\C((x))}$  is split and $\otimes$-multiplicative.
 
\smallskip For $\vert q\vert \neq 1$, the filtration on $\phi$-$Mod_{\C(\{x\})}$ has been considered by Sauloy, who proved (using Adams's lemma \cite{Ad}) that it is induced by the filtration on $\phi$-$Mod_{\C((x))}$, \cf \cite{Sau1}\footnote{in their recent work (\cf \cite{RS}), Ramis and Sauloy changed the convention on the sign of slopes, working with an increasing filtration instead of a decreasing one (this is a mere convention and has nothing to do with the above choice of sign $\delta = \pm v$).}. The  $q$-difference modules of slope $0$ are well understood, \cf \cite{Sau2}.

\smallskip  Recently, the more delicate case $\vert q\vert = 1$ has been tackled by Di Vizio \cite{dV}. Again, under some diophantine conditions (on $q$ and on the so-called exponents), the filtration agrees with the formal one. 
 
  \begin{rem}  The {\it Adams-Sauloy filtration} is exact but not split in general.   For instance, let $(M,\Phi)$ be the  $q$-difference module given by the matrix
   $$ \begin{pmatrix}  1/x & 1/x \\ 0 & 1  \end{pmatrix}$$
   in the canonical basis of $\C(\{x\})^2$. It has breaks $1$ and $0$, and is indecomposable, \cf \cite[2.2.1]{Sau2}. 
  
  On the other hand, for the filtration given by $\delta = +v$ (instead of the Adams-Sauloy filtration given by $\delta= -v$), the same $q$-difference module is semi-stable of slope $1/2$. The tannakian category generated by $M$ is actually equivalent to the one of example \ref{EXX}. \end{rem}

   \medskip For $\vert q\vert > 1$,  $q$-difference modules over ${\C((x))}$ are closely to vector bundles on the elliptic curve $X= \C^\times/q^\Z$.  Let $\gr \,\phi$-$Mod_{\C(\{x\})}\subset \phi$-$Mod_{\C(\{x\})}$  be the tannakian full subcategory consisting of objects such that $M=\gr M$. One has a canonical fiber functor 
 
 \smallskip  \centerline{$\gr \,\phi$-$Mod_{\C(\{x\})}\to Vec_X$}
 
 \smallskip\noindent  which is essentially bijective, \ie induces an isomorphism on skeleta. This functor is compatible with the Adams-Sauloy filtration on the left-hand side and the Harder-Narasimhan filtration on the right-hand side, \cf \cite{PRe}\cite{Sau2}.

    \medskip   \subsection{$\phi$-modules on the Robba ring.}\label{Kf}
       
       Let now $R=\sR$ be the Robba ring over a $p$-adic field $K$, \ie the ring of $K$-holomorphic functions on some open annulus with outer boundary $1$ (such functions are represented by Laurent series with coefficients in $K$ and appropriate convergence conditions). This is a B\'ezout ring (Lazard). 
       
       Let $\sR^{bd}$ be the subring of bounded elements. This is actually a field, which is henselian with respect to the natural ($p$-adic) valuation (which extends in no way to $\sR$ itself), which we normalize to take values in $\Z$. 
        Moreover $\sR^\times \cup \{0\}= \sR^{bd}$.
        
     \smallskip   Let $\phi$ be an injective endomorphism $K$ given by 
       $$\phi(x)= x^{p^m} \;\; \text{\rm or}\;\;  (1+x)^{p^m}-1,$$ 
     and acting via some power of Frobenius on the coefficients (so that $F= K^\phi$). It preserves $\sR^{bd}$.  
             
 \smallskip The slope function
$\,\mu(M) := -\frac{v(\Phi_{\det M})}{\rk M}\,$ 
 gives rise to an integral determinantal slope filtration\footnote{One has $v(\phi(b))= v(b)$ for any $b\in \sR^\times$, but it may happen that for $a\in \sR \setminus \sR^\times$, $\phi(a)/a\in \sR^\times$ and $v(\phi(a)/a)>0$ (example: $a= \log(1+x)$). Hence, in the setting of \ref{phim}, one has to take $\delta= -v$,  not $+v$.}
 on  the quasi-tannakian category $\phi$-$Mod_\sR$. This filtration was introduced and studied by Kedlaya\footnote{actually, he works with the corresponding ascending slope filtration, \cf \ref{asc}.}. 
 
 He proved that any $\phi$-module $M$ of slope $0$ over $\sR$ comes from a (unique) $\phi$-module $M^{bd}$ over $\sR^{bd}$ such that $M^{bd}\otimes_{\sR^{bd}}\widehat{\sR^{bd}}$  is of slope $0$ with respect to the Dieudonn\'e-Manin filtration on $\phi$-$Mod_{\widehat{\sR^{bd}}}$ (a so-called unit-root $\rm F$-isocrystal), \cf \cite{Ke1}.

\smallskip From this, and the fact that the Dieudonn\'e-Manin filtration is $\otimes$-multiplicative, it follows that the {\it Kedlaya filtration} is $\otimes$-multiplicative as well. 
 Unlike the Dieudonn\'e-Manin filtration, however, it is not exact (as Colmez's theory of trianguline representations shows). 
 
 \smallskip Given a $\phi$-module over an ${\sR^{bd}}$, one can consider the Newton polygon of $M\otimes_{\sR^{bd}}\widehat{\sR^{bd}}$ with respect to the Dieudonn\'e-Manin filtration, and the Newton polygon of $M\otimes_{\sR^{bd}} {\sR }$ with respect to the Kedlaya filtration. With our conventions on Newton polygons (\ref{newp}), the former lies above the latter, with the same end-points \cite{Ke1}.

  \medskip
       \subsection{Local $\rm F$-isocrystals, and $(\phi, \Gamma)$-modules.}\label{phigam}
  
  Local $\rm F$-isocrystals are differential modules with Frobenius structure over the Robba ring. More precisely, they are free $\sR$-modules of finite rank which are simultaneously $K\langle x, d/dx\rangle$-modules and 
  $\phi$-modules, in a compatible way: $\Phi$ commutes with the action of $d/dx$.
  
  They form a tannakian category $\rm F$-$Isoc_\sR$ over $F= K^\phi$. The Kedlaya filtration on $\phi$-$Mod_\sR$ induces a $\otimes$-multiplicative slope filtration on this category, which is exact. 
  
  Using this filtration and the characterization of slope $0$ objects, Kedlaya proved the $p$-adic local monodromy theorem (Crew's conjecture) by reduction to the case of unit-root isocrystals on $\sR^{bd}$, which was treated by Tsuzuki, \cf \cite{Ke1}.
  
  \medskip The notion of $(\phi, \Gamma)$-module over $\sR$ is a variant of that of local $\rm F$-isocrystal,  
  which Fontaine introduced in the theory of $p$-adic representations of $p$-adic fields. For simplicity, we take $K=F=\Q_p$. Here $\phi(x)= (1+x)^p-1$.
  $\Gamma $ is the cyclotomic quotient $Gal(\Q_p(\zeta_{p^\infty})/\Q_p)$ of $G_{\Q_p}$ (isomorphic to $\Z_p^\times$ via the cyclotomic character $\chi$); it acts on $\sR$ via $\gamma(x)= (1+x)^{\chi(\gamma)} -1$. The infinitesimal generator of $Lie\, \Gamma$ can be identified with the derivation $ (1+x)\log(1+x)d/dx$\footnote{the factor $\log(1+x)$, which vanishes on the set $\zeta_{p^\infty}-1$, gives rise to  difficulties with ``apparent singularities".}.
  
\smallskip A $(\phi, \Gamma)$-module over $\sR$ (\resp $\sR^{bd}$) is a free $\sR$-module (\resp $\sR^{bd}$-module) of finite rank which is simultaneously a (semilinear continous) $\Gamma$-modules and a
  $\phi$-module, in a compatible way: $\Phi$ commutes with the action of $\Gamma$. 
 $(\phi, \Gamma)$-modules form a quasi-tannakian category  $(\phi, \Gamma)$-$Mod_{\sR}$ (\resp $(\phi, \Gamma)$-$Mod_{\sR^{bd}}$) over $\Q_p$, and the Kedlaya filtration on $\phi$-$Mod_\sR$ induces a $\otimes$-multiplicative slope filtration on $(\phi, \Gamma)$-$Mod_{\sR}$.
  
\medskip  According to Fontaine, Colmez and Cherbonnier \cite{Fo}\cite{CC}, there are $\otimes$-equivalences of tannakian categories: 

   \begin{equation}\label{FCC} Rep^{cont}_{\Q_p} G_{\Q_p} \cong  (\phi, \Gamma)\text{\rm -}Mod_{\sR}(0) \cong (\phi, \Gamma)\text{\rm -}Mod_{\sR^{bd}}(0) \end{equation}
    
    \smallskip\noindent   
where $(0)$ refers to the subcategory of objects of Kedlaya slope $0$ (the fact that $(\phi, \Gamma)$-$Mod_{\sR^{bd}}(0)$ is abelian can also be derived from \ref{STABAB1}).

  \bigskip
  \section{Filtered modules and filtrations of Faltings-Fontaine type.} 
   
     \subsection{Filtered modules.}\label{FRf}
     
    Let $K/F$ be a finite extension of fields of characteristic $0$, and let $n$ be a positive integer. Let  $n$-$Fil_{K/F}$ be the category of finite-dimensional $F$-vector spaces $V$  together with $n $ (separated, exhaustive, decreasing) $\Z$-filtrations
    $F^._\nu$  
    on $V\otimes_F K$.
    
    This is a quasi-tannakian category over $F$. The homomorphism
    
       \smallskip \centerline{$\delta : \;Pic (n$-$Fil_{K/F}) \cong  \Z^n \to \Z$,}
    
    \smallskip\noindent  given by the sum of the coordinates (the notches of the filtrations), gives rise to an integral determinantal slope filtration, which was studied by Faltings and Rapoport \cite{FW}\cite{F2}\cite{Ra} (it occurs in the theory of $p$-adic period mappings).  It is a non-exact (this is easily seen by considering a stable object of rank $>1$). 
    
      \medskip In \cite{FW}, Faltings and W\"{u}stholz relate it to the Harder-Narasimhan filtration, as follows.  Let $X$ be a cyclic covering of ${\bf P}^1$, totally ramified above $n[K:F]$ branch points, at least. To $(V, (F^._\nu))$, they associte a vector bundle $M(V, (F^._\nu))$ on $X$ of rank $\dim V$ and degree $[K:F]\delta(\det (V, (F^._\nu)))$. The construction commutes with $\otimes$. Moreover $M(V, (F^._\nu))$ is semistable if $(V, (F^._\nu))$ is, and conversely provided the degree of the covering $X/{\bf P}^1$ is large enough. The $\otimes$-multiplicativity of the {\it Faltings-Rapoport filtration} thus follows from the $\otimes$-multiplicativity of the Harder-Narasimhan filtration (for other approaches, \cf \cite{F2}\cite{To1}\cite{To2}).

        \subsection{Filtered $\phi$-modules}\label{FFf}  
        
      In the context of \ref{phim}, let $Fil$-$\phi$-$Mod_{R}$ be the category of $\Z$-filtered $\phi$-modules $(V, \Phi, F^.)$ over $R$ (no relation between  $\Phi$ and $F^. $ is imposed). This is a quasi-tannakian category over $F= R^\phi$.  

    It has two natural determinantal slope filtrations: the ``tautological" one induced by $F^.$; and the one given by $\delta = -v$. One can also consider their middle point, \ie the determinantal slope filtration defined by 
 
  \smallskip\centerline{$Pic(Fil$-$\phi$-$Mod_{R})\; \cong \; \Z \times R^\times/\{b/\phi(b),\;b\in R^\times\} $}
  
      \smallskip\noindent  given by $(n\in \Z, c\in R^\times) \mapsto n-v(c)\,$ ($n$ is the notch of the filtration).
    
       \medskip This ``middle filtration" is relevant in the context of \ref{phigam}, where it was considered by Fontaine and others.       
    According to Fontaine and Colmez \cite{CF} (\cf also \cite[V]{Be}), there is a $\otimes$-equivalence of tannakian categories: 

      \begin{equation}\label{CF} Rep^{crys}_{\Q_p} G_{\Q_p} \cong  Fil\text{\rm -}\phi\text{\rm -}Mod_{\Q_p}(0).\end{equation}

    \smallskip\noindent   
where $(0)$ refers to the subcategory of objects of Fontaine slope $0$ (the fact that $Fil$-$\phi$-$Mod_{\Q_p}(0)$ is abelian can also be derived from \ref{STABAB1}), and where the superscript crys refers to crystalline representations.

 \medskip  To close the circle, Berger has constructed a fully faithful $\otimes$-functor of quasi-tannakian categories
 
\begin{equation} Fil\text{\rm -}\phi\text{\rm -}Mod_{\Q_p} \inj  (\phi, \Gamma)\text{\rm -}Mod_{\sR}\end{equation}

    \smallskip\noindent and proven that it preserves the slope filtrations \cite{Be} (this is one of the ways to prove that the {\it Fontaine filtration} is $\otimes$-multiplicative). Via \eqref{FCC} and \eqref{CF}, the embedding of subcategories consisting of objects of slope $0$ (in the sense of Fontaine and Kedlaya, respectively) corresponds to the embedding
    $ \,Rep^{crys}_{\Q_p} G_{\Q_p}\inj Rep^{cont}_{\Q_p} G_{\Q_p}$.

\newpage
    
 \addtocontents{toc}{{\bf IV. A catalogue of $\otimes$-bounded slope filtrations.}\hfill\thepage}
\
\bigskip
\begin{center}
\large\bf IV. A catalogue of $\otimes$-bounded slope filtrations.
\end{center}

\bigskip To produce $\otimes$-bounded slope filtrations on a tannakian category $\sC$ (over a field $F$ of characteristic $0$) is not as easy as to produce determinantal slope filtrations. One way is by constructing a sequence of normal subgroups of the tannakian group as in \ref{p5}. Another way, when $\sC$ is semisimple, is by defining the highest break function and checking the simple conditions of \ref{SSIM}.

We recall that the breaks of a $\otimes$-bounded slope filtration are always non-negative.

\bigskip
  \section{Differential modules and filtrations of Turrittin-Levelt type.}
  
\medskip  \subsection{Formal differential modules in one variable.}\label{TLf} Let $F$ be an algebraically closed field of characteristic zero. Let the derivation $\partial = x\frac{d}{dx}$ acts on $K= {F}((x))$ and respects the $x$-adic valuation $v$ on $K\setminus {F}$.
  
   Let $DMod_{{K}}$ be the category of differential modules $M =(V, \nabla(\partial))$ over ${K}$. This is a tannakian category over ${F}$. 
  
 \smallskip The highest break function associated to the {\it Turrittin-Levelt filtration} is given by the Poincar\'e-Katz rank: 
 \begin{equation} \rho(M)= \max(0, -   v_{sp}(\nabla(\partial))\end{equation} 
involving  the spectral valuation (\cf \eg \cite[2.1]{A5})
 $$v_{sp}(\nabla(\partial)) = \lim \frac{1}{n} v(\nabla(\partial)^n). $$ 
  The conditions of \ref{SSIM} are easily checked using this definition, so that \ref{SSIM} shows that associated (Turrittin-Levelt) filtration is $\otimes$-bounded, as far as one considers semisimple differential modules. 

  To check that it is $\otimes$-bounded on $DMod_{{K}}$, one would have to show that $M= \gr M$ for any $M$ (\cf \ref{extbou}), but this splitting property of formal differential modules is non-trivial. 
  
  It is established by using the fact that any differential module over $K$ is cyclic, \ie of the form $K\langle \partial\rangle/ K\langle \partial\rangle P$, and showing that $P$ admits a unique factorization as in the case of $\phi$-modules (\cf \ref{phiNP}), Newton polygons being defined in a similar way, except that one considers only non-negative slopes.

  \medskip The degree function attached to the Turrittin-Levelt filtration is called the {\it irregularity}, denoted by $ir$. According to G\'erard-Levelt \cite{GL}, it can be computed as follows. Let us consider the following non-decreasing sequence of ${F}[[x]]$-lattices of $V$, starting from an arbitrary one $ V_0$: $\, V_{n+1}= V_n + \nabla(\partial)( V_n).$ 
  Then \begin{equation} ir\, M= \lim \frac{1}{n} \dim_{F} ( V_n/ V_0).\end{equation} 
  It is not clear from this formula that this is an integer. Actually, the integrality of the Turrittin-Levelt filtration follows from the expression of $ir$ in the cyclic case, \cf \eqref{fn}.
    
  \smallskip For $F=\C$, the irregularity is also degree function on the full subcategory $DMod_{{\C}(\{x\})}$ of analytic differential modules, but the associated slope filtration is not the restriction of the Turrittin-Levelt filtration, and does not seem to have any interest (in contrast to the $q$-analog, with the Adams-Sauloy filtration). 
 On the other hand, for  a cyclic analytic differential module $M=  \C(\{x\})\langle \partial \rangle/ \C(\{x\})\langle \partial \rangle\cdot P$, Malgrange has interpreted $ir \,M$ as the index of $P$ acting on $\C[[x]]/\C\{x\}$.

\medskip
  \subsection{Formal differential modules in several variables.}\label{TLf2} 
  Formal (integrable) differential modules in two or more variables are more mysterious, and decisive progress on unveiling their structure is very recent (\cf \cite{Sa}\cite{A5}\cite{Moc}\cite{Ke3}). 
  
  Let us just say a few words about the tannakian category $DMod_{R}$ for $R= \C[[x,y]][\frac{1}{x}]$, which is a non-full subcategory of $DMod_{F((x))}$ ($F=\cup\, \C((y^{1/n}))$ being the algebraic closure of $\C((y))$). 
 
 Let $\langle M\rangle$ (\resp $\langle M_{F((x))}\rangle$) be the tannakian subcategory of $DMod_{R}$ (\resp $DMod_{F((x))}$) generated by $M$ (\resp $M_{F((x))}$). The Turrittin-Levelt filtration of $\langle M_{F((x))}\rangle$ does not induce a filtration on $\langle M\rangle$ in general (a criterion is given in \cite[3.4.1]{A5}). In fact, the irregularity (in the sense of $\langle M_{F((x))}\rangle$) induces a degree function on $\langle M\rangle$, hence a slope filtration, which is not bounded in general however (some formal blow-ups are needed to fix this).

 An example is given by the differential module $M  $ with basis $m_1, m_2$ in which 
 $$\nabla (xd/dx)=  \begin{pmatrix}y/x & 0\\ - 1& 0 \end{pmatrix},\;\; \nabla (yd/dy)= \begin{pmatrix} -y/x & 0\\  1& 0 \end{pmatrix}.$$ The vector $m_1$ generates a differential submodule of slope $1$ (which is the universal destabilizing subobject), but the extension which gives $M$ does not split (whereas $M_{F((x))} = \gr M_{F((x))}$ splits).

  \medskip
    \subsection{Differential modules over the Robba ring.}\label{CMf} 

Let again $\sR$ be the Robba ring over the $p$-adic field $K$, endowed with a Frobenius  $\phi$ as in \ref{Kf}. 

The category $DMod_{\sR}$ of differential modules over $\sR$ is tannakian over $K$. Let $DMod^{(\phi)}_{\sR}$ be the tannakian (full) subcategory of differential modules admitting a Frobenius structure, \ie lying in the essential image of 
 $ F$-$Isoc_\sR  $.  
 
 \smallskip The highest break function associated to the {\it Christol-Mebkhout filtration} is given by the following recipe: $\rho(M)$ is the smallest $\la\in \R_{\geq 0}$ such that for any $r$ sufficiently close to $1$, $M$ admits a basis of  solutions in the open generic disk of radius $r^{1+\la}$. 
 
 \smallskip  The conditions of \ref{SSIM} are easily checked using this definition, which shows that associated (Christol-Mebkhout) filtration is $\otimes$-bounded, as far as one considers semisimple differential modules. 
 
  To check that it is $\otimes$-bounded on $DMod^{(\phi)}_{\sR}$, one would have to show that $M= \gr M$ for any $M$ (\cf \ref{extbou}), but this splitting property is non-trivial.

  \smallskip The degree function attached to the Christol-Mebkhout filtration is called the {\it $p$-adic irregularity}, denoted by $ir_p$. For a cyclic module $M= \sR\langle \partial \rangle/ \sR\langle \partial \rangle\cdot P$, Christol and Mebkhout have interpreted $ir_p\,M$ as a generalized index of $P$ acting on functions in the open unit disk. This interpretation shows that their filtration is integral \cite{CM}. 
  
   \smallskip Recent work by Baldassarri and by Kedlaya suggests that there should be a common framework for the Turrittin-Levelt and the Christol-Mebkhout filtrations, involving Berkovich geometry.

\medskip
\subsection{$q$-difference modules over the Robba ring.}\label{qCMf} 
This has a $q$-analog. Namely, let $q\in K^\phi$ be such that $ {\vert 1-q\vert < p^{-{1\over p-1}}}.$  
 
 The category $q$-$Mod_{\sR}$ of $q$-difference modules over $\sR$ is tannakian over $K$. Let $q$-$Mod^{(\phi)}_{\sR}$ be the tannakian (full) subcategory of differential modules admitting a Frobenius structure, \cf \cite[12.4]{AdV}. 
 
 There is a canonical ``functor of confluence"
 $\; q\text{\rm -}Mod^{(\phi)}_{\sR} \to DMod^{(\phi)}_{\sR} \,$
 which is an {equivalence of tannakian categories}.
 
 This functor is identity on the underlying $\sR$-modules. The differential structure arises as the limit of a canonical sequence $$M_m= (M, \nabla {(d_{q^{p^{m}}})})\in   {q^{p^{m}}}\text{-}Mod^{(\phi)}_\sR$$ related by isomorphisms
$  \Phi_\ast M_{m+1}\cong M_m $.

\smallskip
One can use this equivalence in order to transport the Christol-Mebkhout filtration to $q$-$Mod^{(\phi)}_{\sR}$.  
In order to show that it has the same description in terms of the radius of convergence of solutions in generic disks\footnote{this had been conjectured in \cite[4.3]{A3}; a sketch of the following construction was presented at the French-Nordic conference in Rejkyavik, january 2006.}, we may assume that $M_\infty$ (hence $M_0$) is purely of slope $\la$.
   Note that the convergence of $ \nabla{(d_{q^{p^{m}}})}$ to $\nabla(d/dx)$ implies that $M,  \nabla{(d_{q^{p^{m}}})}$ and $\nabla(d/dx)$ are all defined over some open annulus  $\sA(1-\epsilon , 1)$.

One proceeds in two steps:
 
1) Let $(M, \nabla{(d_q)})$ a $q$-difference module over  $\sA(1-\epsilon , 1)$. Then for $r$ close enough to $1$, the generic radius of convergence of $(M, \nabla{(d_q)})$ and $\Phi_\ast(M, \nabla{(d_q)})$ at $t_r$ and $t_r^{p}$ respectively coincide.
  
\smallskip  2) Let $q_i$ be a sequence converging to $1$, and $(M, \nabla{(d_{q_i})}) $ be a sequence of  $q_i$-difference modules over $\sA(1-\epsilon , 1)$. Then the generic radius of convergence of $(M, \nabla{(d_{q_i})}) $ at $t_r$ converges to the generic radius of convergence of $(M, \nabla{(d/dx)}) $ at $t_r$.
 
\smallskip  One concludes that
for any $i$, and for any $r$ close enough to $1$, the generic radius of convergence of $(M, \nabla{(d_{q_i})}) $ at $t_r$ is $r^{1+\la}$.

\medskip We will not go into further detail about this construction, since Pulita has recently given a more straightforward argument, in greater generality \cite[8.5.4]{Pu} (he relaxed the condition $ {\vert 1-q\vert < p^{-{1\over p-1}}} $, which allows to study other ``confluences " $q\to \zeta \in \mu_{p^\infty}$).

  \bigskip
  
  \section{Galois representations and filtrations of Hasse-Arf type.}  

\medskip
\subsection{Local Galois representations; case of perfect residue field.}

Let $(K, v)$ be a complete discretely valued field with perfect residue characteristic $k$. The {\it Hasse-Arf filtration} on representations of $G_K$ is constructed via a decreasing, left continuous, sequence of open normal subgroups 
 $(G_K^{(\la)})_{\la\in \Q_{\geq 0}}$ of $G_K$. 
The filtration of a representation $M$ is then defined by 
$$F^{\geq \la}M= \Ker (M\to M_{G^{(\la)} })$$
(where $M_{G^{(\la)}}$ stands for the coinvariants). 

\smallskip  Here, ``representation" means ``continuous representation with finite image over some field $F$ of characteristic $0$"; or, with appropriate interpretation, ``$\ell$-adic representation, with $\ell\neq $ car $ k$" \cite{Ka3} (more recently, the case of $p$-adic representations, with $p =$ car $ k$, has also been considered \cite{Co}\cite{Marm}).

\medskip The degree function attached to the Hasse-Arf filtration is called the {\it swan conductor}, denoted by $sw $. It takes values in $\Z$, by the Hasse-Arf theorem (there is also a cohomological interpretation of $sw$ as an index, due to Katz).

\medskip The image $G^{(\la)}$ of the group $G_K^{(\la)}$ in a given finite quotient $Gal(L/K)$ of $G_K$ is described as follows.  For any $i\in \N$, let $G_{(i)}$ be the subgroup of elements $g\in Gal(L/K)$ such that $v_L(g(a)-a)\geq i+1$ for any $a\in \sO_L$. Then the breaks of the locally constant non-increasing sequence $(G^{(\la)})_\la$ are given by
$\,  \la_i = \int_0^i [G_0:G_{-[-t]}]^{-1} dt , \;\;{\text{and}}\;\;
  G^{(\la_i)}= G_{(i)}.$ 
 
\medskip Let us consider the case $K= k((x))$, where $k$ is a perfect field of characteristic $p>0$. Let $F$ be the fraction field of the ring of Witt vectors of $k$, and let $\sR$ be the Robba ring over $F$. Then there are canonical $\otimes$-functors

\begin{equation}Rep^{fin}_F G_K  \to DMod^{(\phi)}_{\sR^{bd} }\to DMod^{(\phi)}_{\sR },  \end{equation}
   which, by the $p$-adic monodromy theorem, induce an equivalence of semisimple tannakian categories
\begin{equation}Rep^{fin}_F G_K \cong (DMod^{(\phi)}_{\sR })_{ssi}. \end{equation}

According to Tsuzuki \cite{Tsu}, this is compatible with the Hasse-Arf and Christol-Mebkhout filtrations respectively. In particular, $sw= ir_p$, and the integrality of $ir_p$ can be deduced from the Hasse-Arf theorem.

\medskip\subsection{Local Galois representations; case of imperfect residue field.} 
Complete discretely valued fields $K$ with imperfect residue field $k$ are more mysterious, and decisive progress on unveiling their higher ramification theory, as defined by Abbes and Saito \cite{AS1}\cite{AS2}, is very recent (\cf  \cite{CP}\cite{X1}\cite{X2}). The idea, launched by Matsuda \cite{Mat} and pursued by Kedlaya and Xiao, is to consider (integrable) differential modules over the Robba ring, with extra derivations acting non-trivially on $F$.

 \medskip \subsection{Local systems over a germ of punctured $p$-adic disk.}\label{Rf}
  Let $R$ be an ind-finite ring such that $R^\times$ contains both $p$ and a subgroup isomorphic to $\mu_{p^\infty}$ (for instance $R=\bar \F_\ell$). 
  
  Ramero \cite{R1}\cite{R2} introduced the $R$-linear abelian category $R$-$Loc_{\Delta^\times}^{br}$ of local system of $R$-modules with ``bounded ramification\footnote{this cohomological condition, which Ramero compares to an $L^1$-condition in harmonic analysis, restricts the wildness of the essential singularity at the puncture.}"  on the germ $\Delta^\times$ of punctured $p$-adic disk, and he endowed $R$-$Loc_{\Delta^\times}^{br}$ with a split slope filtration indexed by $\Q$, \cf  \cite[3.2.17]{R2}. He also gave a cohomological interpretation of the corresponding degree function.

  When $R$ is a field, $R$-$Loc_{\Delta^\times}^{br}$ is a tannakian category  over $R$, and the  {\it Ramero filtration} has all properties of a $\otimes$-bounded filtration (except that   ${\rm char}\, R \neq 0$).

 \newpage
 
 \addtocontents{toc}{{\bf V. Variation of Newton polygons in families.}\hfill\thepage}
\
\medskip
\begin{center}
\large\bf V. Variation of Newton polygons in families.
\end{center}

\bigskip In some situations, one has to consider not just one quasi-tannakian category with a slope filtration, but a whole family parametrized by a fixed space $S$.
Given a global object $M/S$, one can then ask about the variation of the Newton polygon of its fibers $M_s,\; s \in S$. 

We order the set of plane polygons by inclusion (since our Newton polygons are defined by concave functions, $N'\leq N$ if and only if $N'$ lies below $N$). 

\smallskip We review three such situations, which illustrate different behaviours.

\bigskip\section{Families of vector bundles.}

\subsection{} Let $M$ be a flat family of vector bundles over a smooth projective curve $X$ (over a field $F$), parametrized by a $F$-scheme $S$ of finite type. 
 
  For every point $ s$ of $S$, $M$ induces an object $M_s$ of $Vec_{X_{\kappa(s)}}$. Let $NP(M_s)$ denote its Newton polygon (with respect to the Harder-Narasimhan filtration).  
   
  \begin{thm}[Shatz \cite{Sh}] The function $\,s\in S \mapsto NP(M_s)\,$ is \emph{upper semicontinuous\footnote{with respect to the Zariski topology on $S$.}}.
 \end{thm} 
 
(Moreover, if $S$ is connected, the end-points of $NP(M_s)$ are constant).
 
 \smallskip
  This result allows to introduce a constructible stratification of $S$ by Newton polygons, which was studied by Shatz, Atiyah, Bott... 

   \medskip \subsection{} Shatz's proof relies on the possibility of specializing flags on the generic fibers of $M$, can be adapted to the case of families of filtered modules with respect to the Faltings-Rapoport filtration.
 
\medskip
 \section{$\rm F$-isocrystals.}
 
\subsection{}  The right notion of family of Frobenius-modules is Grothendieck's notion of $\rm F$-isocrystal over a noetherian scheme $S$ of characteristic $p$. 

For every geometric point $\bar s$ of $S$, such an $\rm F$-isocrystal $M$ induces an object $ M_{\bar s} $ of $\rm F$-$Isoc_K{(\kappa(\bar s))}$, where $\kappa(\bar s)$ denotes the fraction field of the ring of Witt vectors of $\kappa(\bar s)$. The Newton polygon of $ M_{\bar s}, $ (with respect to the Dieudonn\'e-Manin filtration) depends only on the point $s\in S$ under $\bar s$; we denote it by $NP(M_s)$.

 \begin{thm}[Grothendieck \cite{Gr}] The function $\,s\in S \mapsto NP(M_s)\,$ is \emph{lower semicontinuous}\footnote{Grothendieck uses the (usual) convention on Dieudonn\'e-Manin slopes (which leads in general to an ascending slope filtration and are the opposite of ours, \cf \ref{DMf}). With that convention, the polygons (bordered by convex functions) are upper semicontinuous. The sharp contrast with Shatz's theorem is not a matter of conventions.}.
 \end{thm} 
 
(Moreover, if $S$ is connected, the end-points of $NP(M_s)$ are constant).
 
 \smallskip
  This result  allows to introduce a constructible stratification of $S$ by Newton polygons, which was studied by De Jong and Oort \cite{dJO}. 
 
 \medskip \subsection{} A semicontinuity theorem similar to Grothendieck's is proven in \cite{Li} for Frobenius modules over the Robba ring with coefficients in a reduced $p$-adic affinoid algebra (instead of a $p$-adic field).
 
 On the other hand, Katz's proof \cite{Ka1} of Grothendieck's theorem can be adapted to the case of a family of $q$-difference modules (with no confluence of singularities) with respect to the Adams-Sauloy filtration. 
   
\medskip \section{Families of differential modules.}
  
  \subsection{} Let $f: X\to S$ be a smooth holomorphic family of connected curves parametrized by a complex manifold $S$, and let $Z\subset X$ be a hypersurface of $X$ which is finite etale over $S$. 
  
  Let $M$, be a vector bundle with meromorphic connection relative to $S$, and poles along $Z$ only. 
  
  Then for any $s\in S$, $M_s$ is a differential module on the curve $X_s$, with meromorphic singularities at the finite set of points $z$ such that $f(z)=s$. Let $NP_z(M_s)$ denote the Newton polygon (with respect to the Turrittin-Levelt filtration).

    \begin{thm}\cite[th. A.1]{A5} The function $\,z\in Z \mapsto NP_z(M_{f(z)})\,$ is \emph{lower semicontinuous}.  \end{thm} 
 
(Even if $S$ is connected, the right end-point of $NP(M_s)$ need not be constant: the irregularity may drop by specialization).
 
  \subsection{} If $Z$ is no longer assumed to be etale over $S$ (allowing the possibility of confluence of singularities), the result does not hold: the irregularity may jump by specialization. However, if the relative connection comes from an integrable connection on $X\setminus Z$, then the function 
  $$s\in S \mapsto \sum_{z, \, f(z)= s} \, ir_z(M_s)$$ is {lower semicontinuous} \cite{A5} (as was conjectured by Malgrange).

\newpage

\section*{Appendix: pseudo-$\otimes$-functors and rigidity.}\label{app}

\subsection{}  Let $\sC, \sC'$ be symmetric monoidal categories. 

\smallskip A \emph{pseudo-$\otimes$-functor} $(\varphi, \tilde \varphi, \varphi_\un)$ from $\sC$ to 
$\sC'$ consists of 

- a functor $\phi: \sC\to \sC'$,

- a morphism of functors $\tilde\phi: \, \otimes \circ (\phi, \phi)\Rightarrow \phi\circ \otimes\,$ from $\sC\times \sC$ to $\sC'$,

- an isomorphism $\phi_\un : \un \to \phi(\un)$ in $\sC'$,

\smallskip subject to the usual compatibilities with the monoidal structure. 

\smallskip Thus a pseudo-$\otimes$-functor\footnote{it is called a ``foncteur mono\"{\i}dal unitaire" in \cite{Ay}.} is a $\otimes$-functor (\ie a symmetric monoidal functor) if and only if $\tilde\phi$ is an isomorphism.

   \smallskip
    The aim of this appendix is to reconsider two diagrams of $\otimes$-functors which appear in \cite[I.4.3.3.3,  5.2.3.1]{Saa}, and whose commutativity is asserted there without proof. We shall give an argument which extends to the case of pseudo-$\otimes$-functors.  From there, we shall deduce a criterion for a pseudo-$\otimes$-functor to be a $\otimes$-functor.

    \begin{rem}\label{RA0} The compatibility with units, together with the condition that $\phi_\un$ is an isomorphism, imply that 
$\tilde\phi_{\un, Y}$ and $\tilde\phi_{X, \un}$ are always isomorphisms.

\smallskip The composition of two pseudo-$\otimes$-functors is a pseudo-$\otimes$-functor, with the rule
    $$ \widetilde{\psi\varphi}_{X,Y}= \psi(\tilde\varphi_{X,Y})\circ \tilde \phi_{\varphi(X), \varphi(Y)},\;\;  {\psi\varphi}_{\un}= \psi(\phi_\un)\circ \psi_\un.$$ 
    \end{rem}

     \begin{rem}\label{RA1} $\sC\times \sC$ has a natural structure of symmetric monoidal category. Moreover, $(\otimes, \tilde\otimes, \otimes_\un)$ is a $\otimes$-functor from $\sC\times \sC$ to $\sC$ if one sets
   
   - $\tilde\otimes_{(X_1, X_2),(Y_1, Y_2)} = 1_{X_1}\otimes c_{X_2, Y_1} \otimes 1_{Y_2}$, where $c$ denotes the commutativity constraint in $\sC$ (taking proper account of the associativity constraint). 
     
   - $\otimes_\un =$ the canonical isomorphism $\un \overset{\cong}{\to} \un\otimes \un$ in $\sC$.  
   \end{rem}
   
   \subsection{} A \emph{morphism}  $$u: (\varphi, \tilde \varphi, \varphi_\un) \Rightarrow (\psi, \tilde\psi, \psi_\un)$$ between pseudo-$\otimes$-functors from $\sC$ to 
   $\sC'$ is a natural transformation  
 $ u :  \varphi  \Rightarrow  \psi $ which is compatible with $\tilde\varphi$ and $\tilde\psi$, and with the  constraints. 
 
 \begin{rem}\label{RA2}  Actually, the compatibility with the unit constraints is automatic, and one has 
 $u_\un = \psi_\un
\circ \varphi_1^{-1}$.  
\end{rem} 
  
  \begin{rem}\label{RA3} Given a pseudo-$\otimes$-functor $\phi: \sC\to \sC'$, the functors $\varphi  = \otimes \circ (\phi, \phi)$ and $\psi = \phi\circ \otimes$ have a natural structure of pseudo-$\otimes$-functors from $\sC\times \sC$ to $\sC'$, with 
  $\tilde \phi_{(X_1, X_2), (Y_1,Y_2)} $ (\resp $\tilde \psi_{(X_1, X_2), (Y_1,Y_2)} $) given by 
  $(\tilde\phi_{X_1,Y_1}\otimes \tilde\phi_{X_2,Y_2}) \circ (1_{X_1}\otimes c_{X_2, Y_1}\otimes 1_{Y_2})$
 (\resp $\phi(1_{X_1}\otimes c_{X_2, Y_1}\otimes 1_{Y_2})\circ \tilde\phi_{X_1\otimes Y_2, X_2\otimes Y_2}$) (taking proper account of the associativity constraint). 
 
 Moreover, $\tilde\phi: \varphi  \Rightarrow \psi$ is then a morphism of pseudo-$\otimes$-functors. 
 The compatibility with $\tilde \varphi$ and $\tilde\psi$ and the constraints amounts to some identities which are consequences of the fact that $\phi$ itself is a pseudo-$\otimes$-functor. For instance, the compatibility with $\tilde \varphi$ and $\tilde\psi$ amounts to  
 $$\tilde \phi_{X_1\otimes Y_1, X_2\otimes Y_2}\circ (\tilde\phi_{X_1,Y_1}\otimes \tilde\phi_{X_2,Y_2}) \circ (1_{X_1}\otimes c_{X_2, Y_1}\otimes 1_{Y_2})$$
 $$= \phi(1_{X_1}\otimes c_{X_2, Y_1}\otimes 1_{Y_2})\circ \tilde\phi_{X_1\otimes Y_1, X_2\otimes Y_2}\circ (\tilde\phi_{X_1,X_2}\otimes \tilde\phi_{Y_1,Y_2}).$$

\end{rem}

\subsection{}  Let us now assume that $\sC$ and $\sC'$ are rigid, and let us denote by 
$$D: \sC^{op}\to \sC,\;D': \sC'^{op}\to \sC',$$ or sometimes simply $()^\vee$, the duality equivalences.  

If one identifies a morphism $f: X\to Y$ in $\sC$ with a morphism $Y\to X$ in $\sC^{op}$, $Df$ is the transpose ${}^tf$ of $f$, and is characterized by the commutative square
 \[\begin{CD}     Y^\vee  \otimes X @>{}^tf\otimes 1>>       X ^\vee  \otimes X\\   @V{1\otimes f}VV @VV{ev_X}V     \\      Y^\vee\otimes Y@>>{ ev_Y }>   \un  . \\    \end{CD}\]  
 
The canonical isomorphisms 
$$\tilde D_{X, Y} : X^\vee \otimes Y^\vee \overset{\cong}{\to} (X\otimes Y)^\vee, \;\; D_\un: \un \overset{\cong}{\to} \un^\vee $$ make $D$ into a $\otimes$-equivalence (\cf \cite[I, 5.1.3]{Saa}). Similarly for $D'$.  
 
\subsection{}  Let $(\phi, \tilde \phi, \phi_\un)$ be a pseudo-$\otimes$-functor from $\sC$ to 
$\sC'$. It induces a functor 
$$\phi': \sC^{op} \to \sC'^{op},$$ and there is a canonical natural transformation 
$$\hat\phi: \phi\circ D  \Rightarrow D'\circ \phi',\;\; \hat\phi_X: \phi(X^\vee)\to (\phi(X))^\vee $$
 which is characterized by the commutativity of the square
 \[\begin{CD}    \phi(X^\vee) \otimes \phi(X) @>\hat\phi_X\otimes 1>>      \phi(X)^\vee  \otimes \phi(X)\\   @V{\tilde\phi_{X^\vee, X}}VV @VV{ev_{\phi(X)}}V     \\      \phi(X^\vee\otimes X)@>>{\phi(ev_X)}>   \phi(\un) . \\    \end{CD}\]  
In particular $\hat\phi_\un \bullet 1 = \phi(ev_\un)\circ \tilde\phi_{\un^\vee, \un},$ hence $\hat\phi_\un$ is an isomorphism.

  \begin{rem}\label{RA4} When $\phi$ is a $\otimes$-functor, \ie when $\tilde \phi$ is an isomorphism, then there is a natural $\otimes$-structure on $\phi'$ (given by the transpose of $\tilde\phi^{-1}$) for which $\hat\phi$ becomes an isomorphism of $\otimes$-functors  \cite[I,5.2]{Saa}.   
   \end{rem}
 
 \begin{rem}\label{RA5} For the composition of pseudo-$\otimes$-functors, one has the formula
 $$\widehat{\psi\varphi}_X =  \hat\psi_{\varphi(X)} \circ  \hat\varphi_X $$ which follows from the formula in \ref{RA0} applied to the pair $(X^\vee, X)$.
 
 On the other hand, for $\otimes$ considered as a $\otimes$-functor $\sC\times \sC\to \sC$, one has $\hat(\otimes)_{X_1, X_2}=  \tilde D_{X_1, X_2}$.
 \end{rem}
 
\subsection{} Let $u: (\varphi, \tilde \varphi, \varphi_\un) \Rightarrow (\psi, \tilde\psi, \psi_\un)$ be a morphism of pseudo-$\otimes$-functors from $\sC$ to $\sC'$. 

 \begin{lemma}   One has $\, \hat\varphi = {}^tu \circ \hat\psi \circ (u\ast D)$, \ie for any 
$X$, the composed morphism
 $$ \varphi(X^\vee)\overset{u_{X^\vee}}{\to}\psi(X^\vee) \overset{\hat\psi_X}{\to} (\psi(X))^\vee \overset{{}^tu_{X}}{\to}(\phi(X))^\vee $$
  is $\hat\varphi_X$.
    \end{lemma}

\begin{rem}\label{RA6} In \cite[I.5.2.2.1]{Saa}, it is proven that if $\varphi$ and $\psi$ are $\otimes$-functors, then $\hat\varphi $ and $\hat\psi$ are isomorphisms. In \cite[I.5.2.3.1]{Saa}, it is asserted (without proof) that  $\, \hat\varphi_{X^\vee}^{-1}\circ {}^tu_{X^\vee} \circ \hat\psi_{X^\vee}$ is inverse to $u_X$. The above formula shows that it is left-inverse. To show that it is also right-inverse, one can check that $$  (\hat\varphi^{-1} \circ {}^tu\circ \hat\psi)\ast D: \,(\psi, \tilde\psi, \psi_\un)\Rightarrow(\varphi, \tilde \varphi, \varphi_\un) $$   is a morphism of $\otimes$-functors and apply the lemma to it.  \end{rem}

\proof   We have to show that  
$$ev_{\varphi(X)}\circ ({}^tu_X\otimes 1_{\varphi(X)})\circ (\hat\psi_X\otimes 1_{\varphi(X)})\circ (u_{X^\vee}\otimes 1_{\varphi(X)})= ev_{\varphi(X)}\circ (\hat\varphi_X\otimes 1_{\varphi(X)}).$$
Since $ev_{\varphi(X)}\circ ({}^tu_X\otimes 1_{\varphi(X)})= ev_{\psi(X)}\circ (1_{\psi(X)^\vee}\otimes u_X)$, this amounts to 
$$ev_{\psi(X)}\circ (\hat\psi_X\otimes 1_{\varphi(X)})\circ (u_{X^\vee}\otimes u_X)= ev_{\varphi(X)}\circ (\hat\varphi_X\otimes 1_{\varphi(X)}).$$
 Since $u$ is a morphism of pseudo-$\otimes$-functors, one has a commutative diagram 
 \[\begin{CD}    \varphi(X^\vee) \otimes \varphi(X) @>\tilde\varphi_{X^\vee, X}>>      \varphi(X^\vee  \otimes  X) @>\varphi(ev_X)>> \varphi(\un)\\ 
   @V{u_{X^\vee}\otimes u_X} VV @VV{u_{X^\vee\otimes X}}V  @VV{u_\un}V     \\     
 \psi(X^\vee) \otimes \psi(X) @>\tilde\psi_{X^\vee, X}>>      \psi(X^\vee  \otimes  X) @>\psi(ev_X)>> \psi(\un) . \\    \end{CD}\]  
   On the other hand, the composed morphism in the top row is 
  $\,\varphi_\un^{-1}\circ  ev_{\varphi(X)}\circ (\hat\varphi_X\otimes 1)$, while the composed morphism in the bottom row is 
  $\,\psi_\un^{-1}\circ  ev_{\psi(X)}\circ (\hat\psi_X\otimes 1)$. This establishes the required formula (taking remark \ref{RA3} into account).\qed

Let $(\phi, \tilde \phi, \phi_\un)$ be a pseudo-$\otimes$-functor from $\sC$ to 
$\sC'$.

\begin{lemma} For any $X,Y$, one has a commutative diagram 
   \[\begin{CD}    \phi(X^\vee) \otimes \phi(Y^\vee) @>\tilde\phi_{X^\vee, Y^\vee}>>      \phi(X^\vee  \otimes  Y^\vee)\\   @V{ \hat\phi_X\otimes \hat\phi_Y}VV @VV{\phi(\tilde D_{X,Y})}V     \\      \phi(X)^\vee \otimes \phi(Y)^\vee   &&    \phi((X\otimes Y)^\vee) \\  @V{\tilde D_{\phi(X),\phi(Y)}}VV @VV{\hat\phi_{X\otimes Y}}V     \\     ( \phi(X)\otimes \phi(Y))^\vee@<<{{}^t\tilde\phi_{X,Y}}<   \phi(X\otimes Y)^\vee . \\    \end{CD}\]  
 \end{lemma} 

\begin{rem}\label{RA7} The commutativity of this diagram is asserted without proof in \cite[I.4.3.3.3]{Saa}  (in the case of a $\otimes$-functor).
\end{rem}

\proof  Taking into account remark \ref{RA3}, and with the same notation, we can apply the previous lemma to $u= \tilde \phi: \varphi\Rightarrow \psi$ (with $\sC$ replaced by $\sC\times \sC$). This gives 
$$\hat \varphi_{X,Y} = {}^t\tilde\phi_{X,Y} \circ \hat\psi_{X,Y} \circ  \tilde\phi_{X^\vee, Y^\vee}.$$ 
It remains to identify $\hat \varphi_{X,Y} $ with $\tilde D_{\phi(X),\phi(Y)} \circ ( \hat\phi_X\otimes \hat\phi_Y)$  and  $\hat \psi_{X,Y} $ with $\hat\phi_{X, Y}\circ \phi(\tilde D_{X,Y})$. This follows from remark \ref{RA5}.
  \qed  

\begin{cor}\label{CA} Assume that $\sC'$ is additive and pseudo-abelian (\ie idempotent morphisms have a kernel and a cokernel), and that there is no non-zero object of rank $0$.

Then  $\tilde\phi$ is an isomorphism (\ie $\phi$ is a $\otimes$-functor) if and only if  $\hat\phi$ is an isomorphism.
\end{cor}  

\proof We already know that $\hat\phi$ is an isomorphism if $\phi$ is a $\otimes$-functor (\cite[I.5.2.2.1]{Saa}). Let us prove the converse. We first remark that objects occurring in the above diagram have the same rank. If $\hat\phi$ is an isomorphism, the vertical morphisms are isomorphisms, so that $f =\tilde\phi_{X^\vee, Y^\vee}$ has a left inverse $g$. The kernel and the cokernel of the idempotent $fg$ have rank $0$, hence are $0$.  Therefore $f$ is an isomorphism. 
\qed

  \begin{ex} Let $\sC $ be either the category of polarized Hodge structures, or the category of numerical motives. In both cases, this is a semisimple tannakian category with rational coefficients\footnote{in the latter case, it is assumed that the commutativity constraint is twisted by a sign according to the Koszul rule; in particular, it is assumed that the K\"unneth projectors on the even part of the cohomology are induced by algebraic correspondences.}. Let us consider the Grothendieck {\it coniveau filtration}: for $\la\in \Z$, $F^{\geq \la} M $ is the greatest subobject of $M$ such that the twist $F^{\geq \la} M (\la)$ is effective. This gives a split slope filtration on $\sC$. 
  
  One has 
$F^{\la_1}M_1\otimes  F^{\la_2}M_2\subset F^{\la_1+\la_2}M_2$, so that $\gr$ has a natural structure of pseudo-$\otimes$-endofunctor $(\gr = id , \tilde \gr, \gr_\un =1)$. It is not a $\otimes$-functor.

This can be seen on the example $M_1=M_2=H^1$ of an elliptic curve without complex multiplication: $\tilde\gr_{M,M} $ (\resp $\hat \gr_{M\otimes M}$) is an isomorphism on the $S^2$ component, and zero on the $\displaystyle\wedge^2$ component. 
 In fact, the coniveau filtration is not determinantal. \end{ex}

\medskip
  \end{sloppypar}

\bigskip

\section*{Index}   
{\Small

$\otimes$-bounded (slope filtration)  \dotfill \ref{bsf}, 

break \dotfill \ref{break}, 
 
\medskip category with kernels and cokernels \dotfill \ref{ckc},

conservative (functor) \dotfill \ref{pi}, 

\medskip degree (function) \dotfill \ref{sf}, 

determinantal (slope filtration)  \dotfill \ref{dsf}, 

 \medskip epi, monic   \dotfill \ref{kcok}, 

exact (filtration) \dotfill \ref{exa}, 

exact (functor) \dotfill \ref{kcok} 

extension \dotfill \ref{kcok}, 

\medskip  flag \dotfill \ref{FLAG}, \ref{P3},

\medskip highest break (function) \dotfill \ref{hbf},

\medskip integral (slope filtration)  \dotfill \ref{Int}, 

\medskip kernel, cokernel \dotfill \ref{kcok},

\medskip  left abelian envelop  \dotfill \ref{qac}, 

\medskip  $\otimes$-multiplicative (slope filtration)  \dotfill \ref{msf}, 

multiplicity (of a break)  \dotfill \ref{mul}, 

\medskip  Newton polygon  \dotfill \ref{newp}, 

\medskip  Picard group  \dotfill \ref{deter}, 

proto-abelian (category)  \dotfill \ref{pac}, 

pseudo $\otimes$-functor  \dotfill Appendix, 

pull-back, push-out  \dotfill \ref{pbpo}, 

\medskip  quasi-abelian (category)  \dotfill \ref{qac}, 

   quasi-tannakian (category)  \dotfill \ref{qtrf}, 

\medskip  rank (function) \dotfill \ref{RANK}, \ref{qtr},

Rees deformation \dotfill \ref{Rees}, \ref{remqt},

\medskip  semistable (object)  \dotfill \ref{sst}, 

short exact sequence \dotfill \ref{kcok}, 

slope (filtration) \dotfill \ref{slf},

slope (function) \dotfill \ref{sf}
 
 split (slope filtration) \dotfill \ref{ssf}, 

strict (quotient, subobject, epi, monic)  \dotfill \ref{kcok},  

stability structure (on a triangulated category) \dotfill \ref{sstc}, 

stable (object)  \dotfill \ref{sst},  

 strongly exact (filtration)  \dotfill \ref{exa}, 

\medskip  universal destabilizing subobject \dotfill \ref{uds}. 

\bigskip
\medskip

{\it Special slope filtrations:}  

\medskip
Adams-Sauloy filtration \dotfill \ref{ASVf}, 

Bruasse filtration  \dotfill \ref{bru}, 

Christol-Mebkhout filtration  \dotfill \ref{CMf}, 

coniveau filtration  \dotfill Appendix, 

Dieudonn\'e-Manin filtration \dotfill \ref{dmf},  \ref{exsp}, \ref{EXFF}, \ref{DMf},

Di Vizio filtration  \dotfill \ref{ASVf},

Faltings-Rapoport filtration \dotfill \ref{FRf}, 

Fargues filtration \dotfill \ref{Ff},

Fontaine filtration \dotfill \ref{FFf}, 

 Grayson-Stuhler filtration \dotfill \ref{gsf}, \ref{exs1}, \ref{GSf},

Harder-Narasimhan filtration   \dotfill \ref{hnf}, \ref{exs1}, \ref{exex}, \ref{extri}, \ref{EXF}, \ref{extr}, \ref{EXFF}, \ref{vbc}, \ref{vbhd},

Hasse-Arf filtration \dotfill \ref{haf},  \ref{exsp},  \ref{Int},  \ref{bsf},  \ref{haf}, 

Hoffmann-Jahnel-Stuhler filtration  \dotfill \ref{HJSf},

Kedlaya filtration \dotfill \ref{Kf},

Moriwaki filtration \dotfill \ref{Mf},

Ramero filtration \dotfill \ref{Rf}, 

Turrittin-Levelt filtration \dotfill \ref{tlf},  \ref{exsp}, \ref{bsf}, \ref{TLf}, \ref{TLf2},  

weight filtration \dotfill \ref{weight}.

\bigskip


\begin{thebibliography}{I}  
 
 
\medskip   

     \smallskip \bibitem{AS1} A. Abbes, T. Saito, Ramification of local fields with imperfect residue fields. 
Amer. J. Math., {\bf 124} (5):879Ð920, 2002. 

     \smallskip \bibitem{AS2} A. Abbes, T. Saito, Ramification of local fields with imperfect residue fields, I I. 
Doc. Math., (Extra Vol.):5Ð72 (electronic), 2003. Kazuya KatoÕs fiftieth birthday. 



     \smallskip \bibitem{Ad} C. Adams,  On the linear ordinary $q$-difference equation,  Ann. of Math. (2) {\bf 30} (1928/29), no. 1-4, 195--205.
   
     \smallskip \bibitem{A2}  Y. Andr\'e,  Filtrations de type Hasse-Arf et monodromie $p$-adique,
  Invent. Math.  {\bf 148} (2002),  285-317.
 
  \smallskip \bibitem{A3} Y. Andr\'e, Galois representations, differential equations and 
$q$-difference equations: sketch of a $p$-adic unification, in {\it Complex analysis, dynamical systems, summability of divergent series and Galois theories. I. Volume in honor of Jean-Pierre Ramis},   M. Loday-Richaud ed., S.M.F. Ast\'erisque {\bf 296} (2004), 43-53. 
  
          \smallskip 
 \bibitem{A4}  Y. Andr\'e, Motifs de dimension finie, Exp. 929, S\'em. Bourbaki, in: Ast\'erisque   {\bf 299}  (2005), 115--145. 
 
   \smallskip 
 \bibitem{A5}  Y. Andr\'e, {Structure des connexions m\'eromorphes for\-melles de plusieurs variables et semi-continuit\'e de l'irr\'egularit\'e}, Invent. Math {\bf 170} (2007), 147-198.
 
          \smallskip 
 \bibitem{A6}  Y. Andr\'e, Sur le produit tensoriel de r\'eseaux hermitiens semistables, in preparation.
 
  
      \smallskip \bibitem{AdV} Y. Andr\'e, L. Di Vizio, 
$q$-difference equations and $p$-adic local monodromy, in {\it Complex analysis, dynamical systems, summability of divergent series and Galois theories. I. Volume in honor of Jean-Pierre Ramis},   M. Loday-Richaud ed., S.M.F. Ast\'erisque {\bf 296} (2004), 55-111. 
  
         \smallskip \bibitem{AK} Y. Andr\'e, B. Kahn, Nilpotence, radicaux et structures mono\"{\i}dales, Rendiconti. Sem. Padova {\bf 108} (2002), 107-191.
       
             \smallskip \bibitem{Arf}  C. Arf, Untersuchungen \"uber reinverzweigte Erweiterungen diskret bewerteter perfekter K\"orper, J. reine angew. Math. {\bf 181} (1940), 1-44.
       
         \smallskip \bibitem{ABo} M. Atiyah, R. Bott, The Yang-Mills equations over Riemann surfaces, Phil. Trans. Roy. Soc. London, {\bf 308} (1983), 523-615.
       
    \smallskip \bibitem{Ay} J. Ayoub, Les six op\'erations et le formalisme des cycles \'evanescents dans le monde motivique, Ast\'erisque {\bf 314}, 315 (2007).
   
       \smallskip\bibitem {BRS} L. Barbieri-Viale, A. Rosenschon, M. Saito, Deligne's conjecture on $1$-motives,
 Ann. of Math. (2)  {\bf 158}  (2003),  no. 2, 593--633.
 
  \smallskip \bibitem{Be} L. Berger,  \'Equations diff\'erentielles $p$-adiques et $(\phi, N)$-modules filtr\'es,  Ast\'erisque {\bf 319} (2008), 13-38.

    \smallskip \bibitem{BBD}   A. Be\u\i linson,  J. Bernstein, P. Deligne, Faisceaux pervers, in {\it Analysis and topology on singular spaces}, I (Luminy, 1981),  5--171, Ast\'erisque {\bf 100}, S. M. F. (1982). 
    
      \smallskip 
 \bibitem{BvB} A. Bondal, M. van der Bergh, Generators and representability of functors in commutative and noncommutative algebraic geometry, Mosc. Math. J. {\bf 3} 1, 2003,  1Ð36.
 
 
  \smallskip 
 \bibitem{Bo}  T. Borek, Successive minima and slopes of hermitian vector bundles over number fields,
 J. Number Theory {\bf 113} 2 (2005), 380-388.
 
   \smallskip 
 \bibitem{Bos}  J.-B. Bost, Hermitian vector bundles and stability, talk at Oberwolfach, 
``Algebraische Zahlentheorie, july 1997.
 
    \smallskip 
 \bibitem{BK}  J.-B. Bost, K. K\"unnemann, Hermitian vector bundles and extension groups on arithmetic schemes. I. geometry of numbers, preprint (2007), ArXiv0701343. 

  \smallskip  \bibitem{B} N. Bourbaki, {\it Commutative algebra}, Chapters 1--7. Translated from the French. Reprint of the 1989 English translation. Elements of Mathematics, Springer-Verlag, 1998.
 

   \smallskip 
 \bibitem{Br} T. Bridgeland, Stability conditions on triangulated categories, Annals of Mathematics {\bf 166} (2007), 317Ð345
 
 
 
   \smallskip  \bibitem{Bra} S. Bradlow. Hermitian-Einstein inequalities and Harder-Narasimhan filtrations, Int. Journal of Maths. {\bf 6} (5) (1995), 645Ð656.
 
 
  \smallskip   \bibitem{Bru} L. Bruasse, Filtration de Harder-Narasimhan pour des fibr\'es complexes ou des faisceaux
sans torsion, Annales de lÕInstitut Fourier {\bf 53} 2 (2003), 541-564.

    \smallskip 
 \bibitem{Ch} H. Chen, Harder-Narasimhan categories, preprint (2007), preprint  ArXiv 0706.2648v2.
 
    \smallskip 
 \bibitem{Che} H. Chen, Maximal slope of tensor product of Hermitian vector bundles, preprint (2008), preprint  hal00151961v2.
 
    \smallskip \bibitem{CC}   F. Cherbonnier, P. Colmez, Repr\'esentations p-adiques surconvergentes.
Invent. Math. {\bf 133} no. 3,  (1998) 581Ð611.
     
    
\bibitem{CP} B. Chiarellotto, A. Pulita. Arithmetic and Differential Swan Conductors of rank 
one representations with finite local monodromy, preprint arXiv: math.NT/0711.0701v2. 
 
  \smallskip 
 \bibitem{CM}  G. Christol, Z. Mebkhout, {\it Sur le th\'eor\`eme de l'indice des \'equations diff\'erentielles $p$-adiques III}, Ann. of Maths  {\bf 151} (2000), 385-457.  

 \smallskip   \bibitem{Co} P. Colmez, Conducteur d'Artin d'une repr\'esentaiton de De Rham, preprint IMJ, 2003.
 
 
   \smallskip   \bibitem{CF} P. Colmez, J.M Fontaine, Construction des repr\'esentations $p$-adiques
semi-stables. Inv. Math. {\bf 140} (2000), 1Ð43.
 
  \smallskip   \bibitem{dJO} A. De Jong, F. Oort, Purity of the stratification by Newton polygons, J. A. M. S. {\bf 13} n. 1 (1999), 209-241.
   
   \smallskip    \bibitem {De1} P. Deligne, Th\'eorie de Hodge II, Publ. Math. I.H.E.S. {\bf 40} (1971), 5-57. 

   \smallskip  \bibitem {De2} P. Deligne, Th\'eorie de Hodge III, Publ. Math. I.H.E.S. {\bf 44} (1974), 5-78. 

 
  
  \smallskip   \bibitem{DSP} E. De Shalit, O. Parzanchevski, On tensor products of semistable lattices, preprint Jerusalem, 2006.  

     \smallskip\bibitem{Di}  J. Dieudonn\'e, Lie groups and Lie hyperalgebras over a field of characteristic $p > 0$. IV, 
Amer. J. Math. {\bf 77} (1955), 429Ð452. 
   
  \smallskip\bibitem{dV} L. Di Vizio,  Local analytic classification of q-difference equations with $\vert q\vert =1$,
Journal of Noncommutative Geometry {\bf 3} (1) (2009), 125--149. ArXiv:math/0802.4223.
 
   \smallskip\bibitem{dVRSZ}  L. Di Vizio,  J.-P. Ramis,  J. Sauloy, C. Zhang, \'Equations aux $q$-diff\'erences, Gazette des Math\'ematiciens, Soci\'et\'e Math\'ematique de France  (2003), 20-49.
  
    \smallskip 
 \bibitem{F1} G. Faltings, Stable $G$-bundles and projective connections, J. of Alg. Geom. {\bf 2} (1993), 549-576.
 
   \smallskip 
 \bibitem{F2} G. Faltings, Mumford-Stabilit\"at in der algebraiscen Geometrie, in {\it Proc. Intern. Congress Math. Z\"urich 1994}, 648-655, Birkh\"auser (1995).
 
  \smallskip 
 \bibitem{FW}  G. Faltings, G. W\"ustholz,  Diophantine approximations on projective spaces.  Invent. Math.  {\bf 116}  (1994),  no. 1-3, 109--138.
 
   \smallskip 
 \bibitem{Fa}  L. Fargues, La filtration de Harder-Narasimhan des sch\'emas en groupes finis et plats, preprint Orsay 2008.
 
  \smallskip 
 \bibitem{Fo} J.-M. Fontaine, Repr\'esentations p-adiques des corps locaux I. The Grothendieck
Festschrift, Vol. II, 249Ð309, Progr. Math. 87, Birkh\"auser Boston,
Boston, MA 1990.
 
   \smallskip  \bibitem{Gau} E. Gaudron, Pentes des fibrŽs vectoriels adŽliques sur un corps global,   Rendiconti di Padova (2008).
     
   \smallskip  \bibitem{GL}    R. G\'erard, A. Levelt, Inveraints mesurant l'irr\'egularit\'e en un poinjt sinsgulier d'\'equations diff\'erentielles lin\'eaires, Ann. Inst. Fourier (1973), 157-195.
 
   \smallskip  \bibitem{Gi2} D. Gieseker, Global moduli for surfaces of general type, Invent. Math. {\bf 43} (1977) 233-282.
 


  \smallskip   \bibitem{GKR} A. Gorodentsev, S. Kuleshov, A. Rudakov, Stability data and $t$-structures on a triangulated category. preprint ArXiv:math/0312442. 


    \smallskip 
 \bibitem{Gra} D. Grayson, Higher algebraic $K$-theory II (after Daniel Quillen), in {\it Algebraic $K$-theory}, Proc. Conf. Nothwestern Univ. 1976, 217-240, Springer Lect. Notes {\bf 551} (1976).
 
    \smallskip 
 \bibitem{Gr} A. Grothendieck, {\it Groupes de Barsotti-Tate et cristaux de Dieudonn\'e}, Univ. de Montr\'eal 1974.
 
   \smallskip 
 \bibitem{HRS} D. Happel, I. Reiten, S. Smal\Small{$\phi$}, {\it Tilting in abelian categories and quasitilted algebras}, Memoir A. M. S. {\bf 575} (1996).
 
  \smallskip 
 \bibitem{HN} G. Harder, M. Narasimhan, On the cohomology groups of moduli spaces of vector bundles on curves, Math. Annalen (1974), 215-248.
 
   \smallskip 
 \bibitem{Ha} H. Hasse, F\"uhrer, Diskriminante une Verzweigungsk\"orper relativ Abelscher Zahlk\"orper, J. reine angew. math. {\bf 162} (1930), 169-184.
 
  \smallskip 
 \bibitem{Hi} N. Hitchin, The self-duality equations on a Riemann surface, Proc. London Math. Soc. (3), {\bf 55} (1987), 59-126.
 
  \smallskip 
 \bibitem{HJS}  N. Hoffmann, J. Jahnel and U. Stuhler, Generalized vector bundles on curves,  J. Reine Angew. Math.  {\bf 495}  (1998), 35--60. 
  
     \smallskip   \bibitem{Ju} M. Jurchescu, Theory of categories. (Romanian) 1966, 73--240, Editura Acad. Republicii Socialiste Rom\^ania, Bucharest.
 
   \smallskip   \bibitem{Ka1} N. Katz, {Slope filtrations of $F$-crystals }, in {\it Journ\'ees de G\'eom\'etrie alg\'ebrique de Rennes}, July 1978  Ast\'erisque {\bf 64} (1979).

 \smallskip   \bibitem{Ka2} N. Katz, On the calculation of some differential Galois groups, Invent. Math. {\bf
87} (1987), 13-61.   

 \smallskip   \bibitem{Ka3} N. Katz,  {\it Gauss sums, Kloosterman sums, and monodromy groups}, Annals of Mathematics Studies, {\bf 116} Princeton University Press, Princeton, NJ, 1988.
 
  \smallskip   \bibitem{Ke1} K. Kedlaya, Slope filtrations revisited, Documenta Mathematica {\bf 10} (2005), 447Ð525.

 \smallskip   \bibitem{Ke2} K. Kedlaya, Slope filtrations for relative Frobenius, Ast\'erisque {\bf 319} {\bf 308} (2008), 259Ð301.

 \smallskip   \bibitem{Ke3} K. Kedlaya, Good formal structures on flat meromorphic connections, I: Surfaces, preprint arXiv:0811.0190.
   
 \smallskip   \bibitem{Ke4} K. Kedlaya, {\it $p$-adic differential equations}, book draft 2008..
   
  \smallskip   \bibitem{Kob}  S. Kobayashi, {\it Differential geometry of complex vector bundles}, Princeton University Press (1987).   
  
 \smallskip   \bibitem{Laz} R. Lazarsfeld, {\it Positivity in algebraic geometry.} II. Positivity for vector bundles, and multiplier ideals. Ergeb. der Math. und ihrer Grenzgebiete. 3. Folge.  {\bf 49}. Springer-Verlag, Berlin, (2004).
 
  \smallskip   \bibitem{L}  A. Levelt, Jordan decom\-po\-si\-tion for a class of singular
differential operators,
Ark. Math. {\bf 13} (1975), 1-27.
 
 \smallskip   \bibitem{Li}  R. Liu, Slope filtrations in families, preprint ArXiv 08090331.
 
  \smallskip   \bibitem{ML} S. Mac Lane, {\it Categories for the working mathematician}, 2nd ed., Springer(1988).
 
 \smallskip   \bibitem{Mal} B. Malgrange, {\it \'Equations diff\'erentielles \`a coefficients polynomiaux}, Progress in Mathematics, Birkh\"auser  (1991).
 
  \smallskip   \bibitem{Man}  Yu. Manin, Theory of commutative formal groups over fields of finite characteristic 
(Russian), Uspehi Mat. Nauk {\bf 18} (1963), 3Ð90. 


 \smallskip   \bibitem{Marm} A. Marmora, Irr\'egularit\'e et conducteur de Swan p-adiques, Doc. Math. {\bf 9} (2004), 413-433.
  
    \smallskip   \bibitem{Maru} M. Maruyama, The theorem of Grauert-M\"uhlich-Spindler, Math. Ann. {\bf 255} (1981), 317-333.
    
     \smallskip   \bibitem{Mat} S. Matsuda, Conjecture on Abbes-Saito filtration and Christol-Mebkhout filtration, in {\it Geometric aspects of Dwork theory}. Vol. I, II, pages 845Ð856. Walter de Gruyter, Berlin, 2004. 

  \smallskip   \bibitem{MR1} V. Mehta, A. Ramanathan,    Semistable sheaves on projective varieties and their restriction to curves, Math. Ann., 258 (1982), 213-224. 
   
  \smallskip   \bibitem{Mi}  Y. Miyaoka,   The Chern classes and Kodaira dimension of a minimal variety.  {\it Algebraic geometry, Sendai, 1985},  449--476, Adv. Stud. Pure Math., 10, North-Holland, Amsterdam, (1987).
   
   
   \smallskip   \bibitem{Moc}  T. Mochizuki, Good formal structure for meromorphic flat connections on smooth projective surfaces, preprint arXiv:0803.1346.


 
  \smallskip   \bibitem{Mo}  A. Moriwaki, Subsheaves of a Hermitian torsion-free coherent sheaf on an arithmetic variety, preprint ArXiv:math/0612268.

  \smallskip   \bibitem{Mu}  D. Mumford, Projective invariants of projective structures and applications.   {\it Proc. Internat. Congr. Mathematicians (Stockholm, 1962)}  pp. 526--530 Inst. Mittag-Leffler, Djursholm (1963).

       \smallskip   \bibitem{NS} M. Narasimhan, C. Seshadri, Stable and unitary vector bundles on a compact Riemann surface, Ann. Math. {\bf 82} (1965), 540-567.
     
         \smallskip   \bibitem{OS1}   P. O'Sullivan, Finite dimensionality of reductive envelopes, appendix to \cite{AK}.

  \smallskip   \bibitem{OS2}   P. O'Sullivan, The structure of certain rigid tensor categories.  C. R. Math. Acad. Sci. Paris  {\bf 340}  (2005),  no. 8, 557--562.
     
          \smallskip   \bibitem{Pr} C. Praagman, ÒThe formal classification of linear difference operatorsÓ, Proc. Kon. Ned. Ac. Wet., ser. a  (1983).
 
  \smallskip   \bibitem{Pu} A. Pulita, $p$-adic confluence of $q$-difference equations, to appear in Compos. Math.  
  
  \smallskip   \bibitem{PRe} M. van der Put, M. Reversat, 
 Galois theory of $q$-difference equations,  to appear in Ann. Fac. Sci. Toulouse.
 
    \smallskip   \bibitem{Rai}  D. Ra\u\i kov,  Semiabelian categories, Dokl. Akad. Nauk SSSR {\bf 188} (1969), 1006-1009 (English translation in Soviet Math. Dokl. 1969. V. 10. P. 1242-1245). 
  
   \smallskip   \bibitem{RR}  S. Ramanan, A. Ramanathan, Some remarks on the instability flag, Tohoku Math. J. {\bf 36} (1984), 269-291. 
 
 \smallskip   \bibitem{R1}  L. Ramero, Local monodromy in non-Archimedean analytic geometry, Publ. Math. I.H.E.S. {\bf 102} (2005), 167-280. 
 
  \smallskip   \bibitem{R2}  L. Ramero, Local monodromy in non-Archimedean analytic geometry II, preprint ArXiv/0703225.
 
   \smallskip   \bibitem{RS}  J.-P. Ramis, J. Sauloy,
the $q$-analogue of the wild fundamental group,
RIMS K\^oky\^uroku Bessatsu {\bf B2}, 167-193 (2007). 
  
 \smallskip   \bibitem{Ra}  M. Rapoport, Analogien zwischen den Modulr\"aumen von Vektorb\"undeln und von Flaggen, DMV Tagung Ulm, 1995.
  
     \smallskip   \bibitem{Rud} A. Rudakov, Stability for an abelian category, J. of Algebra {\bf 197} (1997), 231-245.

    \smallskip   \bibitem{Rum} W. Rump, A counterexample to Raikov's conjecture, Bull. London Math. Soc.,  Advance Access Sept. 4, 2008. doi:10.1112/blms/bdn080
   
   
   \smallskip   \bibitem{Sa}  C. Sabbah, {\it Equations diff\'erentielles \`a points singuliers r\'eguliers et ph\'enom\`enes de Stokes en dimension $2$}, Ast\'erisque {\bf 263} (2000).
   
 \smallskip   \bibitem{Saa} N. Savedra Rivano, {\it Cat\'egories tannakiennes}, Lecture Notes in Mathematics {\bf 265}, Springer Verlag  (1972).

  \smallskip   \bibitem{Sau1} J. Sauloy, La filtration canonique par les pentes d'un module aux q-diff\'erences et le gradu\'e associ\'e, 
Annales de l'institut Fourier, {\bf 54} no. 1 (2004), p. 181-210

  \smallskip   \bibitem{Sau2} J. Sauloy, Equations aux $q$-diff\'erences et fibr\'es vectoriels holomorphes sur la courbe elliptique $\C^\ast/q^\Z$, preprint ArXiv:math/07114031.
    
  \smallskip   \bibitem{Sc} J.-P. Schneiders, {\it Quasi-abelian categories and sheaves}, M\'em. S. M. F. {\bf 76} (1991).
  
         \smallskip \bibitem{Se} J.-P. Serre, {\it Corps locaux}, Hermann, 1968. 
         
   \smallskip   \bibitem{Sh} S. Shatz, {The decomposition and specialization of algebraic families of vector bundles}, Compos. Math. {\bf 35} (1977), 163-187.
   
     \smallskip   \bibitem{Si}  C. Simpson,
Higgs bundles and local systems, Publ. Math. IHES {\bf 75} (1992), 5-95. 

      \smallskip   \bibitem{SC} R. Succi Cruciani, Sulle categorie quasi abeliane, Rev. Roumaine Math. Pures Appl.  {\bf 18}  (1973), 105--119. 
      
        \smallskip   \bibitem{St} U. Stuhler, Eine Bemerkung zur Reduktionstheorie quadratischen Formen, Archiv der Math. {\bf 27} (1976), 604-610.
   
   \smallskip   \bibitem{Ta} F. Takemoto, Stable vector bundles on algebraic surfaces.  Nagoya Math. J.  {\bf 47} (1972), 29--48.
   
    \smallskip   \bibitem{Tom} M. Toma, Stable bundles on non-algebraic surfaces giving rise to compact
moduli spaces, C.R. Acad. Sci. Paris, 323 (1996), 501-505.
   
      \smallskip   \bibitem{To1} B. Totaro,  Tensor products of semistables are semistable, in {\it Geometry and Analysis on Complex Manifolds}, Festschrift for Professor S. Kobayashi's 60th Birthday, ed. T. Noguchi, J. Noguchi, and T. Ochiai, World Scientific Publ. Co., Singapore, 1994, pp. 242-250.

    \smallskip   \bibitem{To2} B. Totaro, Tensor products in $p$-adic Hodge theory, Duke Math. J. {\bf 83} (1996), 79-104.
    
     \smallskip  \bibitem{Tsu} N. Tsuzuki,  The local index and the Swan conductor, Compos. Math. {\bf 111} (1998), 245-288.
     
     \smallskip  \bibitem{Tu}    H. Turrittin, Convergent solutions of ordinary differential equations in the neighborhood of an irregular point, Acta Math. {\bf 93} (1955), 27-66.
     
  \smallskip  \bibitem{X1}  L. Xiao. On ramification filtrations and p-adic differential equations, I: equal characteristic 
case, preprint  arXiv: 0801.4962. 

   \smallskip  \bibitem{X2}  L. Xiao. On ramification filtrations and p-adic differential equations, II: mixed characteristic  
case, preprint arXiv: 08113792. 

   \smallskip  \bibitem{Y}  N. Yoneda, On Ext and exact sequences, J. Fac. Sci. Univ. Tokyo {\bf 18} (1960), 507-576.
 \end{thebibliography}
 \end{document}